\documentclass[11pt]{article}
\usepackage{amsfonts}
\usepackage{color}
\newtheorem{theo}{Theorem}[section]
\newtheorem{lem}[theo]{Lemma}
\newtheorem{cor}[theo]{Corollary}

\newtheorem{defi}[theo]{Definition}

\newcommand{\mysection}[1]{\section{#1} \setcounter{equation}{0}}
\newcommand{\proof}{{\sc Proof.} \quad}
\newcommand{\proofc}{{\sc Proof} \ }
\newcommand{\be}{\begin{equation} \label}
\newcommand{\ee}{\end{equation}}
\newcommand{\bea}{\begin{eqnarray}\label}
\newcommand{\eea}{\end{eqnarray}}
\newcommand{\bas}{\begin{eqnarray*}}
\newcommand{\eas}{\end{eqnarray*}}
\newcommand{\bit}{\begin{itemize}}
\newcommand{\eit}{\end{itemize}}
\newcommand{\qed}{\hfill$\Box$ \vskip.2cm}
\newcommand{\nn}{\nonumber}
\newcommand{\R}{\mathbb{R}}
\newcommand{\N}{\mathbb{N}}
\newcommand{\pO}{\partial\Omega}

\newcommand{\eps}{\varepsilon}

\newcommand{\supp}{{\rm supp} \, }

\newcommand{\wsto}{\stackrel{\star}{\rightharpoonup}}
\newcommand{\hra}{\hookrightarrow}
\newcommand{\io}{\int_\Omega}
\newcommand{\mint}{- \hspace*{-3.3mm} \int}
\newcommand{\Mint}{- \hspace*{-4mm} \int}
\newcommand{\mult}{\otimes}

\newcommand{\abs}{\\[5pt]}

\newcommand{\proj}{{\cal P}}
\newcommand{\neps}{n_\eps}
\newcommand{\ceps}{c_\eps}
\newcommand{\ueps}{u_\eps}

\newcommand{\Peps}{P_\eps}
\newcommand{\yeps}{Y_\eps}
\newcommand{\tme}{T_{max,\eps}}
\newcommand{\tz}{T_0}
\newcommand{\tu}{\widetilde{u}}
\newcommand{\hu}{\widehat{u}}
\newcommand{\onz}{\overline{n_0}}
\newcommand{\one}{{\bf{1}}}
\newcommand{\set}{{\cal S}_{m,M,L,\tz}}
\newcommand{\F}{{\cal F}_\kappa}

\newcommand{\xitz}{\xi_{t_0}}
\begin{document}
\enlargethispage{10mm}
\title{How far do chemotaxis-driven forces influence regularity\\
in the Navier-Stokes system?}
\author{
Michael Winkler\footnote{michael.winkler@math.uni-paderborn.de}\\
{\small Institut f\"ur Mathematik, Universit\"at Paderborn,}\\
{\small 33098 Paderborn, Germany} }
\date{}
\maketitle
\begin{abstract}
\noindent 
  The chemotaxis-Navier-Stokes system
    \begin{equation}\label{00}
    \left\{ \begin{array}{rcll}
    n_t + u\cdot\nabla n &=& \Delta n - \nabla \cdot (n\chi(c)\nabla c),\\[1mm]
    c_t + u\cdot\nabla c &=& \Delta c-nf(c),    \\[1mm]
    u_t + (u\cdot\nabla)u &=& \Delta u + \nabla P + n \nabla \Phi,  \\[1mm]
    \nabla \cdot u &=& 0,
    \end{array} \right.
	\qquad \qquad (\star)
  \end{equation}
  is considered 
  under boundary conditions of homogeneous Neumann type for $n$ and $c$, and Dirichlet type for $u$,
  in a bounded convex domain $\Omega\subset \R^3$ with smooth boundary, 
  where $\Phi \in W^{1,\infty}(\Omega)$ and $\chi$ and $f$ 
  are sufficiently smooth given functions generalizing the prototypes $\chi\equiv const.$ and $f(s)=s$ for $s\ge 0$.\abs
  It is known that for all suitably regular initial data $n_0, c_0$ and $u_0$ satisfying $0\not\equiv n_0\ge 0$, $c_0\ge 0$
  and $\nabla\cdot u_0=0$, a corresponding initial-boundary
  value problem admits at least one global weak solution which can be obtained as the pointwise limit
  of a sequence of solutions to appropriately regularized problems.
  The present paper shows that after some relaxation time, this solution enjoys further regularity properties
  and thereby complies with the concept of {\em eventual energy solutions} which is newly introduced here,
  and which inter alia requires that two quasi-dissipative inequalities are ultimately satisfied.\abs
  Moreover, it is shown that actually for any such eventual energy solution $(n,c,u)$
  there exists a waiting time $T_0\in (0,\infty)$
  with the property that $(n,c,u)$ is smooth in $\bar\Omega\times [T_0,\infty)$, and that 
  \bas
	n(x,t)\to \overline{n_0},
	\qquad
	c(x,t)\to 0
	\qquad \mbox{and} \qquad
	u(x,t)\to 0
  \eas
  hold as $t\to\infty$, uniformly with respect to $x\in\Omega$.\abs
  This resembles a classical result on the three-dimensional Navier-Stokes system, asserting eventual smoothness
  of arbitrary weak solutions thereof which additionally fulfill the associated natural energy inequality.
  In consequence, our results inter alia indicate that 
  under the considered boundary conditions, the possibly destabilizing action of chemotactic cross-diffusion
  in ($\star$) does not substantially affect the regularity properties of the fluid flow at least on large time scales.\abs
\end{abstract}
\newpage
\section{Introduction}\label{intro}
{\bf Chemotaxis and blow-up.} \quad
When primitive microorganisms interact with their environment, their individually unstructured behavior
may switch to quite complicated dynamics at macroscopic levels.
Prototypical situations include spontaneous formation of aggregates such as in slime mold formation
processes e.g.~of {\em Dictyostelium discoideum}, organization of cell positioning during embryonic development,
or also the invasion of tumors into healthy tissue. 
An important role in numerous structure-generating processes is known to be played by various types of 
biased cell movement in response to external cues such as chemical signal substances, mechanical stimuli, or gradients in 
voltage or acidic concentration, for instance. 
Such {\em taxis} mechanisms have been thoroughly studied in various contexts, also at a theoretical level, 
with the celebrated Keller-Segel system of chemotaxis consituting the apparently most paradigmatic
representative in the field of macroscopic mathematical models (\cite{keller_segel}; see also
\cite{hillen_painter2009} for a comprehensive survey on modeling aspects).
Indeed, intense analysis on the latter has confirmed the conjecture that spontaneous formation
of aggregates, in the extreme mathematical sense of finite-time blow-up of solutions, may arise even in the simple 
two-component framework containing a population of cells moving chemotactically upward gradients of a signal substance,
provided that the system reinforces itself in that cells produce the chemical in question
(\cite{herrero_velazquez}, \cite{jaeger_luckhaus}, \cite{win_JMPA}, \cite{mizoguchi_win}).\abs
{\bf Chemotaxis-fluid interaction.} \quad
In the case of even more primitive organisms, chemotactically moving toward a nutrient which they consume rather
than produce, 
the correspondingly modified chemotaxis system 
possesses global bounded smooth solutions in the spatially two-dimensional setting, whereas
in the three-dimensional counterpart at least global weak solutions can be constructed which eventually become smooth
and bounded (\cite{taowin_JDE2}).
On the other hand, more recent findings
indicate that also populations of such simple	
individuals may exhibit quite colorful collective behavior:
As suggested by striking experiments revealing spontaneous formation of plume-like aggregates in populations
of {\em Bacillus subtilis} suspended in sessile water drops,
in such situations it may be necessary to take into account
the mutual interaction of cells and their movement on the one hand, 
and of the surrounding medium on the other.
Accordingly, the modeling approach in \cite{goldstein2004} and \cite{goldstein2005}
in particular assumes that besides chemotactic movement, signal consumption 
and transport of both cells and signal through the fluid, 
there is a significant buoyancy-driven effect of cells on the fluid dynamics.
One is thereby led to considering the coupled chemotaxis-Navier-Stokes system
\be{0}
    \left\{ \begin{array}{rcll}
    n_t + u\cdot\nabla n &=& \Delta n - \nabla \cdot (n\chi(c)\nabla c),
	\qquad & x\in\Omega, \ t>0,\\[1mm]
    c_t + u\cdot\nabla c &=& \Delta c-nf(c), \qquad & x\in\Omega, \ t>0,   \\[1mm]
     u_t + (u\cdot\nabla )u  &=& \Delta u + \nabla P + n \nabla \Phi, \qquad & x\in\Omega, \ t>0, \\[1mm]
    \nabla \cdot u &=& 0, \qquad & x\in\Omega, \ t>0,
    \end{array} \right.
\ee
for the unknown $(n,c,u,P)$ in the physical domain $\Omega \subset \R^N$, 
where the chemotactic sensitivity $\chi$, the signal consumption rate $f$ and the gravitational potential $\Phi$ 
are given parameter functions.\abs
From a 		
viewpoint of mathematical analysis, 
this system couples the well-known obstacles from the theory of the
Navier-Stokes equations to the typical difficulties arising in the study of chemotaxis systems.
Up to now, in the case $N=3$ it is not only unknown whether the incompressible Navier-Stokes equations possess
global smooth solutions for arbitrarily large smooth initial data (e.g.~under Dirichlet boundary conditions in bounded
domains, cf.~\cite{leray} or also \cite{wiegner} and \cite{sohr}); 
also the chemotaxis-only subsystem of (\ref{0}) obtained on letting $u\equiv 0$
is far from understood in this framework, with no answer available yet e.g.~to the question whether 
the global weak
solutions, known to exist for any reasonably regular initial data in a corresponding Neumann-type initial-boundary
value problem in bounded convex domains, may blow up in finite time before becoming 
ultimately smooth (\cite{taowin_JDE2}).\abs
Accordingly, the knowledge on such coupled chemotaxis-fluid systems is at a rather early stage yet, with most
previous works focusing on the basic issues of global solvability in various functional frameworks. 
For instance, global existence of uniquely determined smooth solutions is known for an initial-boundary value problem 
associated with the two-dimensional version of (\ref{0}), under structural
assumptions on the parameter functions $\chi$ and $f$ which are mild enough so as to include
the prototypical choices
\be{proto}
	\chi\equiv const.
	\qquad \mbox{and} \qquad
	f(s)=s, \quad s\ge 0,
\ee
and for all reasonably smooth initial data (\cite{win_CPDE}).\abs
In the context of three-dimensional frameworks, large bodies of the existing literature address variants of (\ref{0})
involving diverse types of regularizing modifications.
For the 
chemotaxis-Stokes system obtained from (\ref{0})
on neglecting the convective term $(u\cdot\nabla)u$ in the fluid evolution, 
global existence results 	
have been derived for the Cauchy problem in $\R^3$ under certain additional requirements on $\chi$ and $f$
and a smallness assumption e.g.~on $c$ (\cite{DLM}), and also for a corresponding boundary-value problem
without any such further restrictions (\cite{win_CPDE}).
Even in this simplified setting, all these solutions constructed so far 
are merely weak solutions, with widely unknown boundedness and 
regularity properties which in fact might be poor so as to be consistent with several conceivable types of blow-up 
phenomena, for example of finite-time blow up of $n$ with respect to the norm in $L^\infty(\Omega)$ 
(cf.~also \cite{chae_kang_lee_CPDE} for a detailed discussion of extensibility criteria for local-in-time
smooth solutions).\abs
As an additional regularizing mechanism, numerous works study the enhancement of cell diffusion at large densities,
modeled by replacing $\Delta n$ with the porous medium-type diffusion term $\Delta n^m$ for $m>1$.
Under 	
mild assumptions on $\chi$ and $f$, 
the corresponding three-dimensional chemotaxis-Stokes system then again admits global weak solutions
whenever $m>1$ 	
(\cite{duan_xiang_IMRN2012}; cf.~also \cite{liu_lorz} for a precedent);
in the case when moreover $m>\frac{8}{7}$ and $\Omega$ is a bounded convex domain in $\R^3$, the first component
$n$ of such a solution is in fact locally bounded in $\bar\Omega\times [0,\infty)$ (\cite{taowin_ANNIHP}), and
if even $m>\frac{7}{6}$, then in the latter situation $n$ actually remains bounded in all of $\Omega\times (0,\infty)$
(\cite{win_ctf_3d_nonlinear_general}).
In the case $m>\frac{4}{3}$, global existence of -- possibly unbounded -- weak solutions has been established 
in \cite{vorotnikov} even for the associated full chemotaxis-Navier-Stokes system;
results on global existence and boundedness in two-dimensional chemotaxis-fluid systems 
with nonlinear cell diffusion can be found
in \cite{DiFLM}, \cite{taowin_DCDSA} and \cite{ishida_2d}.
Examples of further regularizations, as discussed in the literature with regard to global weak solvability,
consist in considering saturation effects in the
cross-diffusive term at large cell densities (\cite{cao_wang}), or also including logistic-type cell proliferation
and death (\cite{vorotnikov}).\abs
{\bf Solvability and asymptotics in the three-dimensional chemotaxis-Navier-Stokes system.} \quad
Concerning the full three-dimensional chemotaxis-Navier-Stokes system (\ref{0}) with linear cell diffusion, the question
of global solvability is apparently more delicate, and
accordingly the first result in this direction resorted to the construction of global solutions to a corresponding
Cauchy problem in $\R^3$ which emenate from initial data suitably close to one of the constant equilibria
$(a,0,0)$ for $a>0$ (\cite{DLM}).
As for general -- and, in particular, large -- initial data, in view of the limited knowledge on global regularity in the
Navier-Stokes subsystem of (\ref{0}) only weak solutions can currently be expected.
A result optimal in this respect has recently been achieved in \cite{win_ct_nasto_exist}, where it has been shown that
under the assumptions (\ref{reg_coeff}) and (\ref{struct}) below,
for any initial data fulfilling (\ref{init}) the problem (\ref{0}) possesses at least one global weak solution.\\
Concerning the large time behavior of solutions, only very few seems known even in simplified situations: 
Except for some results on convergence of solutions satisfying certain smallness conditions
(\cite{DLM}, \cite{chae_kang_lee_CPDE}), the only statememts we are aware of which cover arbitrarily large initial data
in (\ref{0}) address its two-dimensional version in bounded convex domains 
in which any solution approaches the constant state $(\onz,0,0)$ in the large time limit at an exponential rate, 
where $\onz:=\mint_\Omega n_0>0$ (\cite{win_ARMA}, \cite{zhang_li}).
In the three-dimensional counterpart, a similar stabilization result so far could be obtained only for the 
chemotaxis-Stokes variant in which additional regularity is enforced by 
the presence of porous medium-type cell diffusion (\cite{win_ctf_3d_nonlinear_general}).\abs
{\bf Main results: Eventual smoothness and stabilization.} \quad
The objective of the present work is to undertake a further step toward a qualitative understanding of the
chemotaxis-fluid interaction modeled in (\ref{0}),
especially with regard to possible effects of the chemotaxis-driven forcing on the fluid motion, 
and of the latter on the distribution of cells.
Our main results in this direction will reveal that any such mutual influence		
will in fact disappear asymptotically
in that the large time behavior of solutions is essentially governed by the
decoupled chemotaxis-only and Navier-Stokes subsystems obtained on neglecting the components $u$ and $(n,c)$, respectively.
In particular, as in the unforced Navier-Stokes equations (\cite{wiegner}), the fluid velocity $u$ will become 
smooth eventually and decay uniformly in the large time limit;
likewise, the couple $(n,c)$ enjoys enjoys a similar ultimate smoothness property and approaches the spatially homogeneous
limit $(\onz,0)$ associated with the respective mass level, thus resembling the behavior in the associated
fluid-free chemotaxis system (\cite{taowin_JDE2}).\abs
In order to make this more precise, let us consider (\ref{0}) in a bounded convex domain $\Omega\subset \R^3$ 
with smooth boundary, along with the initial conditions
\be{0i}
    n(x,0)=n_0(x), \quad c(x,0)=c_0(x) \quad \mbox{and} \quad u(x,0)=u_0(x), \qquad x\in\Omega,
\ee
and under the boundary conditions
\be{0b}
	\frac{\partial n}{\partial\nu}=\frac{\partial c}{\partial\nu}=0 
	\quad \mbox{and} \quad u=0
	\qquad \mbox{on } \pO.
\ee
Here we shall require that
\be{init}
    \left\{
    \begin{array}{l}
    n_0 \in L \log L(\Omega) \quad \mbox{is nonnegative with $n_0\not\equiv 0$,} \quad
	\mbox{that} \\
    c_0 \in L^\infty(\Omega) \quad \mbox{ is nonnegative and such that $\sqrt{c_0} \in W^{1,2}(\Omega)$, \quad and that}\\
    u_0 \in L^2_\sigma(\Omega),
    \end{array}
    \right.
\ee
where as usual, $L\log L(\Omega)$ denotes the Orlicz space corresponding to the Young function $[0,\infty) \ni z \mapsto
z\ln (1+z)$, and where for $p>1$, by 
$L^p_\sigma(\Omega):=\{ \varphi\in L^p(\Omega;\R^3) \ | \ \nabla \cdot \varphi=0 \}$ 
we abbreviate the space of all solenoidal vector fields in $L^p(\Omega)$.\\
Throughout this paper, the chemotactic sensitivity $\chi$, the signal consumption rate $f$ in (\ref{0}) and the gravitational
potential $\Phi$ are assumed to be such that
\be{reg_coeff}
	\left\{ \begin{array}{l}
	\chi \in C^2([0,\infty)) \quad \mbox{is positive on $[0,\infty)$},\\
	f\in C^1([0,\infty)) \quad \mbox{is positive on $(0,\infty)$ with $f(0)=0$, that} \\
	\Phi\in W^{1,\infty}(\Omega),
	\end{array} \right.
\ee
and that the structural requirements
\be{struct}
	\Big(\frac{f}{\chi}\Big)'>0,
	\quad 
	\Big(\frac{f}{\chi}\Big)'' \le 0
	\quad \mbox{and} \quad
	(\chi\cdot f)'' \ge 0
	\qquad \mbox{on } [0,\infty)
\ee
are fulfilled, noting that the latter hypotheses are mild enough so as to allow e.g.~for the choices in (\ref{proto}).\abs
Within this framework, in view of the unsolved uniqueness problem for the Navier-Stokes equations we cannot
expect weak solutions of (\ref{0}) to be unique; accordingly, 
it seems desirable to derive 
results on qualitative behavior which are independent of a particular construction of solutions.
Inspired by analogues from the analysis of the Navier-Stokes system, we shall thus consider rather arbitrary weak solutions
enjoying certain additional properties which are essentially linked to strutural features of (\ref{0}).
Here besides the natural energy inequality (\ref{energy1}) 
associated with the Navier-Stokes system, singling out the so-called
turbulent solutions among all weak solutions of the latter, a central role will be played by a second energy-like
functional $\F$ which for given $\kappa>0$ is defined by
\be{def_energy}
	\F[n,c,u] := \io n\ln n 
	+ \frac{1}{2} \io \frac{\chi(c)}{f(c)} |\nabla c|^2 
	+ \kappa \io |u|^2
\ee
whenever 
$n\in L\log L(\Omega)$ and $c\in W^{1,2}(\Omega)$ are nonnegative and such that 
$\frac{\chi(c)}{f(c)}|\nabla c|^2 \in L^1(\Omega)$, and $u\in L^2(\Omega;\R^3)$ (\cite{win_ct_nasto_exist}).\abs
We now select a subclass of weak solutions to (\ref{0}) as follows.
\begin{defi}\label{defi_ees}
  Suppose that $(n,c,u)$ is a global weak solution of (\ref{0}) in the sense of Definition \ref{defi_weak} below. 
  Then we call $(n,c,u)$ an {\em eventual energy solution} of (\ref{0}) if there exists $T>0$ such that
  \be{reg_ees}
	\begin{array}{l}
	n\in L^4_{loc}(\bar\Omega\times [T,\infty)) \cap L^2_{loc}([T,\infty);W^{1,2}(\Omega))
	\quad \mbox{with} \quad 
	n^\frac{1}{2} \in L^2_{loc}([T,\infty);W^{1,2}(\Omega)), \\
	c\in L^\infty_{loc}(\Omega\times [T,\infty))
	\quad \mbox{with} \quad
	c^\frac{1}{4} \in L^4_{loc}([T,\infty);W^{1,4}(\Omega))
	\quad \mbox{and} \\
	u\in L^\infty_{loc}([T,\infty);L^2_\sigma(\Omega))
	\cap L^2_{loc}([T,\infty);W_0^{1,2}(\Omega)),
	\end{array}
  \ee
  if 
  \be{energy1}
	\frac{1}{2} \io |u(\cdot,t)|^2 + \int_{t_0}^t \io |\nabla u|^2
	\le \frac{1}{2} \io |u(\cdot,t_0)|^2 
	+\int_{t_0}^t \io nu\cdot\nabla \Phi
	\qquad \mbox{for a.e.~$t_0>T$ and all } t>t_0,
  \ee
  and if there exist $\kappa>0$ and $K>0$ such that
  \bea{energy}
	\frac{d}{dt} \F[n,c,u](t) 
	+ \frac{1}{K} \io \bigg\{ \frac{|\nabla n|^2}{n} + \frac{|\nabla c|^4}{c^3} + |\nabla u|^2 \bigg\}
	\le K
	\qquad \mbox{in } {\mathcal D}'((T,\infty)).
  \eea
\end{defi}
{\bf Remark.} \quad
  The regularity assumptions in (\ref{reg_ees}) warrant that $\frac{|\nabla n|^2}{n}$, $\frac{|\nabla c|^4}{c^3}$ and
  $|\nabla u|^2$ belong to $L^1_{loc}(\bar\Omega\times (T,\infty))$,
  and that moreover $(T,\infty)\ni t\mapsto \F[n,c,u](t) \in L^1_{loc}((T,\infty))$ (cf.~Lemma \ref{lem333}),
  implying that (\ref{energy}) indeed is meaningful.\abs
It has been shown in \cite{win_ct_nasto_exist} that under the assumptions (\ref{reg_coeff}) and (\ref{struct}),
for any initial data fulfilling (\ref{init}) the problem (\ref{0}) possesses at least one global weak solution
in the natural sense specified in Definition \ref{defi_weak} below. 
The first of our results asserts that this solution actually enjoys all the above properties of
an eventual energy solution; we shall thereby prove the following.
\begin{theo}\label{theo_exist_ees}
  Let (\ref{reg_coeff}) and (\ref{struct}) hold, and assume that
  $n_0, c_0$ and $u_0$ comply with (\ref{init}). 
  Then there exists at least one eventual energy solution of (\ref{0}).
\end{theo}
We can secondly prove that in fact any such eventual energy solution becomes smooth ultimately, and that it approaches
the unique spatially homogeneous steady state compatible with the preserved total mass $\io n_0>0$.
\begin{theo}\label{theo_eventual}
  Let (\ref{reg_coeff}) and (\ref{struct}) hold, 
  and suppose that $(n,c,u)$ is an eventual energy solution of (\ref{0}) with some initial data 
  $n_0, c_0$ and $u_0$ satisfying (\ref{init}).
  Then there exist $T>0$ and $P\in C^{1,0}(\bar\Omega\times [T,\infty))$ such that
  \be{reg}
	\left\{ \begin{array}{l}
	n\in C^{2,1}(\bar\Omega\times [T,\infty)), \nn\\
	c\in C^{2,1}(\bar\Omega\times [T,\infty)) \qquad \mbox{and} \nn\\
	u\in C^{2,1}(\bar\Omega\times [T,\infty);\R^3),
	\end{array} \right.
  \ee
  and such that $(n,c,u,P)$ solves the boundary value problem in (\ref{0}) classically in $\bar\Omega\times [T,\infty)$.
  Furthermore,
  \bea{conv}
	n(\cdot,t)\to \onz
	\quad \mbox{in } L^\infty(\Omega),
	\qquad
	c(\cdot,t)\to 0
	\quad \mbox{in } L^\infty(\Omega)
	\qquad \mbox{and} \qquad
	u(\cdot,t)\to 0
	\quad \mbox{in } L^\infty(\Omega)
  \eea
  as $t\to\infty$, where $\onz:=\mint_\Omega n_0$.
\end{theo}
Theorem \ref{theo_eventual} may be interpreted as rigorously reflecting that despite the possibly disordering influence of the
fluid, the signal consumption process in (\ref{0}) occurs in such a regular manner that ultimately even the gradients
of the chemical become irrelevant with regard to their chemoattractive impact,
and that in consequence the cell population homogenizes efficiently enough so as to let any substantial destabilizing effect
on the fluid vanish asymptotically.
This may become substantially different in situations when 
different types of boundary conditions are considered, possibly accounting for signal influx,
or when signal absorption is replaced with mechanisms of
signal production by cells, as recently proposed and studied in contexts involving fluid interaction 
in \cite{kiselev_ryzhik1}, \cite{kiselev_ryzhik2}, \cite{espejo_suzuki} and \cite{tao_ctf}. \abs
{\bf Main ideas. Organization of the paper.} \quad
The overall strategy pursued in the course of our reasoning
consists in showing that firstly the solution component $c$ must decay with
respect to the norm in $L^\infty(\Omega)$ as $t\to\infty$, and that secondly appropriate smallness of this component in
$L^\infty(\Omega\times (\tz,\infty))$ for some $\tz>0$ entails smoothness of $(n,c,u)$ in $\bar\Omega \times (T_1,\infty)$
for some $T_1>\tz$.\abs
The first of these properties will be a consequence of some basic dissipative features of (\ref{0}) combined with
suitable uniform-in-time regularity estimates implied by the energy inequality (\ref{energy}) (Sections \ref{sect3} and 
\ref{sect4}).
In accomplishing the second of the mentioned steps,
we will generalize a related procedure pursued in \cite{win_ARMA} for smooth solutions of the two-dimensional version
of (\ref{0}), where a similar conclusion was derived on the basis of the observation that for any given $p\ge 2$
the functional 
$\io \frac{n^p}{\delta-c}$ acts as an entropy, provided that $\delta=\delta(p)$ is suitably 
small and $c$ remains below the threshold $\delta$ throughout evolution.
Since in the present case we intend to address arbitrary eventual energy solutions,
the lack of a priori knowledge on regularity properties beyond those listed in Definition \ref{defi_ees}
will require the use of a substantially more subtle testing technique to track the time evolution of functionals
of the above type, simultaneously involving the first two equations in (\ref{0}).
Moreover, this limited information on regularity will force us to firstly restrict our key statement in this direction, 
presented in the extensive Lemma \ref{lem50}, to functionals of the form $\io \psi(n)\rho(c)$ with convex
$\psi$ and $\rho$ subject to technical assumptions which inter alia require that $\psi(n)$ does not increase 
faster than $n^\frac{9}{5}$ as $n\to\infty$; only in a second step we will see in Section \ref{sect6}
by means of an approximation argument that actually any algebraic growth of $\psi$ can be achieved.\abs
In Section \ref{sect7} we shall infer from the correspondingly gained entropy-dissipation inequalities that
$n$ stabilizes in a sense yet weaker than claimed in Theorem \ref{theo_eventual}, but strong enough to assert
a decay property of the forcing term in the Navier-Stokes system in (\ref{0}) which is sufficient to imply 
decay of $u$ with respect to the norm in $L^p(\Omega)$ for any finite $p\ge 1$ (Lemma \ref{lem57}).
The eventual integrability property of $u$ thereby implied will enable us to perform 
a series of arguments based on maximal Sobolev regularity in the Stokes evolution system and 
inhomogeneous linear heat equations to successively obtain further ultimate regularity properties of $u$, $c$ and $n$
which by standard Schauder theory imply eventual smoothness
(Lemma \ref{lem58}--Lemma \ref{lem556}).
This improved knowledge on regularity 
thereupon allows for turning the weak decay information previously gathered into the desired uniform
convergence statements and thereby complete the proof of Theorem \ref{theo_eventual} in Section \ref{sect8}.\abs
Finally, in proving Theorem \ref{theo_exist_ees} in Section \ref{sect9}
we shall make use of the fact that our arguments in Section \ref{sect3} through Section \ref{sect6} are formulated
in a manner slightly more general than used in the mere analysis of eventual energy solutions, namely simultaneously
covering also all of the solutions to the approximate systems (\ref{0eps}), uniformly with respect to the
regularization parameter $\eps\in (0,1)$. Since an appropriate sequence of such solutions is known to approach
a weak solution of (\ref{0}) satisfying (\ref{energy1}) and (\ref{energy}) for $\tz:=0$, 
the additional properties thereby obtained assert that this limit in fact is an eventual energy solution.\abs
Throughout the paper, we let $A:=-\proj \Delta$ denote the Stokes operator which for any $p\in (1,\infty)$
is sectorial in $L^p_\sigma(\Omega)$ when considered with domain 
$D(A)\equiv D(A_p)=W^{2,p}(\Omega) \cap W_0^{1,p}(\Omega)\cap L^p_\sigma(\Omega)$, 
and hence generates the analytic Stokes semigroup $(e^{-tA})_{t\ge 0}$.
Here, by $\proj$ we mean the associated Helmholtz projection mapping $L^p(\Omega)$ onto $L^p_\sigma(\Omega)$
(cf.~\cite{giga1986}, \cite{giga_sohr}, \cite{giga1981_the_other}).
\mysection{Weak solutions}
The following notion of weak solutions to (\ref{0}) is taken from \cite{win_ct_nasto_exist}.
Here and in the sequel, for vectors $v\in \R^3$ and $w\in\R^3$ we let $v\mult w$ denote the matrix 
$(a_{ij})_{i,j\in \{1,2,3\}}\in\R^{3\times 3}$ defined on setting $a_{ij}:=v_i w_j$ for $i,j\in \{1,2,3\}$.
\begin{defi}\label{defi_weak}
  By a {\em global weak solution} of (\ref{0}), (\ref{0i}), (\ref{0b}) we mean a triple $(n,c,u)$ of functions
  \bea{reg_weak}
	n\in L^1_{loc}([0,\infty);W^{1,1}(\Omega)), 
	\quad
	c\in L^1_{loc}([0,\infty);W^{1,1}(\Omega)), 
	\quad
	u\in L^1_{loc}([0,\infty);W_0^{1,1}(\Omega;\R^3)), 
  \eea
  such that $n\ge 0$ and $c\ge 0$ a.e.~in $\Omega\times (0,\infty)$,
  \bea{reg_weak2}
	& & nf(c)\in L^1_{loc}(\bar\Omega\times [0,\infty)),
	\qquad 
	u\mult u \in L^1_{loc}(\bar\Omega\times [0,\infty);\R^{3\times 3}),
	\qquad \mbox{and} \nn\\
	& & n\chi(c)\nabla c, nu \mbox{ and } cu
	\ \mbox{belong to } L^1_{loc}(\bar\Omega\times [0,\infty);\R^3),
  \eea
  that $\nabla \cdot u=0$ a.e.~in $\Omega\times (0,\infty)$, and that
  \bea{w1}
	-\int_0^\infty \io n\phi_t - \io n_0\phi(\cdot,0)
	= - \int_0^\infty \io \nabla n\cdot \nabla \phi
	+ \int_0^\infty \io n\chi(c) \nabla c \cdot\nabla\phi
	+ \int_0^\infty \io n u \cdot \nabla\phi
  \eea
  for all $\phi\in C_0^\infty(\bar\Omega\times [0,\infty))$,
  \bea{w2}
	-\int_0^\infty \io c\phi_t - \io c_0\phi(\cdot,0)
	= - \int_0^\infty \io \nabla c\cdot \nabla \phi
	- \int_0^\infty \io nf(c) \phi
	+ \int_0^\infty \io c u\cdot\nabla \phi
  \eea
  for all $\phi\in C_0^\infty(\bar\Omega\times [0,\infty))$ as well as
  \bea{w3}
	-\int_0^\infty \io u\cdot\phi_t - \io u_0\cdot \phi(\cdot,0)
	= - \int_0^\infty \io \nabla u \cdot\nabla \phi
	+ \int_0^\infty u\mult u \cdot \nabla \phi
	+ \int_0^\infty \io n\nabla \Phi\cdot \phi
  \eea
  for all $\phi\in C_0^\infty(\Omega\times [0,\infty);\R^3)$ satisfying $\nabla\cdot \phi\equiv 0$.
\end{defi}
The above solution concept meets the basic natural requirement that solutions preserve mass during evolution.
\begin{lem}\label{lem_mass}
  Suppose that $(n,c,u)$ is a global weak solution of (\ref{0}). Then $n\in L^\infty((0,\infty);L^1(\Omega))$ with
  \be{mass_ees}
	\io n(\cdot,t) = \io n_0
	\qquad \mbox{for a.e.~} t>0.
  \ee
\end{lem}
\proof
  Let $t_0>0$ be a Lebesgue point of $(0,\infty)\ni t \mapsto \io n(x,t)dx$.
  For $\delta\in (0,1)$ we then approximate
  \be{zeta}
	\zeta_\delta(t):=\left\{ \begin{array}{ll}
	1 \qquad & \mbox{if } t \le t_0, \\[1mm]
	\frac{t_0+\delta-t}{\delta}, \qquad & t\in (t_0,t_0+\delta), \\[1mm]
	0 \qquad & \mbox{if } t\ge t_0+\delta,
	\end{array} \right.
  \ee
  by taking any sequence $(\zeta_{\delta j})_{j\in\N} \subset
  C^\infty([0,\infty))$ fulfilling $\zeta_{\delta j} \equiv 1$ in $[0,t_0)$,
  $\zeta_{\delta j} \equiv 0$ in $(t_0+1,\infty)$ and
  $\zeta_{\delta j} \wsto \zeta_\delta$ in $W^{1,\infty}((0,t_0+1))$ as $j\to\infty$. 
  For each $\delta\in (0,1)$ and $j\in\N$, we may then use 
  $\phi(x,t):=\zeta_{\delta j}(t)$, $(x,t)\in\bar\Omega\times [0,\infty)$, as a test function in (\ref{w1}).
  In the correspondingly obtained identity
  \bas
	- \int_{t_0}^{t_0+1} \io \zeta_{\delta j}'(t) n(x,t) dxdt = \io n_0(x) dx,
  \eas
  we first let $j\to\infty$ to obtain
  \bas
	\frac{1}{\delta} \int_{t_0}^{t_0+\delta} \io n(x,t)dxdt = \io n_0(x) dx,
  \eas
  for all $\delta\in (0,1)$, whereupon we take
  $\delta\searrow 0$ to infer from the assumed Lebesgue point property of $t_0$ 
  that $\io n(\cdot,t_0)=\io n_0$.	
  Since the complement in $(0,\infty)$ of the set of all such $t_0$ has measure zero, this proves 
  (\ref{mass_ees}).
\qed
\mysection{A family of chemotaxis problems with prescribed convection}\label{sect3}
From \cite{win_ct_nasto_exist} 
we already know that (\ref{0}) possesses a global weak solution, and that this solution can be obtained
as the limit of smooth solutions
to certain regularized problems (cf.~(\ref{0eps}) and Lemma \ref{lem_limit} below).
Verifying Theorem \ref{theo_exist_ees} thus amounts to showing that these approximate solutions in fact enjoy further
regularity features which ensure that the limit in fact will be an eventual energy solution. 
In view of the circumstance that also Theorem \ref{theo_eventual} requires proving regularity,
to avoid repetitions we find it convenient to organize our line of arguments in such a way that 
in the first part of our analysis we consider a generalized variant of the first two equations in (\ref{0})
which includes both the original version appearing in (\ref{0}) and also the regularized subsystem thereof
considered later in Section \ref{sect9}.\abs
Moreover, in order to underline that several asymptotic solution properties, including decay of the component $c$, are
widely independent of the particular structure of the fluid flow, 
let us in this and the following sections consider the boundary value problem
\be{d3.1}
    	\left\{ \begin{array}{rcll}
    	n_t + \tu\cdot\nabla n &=& \Delta n - \nabla \cdot (nF'(n)\chi(c)\nabla c),
	\qquad & x\in\Omega, \ t>\tz,\\[1mm]
    	c_t + \tu\cdot\nabla c &=& \Delta c-F(n)f(c), \qquad & x\in\Omega, \ t>\tz,   \\[1mm]
	& & \hspace*{-28mm}
	\frac{\partial n}{\partial\nu}=\frac{\partial c}{\partial\nu}=0,  
	& x\in\pO, \ t>\tz,
	\end{array} \right.
\ee
where
$\tu\in L^\infty_{loc}((\tz,\infty);L^2_\sigma(\Omega)) \cap L^2_{loc}((\tz,\infty);W_0^{1,2}(\Omega))$
is a given function, and where
\be{F1}
	F\in C^1([0,\infty))
	\mbox{ is such that } 
	F(0)=0
	\qquad \mbox{and} \qquad
	0 \le F'(s)\le 1
	\quad \mbox{for all } s\ge 0
\ee
as well as
\be{F2}
	F(s) \ge \frac{s}{2}
	\qquad \mbox{for all } s \in [0,1].
\ee
For proving Theorem \ref{theo_eventual} it would be sufficient to concentrate throughout
on the case $F(s):=s, \ s\ge 0$; in the proof of Theorem \ref{theo_exist_ees}, however, we will apply
some of the results obtained for (\ref{d3.1}) upon choosing $F(s):=F_\eps(s):=\frac{1}{\eps} \ln (1+\eps s)$
for $s\ge 0$ and $\eps\in (0,1)$, which is as well consistent with (\ref{F1}) and (\ref{F2}) 
(see Section \ref{sect9}).\abs
We shall study (\ref{d3.1}) in the framework of solutions fulfilling the regularity requirements in 
Definition \ref{defi_ees}:
\begin{defi}\label{defi3}
  Let $\tz\ge 0$ and
  $\tu\in L^\infty_{loc}((\tz,\infty);L^2_\sigma(\Omega)) \cap L^2_{loc}((\tz,\infty);W_0^{1,2}(\Omega))$,
  and suppose that $F$ satisfies (\ref{F1}).
  Then a couple $(n,c)$ of nonnegative functions defined a.e.~in $\Omega\times (\tz,\infty)$ will
  be called a {\em strong solution} of the boundary value problem (\ref{d3.1}) in $\Omega\times (\tz,\infty)$ if
  \be{d3.2}
	\begin{array}{l}
	n \in L^\infty_{loc}((\tz,\infty);L^1(\Omega))
	\cap L^4_{loc}(\bar\Omega\times (\tz,\infty)) 
	\cap L^2_{loc}((\tz,\infty);W^{1,2}(\Omega))
	\quad \mbox{and} \\
	c \in L^\infty_{loc}(\bar\Omega \times (\tz,\infty))
	\cap L^4_{loc}((\tz,\infty);W^{1,4}(\Omega)),
	\end{array}
  \ee
  and if 
  \be{d3.3}
	-\int_{\tz}^\infty \io n\phi_t 		
	= - \int_{\tz}^\infty \io \nabla n \cdot \nabla \phi
	+ \int_{\tz}^\infty \io nF'(n)\chi(c)\nabla c\cdot\nabla \phi
	- \int_{\tz}^\infty \io \tu\cdot \nabla n \, \phi
  \ee
  as well as
  \be{d3.4}
	-\int_{\tz}^\infty \io c\phi_t 		
	= - \int_{\tz}^\infty \io \nabla c \cdot \nabla \phi
	- \int_{\tz}^\infty \io F(n)f(c) \phi
	- \int_{\tz}^\infty \io \tu\cdot \nabla c \, \phi
  \ee
  hold for all $\phi\in C_0^\infty(\bar\Omega\times (\tz,\infty))$.
\end{defi}
The following lemma inter alia asserts that under the assumptions in (\ref{d3.2})
the integral identities (\ref{d3.3}) and (\ref{d3.4}) are indeed meaningful, and beyond this
it provides some further regularity properties of the sources, fluxes and transport terms therein. 
These properties will become essential in the proof of Lemma \ref{lem50}.
\begin{lem}\label{lem51}
  Let $\tz\ge 0$, 
  $\tu\in L^\infty_{loc}((\tz,\infty);L^2_\sigma(\Omega)) \cap L^2_{loc}((\tz,\infty);W_0^{1,2}(\Omega))$
  and $F$ be such that (\ref{F1}) holds.\abs
  i) \ If $n$ and $c$ are nonnegative and satisfy (\ref{d3.2}), then 
  \be{51.1}
	\begin{array}{l}
	nF'(n)\chi(c)\nabla c \in L^2_{loc}(\bar\Omega\times (\tz,\infty)),
	\qquad
  	\tu\cdot \nabla n \in L^\frac{5}{4}_{loc}(\bar\Omega\times (\tz,\infty)), \nn\\[2mm]
  	F(n)f(c) \in L^4_{loc}(\bar\Omega\times (\tz,\infty))
  	\qquad \mbox{and} \qquad
  	\tu\cdot\nabla c \in L^\frac{20}{11}_{loc}(\bar\Omega\times (\tz,\infty)).
	\end{array} 
  \ee
  In particular, all integrals in (\ref{d3.3}) and (\ref{d3.4}) are well-defined.\abs
  ii) \ If $(n,c)$ is a strong solution of (\ref{d3.1}) in $\Omega\times (\tz,\infty)$, then
  the identity (\ref{d3.3}) is actually valid for any $\phi\in L^5(\Omega\times (\tz,\infty))$ 
  which has compact support in $\bar\Omega\times (\tz,\infty)$, and for which
  $\nabla\phi\in L^2(\Omega\times (\tz,\infty))$ and $\phi_t\in L^\frac{4}{3}(\Omega\times (\tz,\infty))$.\abs
  iii) \ Whenever $(n,c)$ is a strong solution of (\ref{d3.1}) in $\Omega\times (\tz,\infty)$,
  the equation (\ref{d3.4}) continues to hold
  for any $\phi\in L^\frac{20}{9}(\Omega\times (\tz,\infty))$ with $\nabla\phi \in L^\frac{4}{3}(\Omega\times (\tz,\infty))$
  and $\phi_t\in L^1(\Omega\times (\tz,\infty))$, for which $\supp \phi$ is a compact subset of
  $\bar\Omega\times (\tz,\infty)$.
\end{lem}
\proof
  i) \ In view of Lemma \ref{lem511}, the regularity hypotheses on $\tu$ entail that 
  $\tu\in L^\frac{10}{3}_{loc}(\bar\Omega\times (\tz,\infty))$.
  Since $\frac{1}{4}+\frac{1}{4}=\frac{1}{2}$ and $\frac{3}{10}+\frac{1}{2}=\frac{4}{5}$ as well as
  $\frac{3}{10}+\frac{1}{4}=\frac{11}{20}$, and since $n\in L^4_{loc}(\bar\Omega\times (\tz,\infty))$,
  $c\in L^\infty_{loc}(\bar\Omega\times (\tz,\infty))$
  and $\nabla c \in L^4_{loc}(\bar\Omega\times (\tz,\infty))$ again by Definition \ref{defi3}, 
  several applications of the H\"older inequality and (\ref{F1}) readily yield the claimed integrability
  properties.\abs
  ii) and iii) \ In view of i), both statements
  can immediately be obtained upon performing standard approximation procedures.
\qed
The proof of the next lemma on a basic dissipative property of the second equation in (\ref{d3.1})
follows a testing procedure which is well-established in the context of related parabolic problems 
in their weak formulation (\cite{alt_luckhaus}), relying on the convexity of $[0,\infty)\ni s \mapsto s^p$ for $p\ge 1$.
As we are not aware of a reference precisely covering the present situation, let us include the main arguments for
completeness.
Upon slight modification, the argument can be adapted so as to extend the result to any $p\ge 1$;
for simplicity in presentation, however, we restrict ourselves to the cases $p=1$ and $p\ge 2$ relevant below,
using that then $0\le s\mapsto s^p$ belongs to $C^2([0,\infty))$.\\
For the following proof, as well as for the reasoning in Lemma \ref{lem50}, 
let us separately introduce a variant of
the cut-off function in (\ref{zeta}) defined by
\be{zeta2}
	\zeta_\delta(t):=\left\{ \begin{array}{ll}
	0 \qquad & \mbox{if } t\le t_0-\delta \mbox{ or } t\ge t_1+\delta, \\[1mm]
	\frac{t-t_0+\delta}{\delta} \qquad & \mbox{if } t\in (t_0-\delta,t_0), \\[1mm]
	1 \qquad & \mbox{if } t\in [t_0,t_1], \\[1mm]
	\frac{t_1+\delta-t}{\delta} \qquad & \mbox{if } t\in (t_1,t_1+\delta),
	\end{array} \right.
\ee
for given $t_0\in\R, t_1>t_0$ and $\delta>0$.
\begin{lem}\label{lem31}
  Let $\tz\ge 0$ and
  $\tu\in L^\infty_{loc}((\tz,\infty);L^2_\sigma(\Omega)) \cap L^2_{loc}((\tz,\infty);W_0^{1,2}(\Omega))$
  and $F$ be such that (\ref{F1}) holds,
  and suppose that $(n,c)$ is a strong solution of (\ref{d3.1}) in $\Omega\times (\tz,\infty)$.
  Then for each $p\in \{1\}\cup [2,\infty)$ there exists a null set $N(p)\subset (\tz,\infty)$ such that 
  \bea{31.1}
	\io c^p(\cdot,t)
	+ p(p-1) \int_{\tz}^t \io c^{p-2} |\nabla c|^2
	+ p \int_{\tz}^t \io F(n)f(c) c^{p-1}
	&\le& \io c^p(\cdot,t_0) \nn\\[2mm]
	& & \hspace*{-60mm} \mbox{for all $t_0\in (\tz,\infty)\setminus N(p)$ and any $t\in (t_0,\infty)\setminus N(p)$.}
  \eea
\end{lem}
\proof
  For $p \ge 1$, 
  we let $N(p)\subset(\tz,\infty)$ denote the complement of the 
  set of all Lebesgue points of $(\tz,\infty) \ni t \mapsto \io c^p(\cdot,t)$,
  and in order to prove the inequality in (\ref{31.1}) for all $t_0\in (\tz,\infty)\setminus N(p)$ and
  $t=t_1\in (t_0,\infty)\setminus N(p)$, given any such $t_0$ and $t_1$ we let $\zeta_\delta$ denote the cut-off
  function defined in (\ref{zeta2})
  for $\delta\in (0,\delta_0)$ with $\delta_0:=\min\{1,t_0-\tz\}$.
  Then according to Lemma \ref{lem51} iii), for all $h\in (0,\delta_0-\delta)$ we may apply (\ref{d3.4}) to
  $\phi(x,t):=\zeta_\delta(t) \cdot S_h[c^{p-1}](x,t)$, $(x,t)\in\Omega\times (\tz,\infty)$,
  with $p\in \{1\}\cup [2,\infty)$ and the Steklov average operator $S_h$ being defined in (\ref{S}). 
  This leads to the identity
  \bea{31.2}
	I_1(\delta,h)+I_2(\delta,h)+I_3(\delta,h)
	&:=& - \frac{1}{\delta} \int_{t_0-\delta}^{t_0} \io c \cdot S_h[c^{p-1}]
	+ \frac{1}{\delta} \int_{t_1}^{t_1+\delta} \io c\cdot S_h[c^{p-1}] \nn\\
	& & - \int_{\tz}^{t_1+1} \io \zeta_\delta(t) c(x,t) \cdot \frac{c^{p-1}(x,t+h)-c^{p-1}(x,t)}{h} dxdt \nn\\
	&=&
	- \int_{\tz}^{t_1+1} \io \zeta_\delta(t) \nabla c(x,t) \cdot \nabla S_h[c^{p-1}](x,t) dxdt \nn\\
	& & - \int_{\tz}^{t_1+1} \io \zeta_\delta(t) F(n(x,t)) f(c(x,t)) S_h[c^{p-1}](x,t) dxdt \nn\\
	& & - \int_{\tz}^{t_1+1} \io \zeta_\delta(t) (\tu(x,t)\cdot \nabla c) S_h[c^{p-1}](x,t) dxdt \nn\\[2mm]
	&=:& I_4(\delta,h)+I_5(\delta,h)+I_6(\delta,h)
  \eea
  for all $\delta\in (0,\delta_0)$ and $h\in (0,\delta_0-\delta)$, in which we observe that since 
  $(p-1) c^{p-2}$ is bounded in both cases $p=1$ and $p\ge 2$, it follows from the definition of $S_h$, the inclusion
  $\nabla c \in L^2_{loc}(\bar\Omega\times (\tz,\infty))$ and Lemma \ref{lem60} that
  \bas
	\nabla S_h[c^{p-1}]=(p-1) S_h[c^{p-2}\nabla c] \to (p-1) c^{p-2}\nabla c
	\quad \mbox{in } L^2_{loc}(\bar\Omega\times (\tz,\infty))
	\qquad \mbox{as } h\searrow 0.
  \eas
  Since Lemma \ref{lem60} also warrants that $S_h[c^{p-1}] \to c^{p-1}$ in 
  $L^\frac{20}{9}_{loc}(\bar\Omega\times (\tz,\infty))$
  as $h\searrow 0$, and since the required regularity properties of $n$ and $c$ along with (\ref{F1}) readily ensure that
  $F(n)f(c) \in L^4_{loc}(\bar\Omega\times (\tz,\infty)) 
  \subset L^\frac{20}{11}_{loc}(\bar\Omega\times (\tz,\infty))$  and also 
  $\tu\cdot\nabla c \in L^\frac{20}{11}_{loc}(\bar\Omega\times (\tz,\infty))$,
  in (\ref{31.2}) we obtain
  \bea{31.3}
	I_4(\delta,h)+I_5(\delta,h)+I_6(\delta,h)
	&\to& - (p-1) \int_{\tz}^{t_1+1} \io \zeta_\delta(t) c^{p-2} |\nabla c|^2
	-\int_{\tz}^{t_1+1} \io \zeta_\delta(t) F(n) c^{p-1} f(c) \nn\\
	& & - \int_{\tz}^{t_1+1} \io \zeta_\delta(t) c^{p-1} \tu \cdot \nabla c 
  \eea
  as $h\searrow 0$, and likewise
  \be{31.4}
	I_1(\delta,h)+I_2(\delta,h)
	\to - \frac{1}{\delta} \int_{t_0-\delta}^{t_0} \io c^p + \frac{1}{\delta} \int_{t_1}^{t_1+\delta} \io c^p
	\qquad \mbox{as } h\searrow 0.
  \ee
  As for the third integral on the left of (\ref{31.2}), we estimate using Young's inequality to find upon a substitution that
  \bas
	I_3 &\ge& - \frac{1}{h} \cdot \bigg\{ 
	\frac{1}{p} \int_{\tz}^{t_1+1} \io \zeta_\delta(t) c^p(x,t) dxdt
	+ \frac{p-1}{p} \int_{\tz}^{t_1+1} \io \zeta_\delta(t) c^p(x,t+h) dxdt \bigg\} \nn\\
	& & + \frac{1}{h} \int_{\tz}^{t_1+1} \io \zeta_\delta(t) c^p(x,t) dxdt \nn\\
	&=& \frac{p-1}{p} \int_{\tz}^{t_1+1} \io \frac{\zeta_\delta(t)-\zeta_\delta(t-h)}{h} \cdot c^p(x,t) dxdt
  \eas
  for all $\delta\in (0,\delta_0)$ and $h\in (0,\delta_0-\delta)$, so that by the dominated convergence theorem we conclude 
  that
  \bea{31.5}
	\liminf_{h\searrow 0} I_3(\delta,h) 
	\ge \frac{p-1}{p} \int_{\tz}^{t_1+1} \io \zeta_\delta'(t) c^p
	= \frac{p-1}{p\delta} \int_{t_0-\delta}^{t_0} \io c^p 
	- \frac{p-1}{p\delta} \int_{t_1}^{t_1+\delta} c^p
  \eea
  as $h\searrow 0$. Since $\tu$ is solenoidal and hence 
  \bas
	- \int_{\tz}^{t_1+1} \io \zeta_\delta(t) c^{p-1} \tu \cdot \nabla c 
	= - \frac{1}{p} \int_{\tz}^{t_1+1} \io \zeta_\delta(t) \tu \cdot \nabla c^p = 0,
  \eas
  combining (\ref{31.2})-(\ref{31.5}) and rearranging shows that
  \bas
	\frac{1}{p\delta} \int_{t_1}^{t_1+\delta} \io c^p
	+ (p-1) \int_{\tz}^{t_1+1} \io \zeta_\delta(t) c^{p-2} |\nabla c|^2
	+ \int_{\tz}^{t_1+1} \io \zeta_\delta(t) F(n) c^{p-1} f(c)
	\le \frac{1}{p\delta} \int_{t_0-\delta}^{t_0} c^p
  \eas
  for all $\delta\in (0,\delta_0)$.
  Since $\zeta_\delta\equiv 1$ in $(t_0,t_1)$, thanks to the assumed Lebesgue point properties of $t_0$ and $t_1$ 
  this readily yields the desired inequality for such $t_0$ and $t_1$ on taking $\delta\searrow 0$.
\qed
Evaluating (\ref{31.1}) for $p=1$ and $p=2$ as well as in the limit case $p\to\infty$ we obtain
the following consequence which provides some first, still quite weak, information on decay of $c$,
at least under the assumption that $n$, and hence $F(n)$ by (\ref{F2}), remains positive in an 
appropriate sense.
\begin{cor}\label{cor32}
  Let $\tz\ge 0$ and
  $\tu\in L^\infty_{loc}([\tz,\infty);L^2_\sigma(\Omega)) \cap L^2_{loc}([\tz,\infty);W_0^{1,2}(\Omega))$
  and $F$ be such that (\ref{F1}) holds,
  and suppose that $(n,c)$ is a strong solution of (\ref{d3.1}) in $\Omega\times (\tz,\infty)$.
  Then 
  \be{32.1}
	\int_{\tz}^\infty \io F(n)f(c) \le {\rm{ess}}\liminf_{\hspace*{-5mm} t\searrow \tz} \io c(\cdot,t)
  \ee
  and
  \be{32.2}
	\int_{\tz}^\infty \io |\nabla c|^2 \le {\rm{ess}}\liminf_{\hspace*{-5mm} t\searrow \tz} 
	\frac{1}{2} \io c^2(\cdot,t),
  \ee
  and there exists a null set $N\subset (\tz,\infty)$ such that
  \be{cinfty}
	\|c(\cdot,t)\|_{L^\infty(\Omega)} \le \|c(\cdot,t_0)\|_{L^\infty(\Omega)}
	\qquad \mbox{for all $t_0\in (\tz,\infty)\setminus N$ and any $t\in(t_0,\infty)\setminus N$.}
  \ee
\end{cor}
\proof
  The inequalities in (\ref{32.1}) and (\ref{32.2}) immediately result from applying Lemma \ref{lem31}
  to $p:=1$ and $p:=2$, respectively. 
  We next invoke Lemma \ref{lem31} for $p_j:=j$ for $j\in\N$ to obtain null sets $N(p_j)$ with the properties
  listed there. For the null set $N:=\bigcup_{j\in\N} N(p_j) \subset (\tz,\infty)$ we thus obtain that
  $\|c(\cdot,t)\|_{L^{p_j}(\Omega)} \le \|c(\cdot,t_0)\|_{L^{p_j}(\Omega)}$ whenever $j\in\N$ as well as
  $t_0\in [\tz,\infty)\setminus N$ and $t\in (t_0,\infty)\setminus N$, because $F$ and $f$ are both nonnegative. 
  Taking $j\to\infty$ here shows that $\|c(\cdot,t)\|_{L^\infty(\Omega)} \le \|c(\cdot,t_0)\|_{L^\infty(\Omega)}$
  for any such $t_0$ and $t$, which yields (\ref{cinfty}).
\qed
\mysection{A doubly uniform decay property of solutions to (\ref{d3.1})}\label{sect4}
In this section we shall make sure that any solution $(n,c)$ of (\ref{d3.1}) satisfies
$c(\cdot,t)\to 0$ in $L^\infty(\Omega)$ as $t\to\infty$, which will be a fundamental ingredient
for our later regularity arguments. 
Since apart from considering any fixed eventual energy solution of (\ref{0}) 
we wish to address the entire family of solutions to the approximate problems (\ref{0eps}) for $\eps\in (0,1)$, 
our plan will be to make sure that this convergence is actually uniform with respect to the choice of
the considered solution, as well as of $\tu$ and $F$, in an appropriate sense.
In order to make this more precise, let us introduce the following notation.
\begin{defi}\label{defi_set}
  Given $m>0, M>0, L>0$ and $\tz\ge 0$, we let 
  \be{set}
	\set
  \ee
  denote the set of all triples $(n,c,F)$
  of functions $n:\Omega\times (\tz,\infty)\to\R, c: \Omega\times (0,\infty)\to \R$ and $F:[0,\infty)\to \R$
  such that $F$ satisfies (\ref{F1}) and (\ref{F2}), and that for some
  $\tu\in L^\infty_{loc}((\tz,\infty);L^2_\sigma(\Omega)) \cap L^2_{loc}((\tz,\infty);W_0^{1,2}(\Omega))$, 
  the pair $(n,c)$ is a strong solution of (\ref{d3.1}) in $\Omega\times (\tz,\infty)$
  satisfying
  \be{set1}
	\io n(\cdot,t)=m
	\quad \mbox{and} \quad
	\|c(\cdot,t)\|_{L^\infty(\Omega)} \le M
	\qquad \mbox{for a.e.~} t>\tz
  \ee
  as well as
  \be{set2}
	\int_t^{t+1} \io \Big\{\frac{|\nabla n|^2}{n} + |\nabla c|^4 \Big\} \le L 
	\qquad \mbox{for all } t>\tz.
  \ee
\end{defi}
In this framework, an adequate interpretation of (\ref{32.1}) and (\ref{32.2}) in Corollary \ref{cor32} yields the following.
\begin{lem}\label{lem43}
  Let $m>0, M>0, L>0$ and $\tz\ge 0$. Then the set $\set$ in (\ref{set}) has the property that	
  \be{43.1}
	\sup_{(n,c,F)\in \set} \inf_{t\in [\tz,\tz+\tau]} \int_t^{t+1} \io \Big\{ F(n)f(c) + |\nabla c|^2 \Big\} 
	\to 0 \qquad \mbox{as } \tau\to\infty.
  \ee
\end{lem}
\proof
  In order to verify (\ref{43.1}) we let $\delta>0$ be given and pick any integer $k$ fulfilling $k>\frac{C_1}{\delta}$, where
  $C_1:=M|\Omega| + \frac{M^2|\Omega|}{2}$.
  We claim that then for each $(n,c,F) \in \set$ we have
  \be{43.2}
	\inf_{t\in [\tz,\tz+k]} \int_t^{t+1} \io \Big\{ F(n)f(c) + |\nabla c|^2 \Big\} < \delta.
  \ee
  To see this, given any such $(n,c)$ we first apply Corollary \ref{cor32} to obtain
  \bas
	\int_{\tz}^\infty \io F(n)f(c) \le {\rm{ess}}\liminf_{\hspace*{-5mm} t\searrow \tz} \io c(\cdot,t) \le M|\Omega|
  \eas
  and
  \bas
	\int_{\tz}^\infty \io |\nabla c|^2 \le {\rm{ess}}\liminf_{\hspace*{-5mm} t\searrow \tz} \frac{1}{2} \io c^2(\cdot,t)
	\le \frac{M^2|\Omega|}{2},
  \eas
  whence for $h(t):=\io F(n(\cdot,t))f(c(\cdot,t)) + \io |\nabla c(\cdot,t)|^2$, $t>\tz$, we obtain
  \be{43.3}
	\int_{\tz}^\infty h(t)dt \le C_1
  \ee
  by definition of $C_1$.\\
  Now if (\ref{43.2}) was false, then (\ref{43.3}) would imply that
  \bas
	C_1 \ge \int_{\tz}^{\tz+k} h(t)dt
	= \sum_{j=1}^k \int_{\tz+j-1}^{\tz+j} h(t)dt
	\ge \sum_{j=1}^k \delta =k \delta,
  \eas
  which is absurd in view of our choice of $k$.
  As $(n,c,F)\in\set$ and $\delta>0$ were arbitrary, this establishes (\ref{43.1}).
\qed
Now a crucial point appears to consist in deriving more substantial decay properties from this without any 
further knowledge on possible lower bounds for $n$ beyond the weak information that its mass $\io n$ remains
constantly positive by definition of $\set$. 
In order to prepare a first step in this direction, we state an elementary observation which will below
be related to a lower estimate for the first integral appearing in (\ref{43.1}).
\begin{lem}\label{lem433}
  Assume that $F$ satisfies (\ref{F1}) and (\ref{F2}). Then for all $m>0$ and $B\ge \frac{m}{8}$,
  each nonne\-gative $\varphi\in L^3(\Omega)$ fulfilling
  \be{433.1}
	\io \varphi=m
	\qquad \mbox{and} \qquad
	\io \varphi^3 \le B
  \ee
  has the property that
  \be{433.2}
	\io F(\varphi) \ge \sqrt{\frac{m^3}{128B}}.
  \ee
\end{lem}
\proof
  As $B\ge \frac{m}{8}$, the number $C_1:=\sqrt{\frac{8B}{m}}$ satisfies $C_1\ge 1$, 
  so that combining (\ref{F2}) with (\ref{F1}) shows that
  \bas
	F(s) \ge \min \Big\{\frac{s}{2} \, , \, \frac{1}{2} \Big\}
	\ge \frac{s}{2C_1}
	\qquad \mbox{for all } s\in [0,C_1].
  \eas
  Hence, given any nonnegatve $\varphi\in L^3(\Omega)$ fulfilling (\ref{433.1}), we have
  \be{433.3}
	\io F(\varphi) \ge \int_{\{\varphi\le C_1\}} F(\varphi)
	\ge \frac{1}{2C_1} \int_{\{\varphi\le C_1\}} \varphi.
  \ee
  In order to further estimate the latter integral, we use the H\"older inequality and (\ref{433.1}) to obtain
  \bas
	\int_{\{\varphi>C_1\}} \varphi
	\le \Big( \io \varphi^3 \Big)^\frac{1}{3} \cdot \Big| \{\varphi>C_1\} \Big|^\frac{2}{3}
	\le B^\frac{1}{3} \cdot \Big|\{\varphi>C_1\}\Big|^\frac{2}{3},
  \eas
  so that since $|\{\varphi>C_1\}| \le \frac{m}{C_1}$ by the Chebyshev inequality, we conclude that
  \bas
	\int_{\{\varphi\le C_1\}} \varphi = m - \int_{\{\varphi>C_1\}} \varphi
	\ge m-B^\frac{1}{3} \cdot \Big|\{\varphi>C_1\}\Big|^\frac{2}{3}
	\ge m - B^\frac{1}{3} \cdot \Big(\frac{m}{C_1}\Big)^\frac{2}{3}
	= m-B^\frac{1}{3} \cdot \bigg( \frac{m}{\sqrt{\frac{8B}{m}}} \bigg)^\frac{2}{3}
	= \frac{m}{2}.
  \eas
  Therefore, (\ref{433.2}) results from (\ref{433.3}).
\qed
By means of appropriate interpolation
making use of the regularity property (\ref{set2}) jointly shared by all elements of $\set$, 
we can thereby turn Lemma \ref{lem43}
into a statement on decay of $f(c)$ which is uniform with respect to functions from this set.
\begin{lem}\label{lem44}
  Let $m>0, M>0, L>0$ and $\tz\ge 0$. Then
  \bas
	\sup_{(n,c,F)\in\set} 
	\inf_{\begin{array}{c} \scriptstyle S\subset (T_0,T_0+\tau) \\[-1mm]	
	\scriptstyle S \mbox{ \rm \small is measurable with } |S|\ge \frac{1}{2} \\[-1mm]
	\scriptstyle \rm{ and \ diam} \, S \le 1
	\end{array} }
	\int_S \io f(c(x,t))dxdt \to 0
	\qquad \mbox{as } \tau\to\infty,
  \eas
  where $\set$ is as defined in (\ref{set}).
\end{lem}
\proof
  We need to show that for each fixed $m,M,L$ and $\tz$, given $\delta>0$ we can find $\tau>0$ with the property that for any
  $(n,c,F)\in \set$ there exists a measurable set $S\subset (\tz,\tz+\tau)$ such that $|S|\ge\frac{1}{2}$ and
  $\rm{diam} \, S \le 1$ as well as
  \be{44.1}
	\int_S \io f(c(x,t)) dxdt < \delta.
  \ee
  In order to prepare our definition of $\tau$, let us first make use of the embedding $W^{1,2}(\Omega)\hra L^6(\Omega)$, 
  which in conjunction with the Poincar\'e inequality yields $C_1>0$ and $C_2>0$ such that
  \be{44.01}
	\|\varphi\|_{L^6(\Omega)}^2
	\le C_1\|\nabla \varphi\|_{L^2(\Omega)}^2 
	+ C_1 \|\varphi\|_{L^2(\Omega)}^2
	\qquad \mbox{for all } \varphi\in W^{1,2}(\Omega)
  \ee
  and
  \be{44.2}
	\|\varphi-\overline{\varphi}\|_{L^6(\Omega)}
	\le C_2\|\nabla \varphi\|_{L^2(\Omega)}
	\qquad \mbox{for all } \varphi\in W^{1,2}(\Omega),
  \ee
  where again we have set $\overline{\varphi}:=\mint_{\Omega} \varphi$ for $\varphi\in L^1(\Omega)$.
  Next, 		
  the Gagliardo-Nirenberg inequality provides $C_3>0$ fulfilling
  \be{44.3}
	\|\varphi\|_{L^\frac{12}{5}(\Omega)}^4 \le C_3\|\nabla\varphi\|_{L^2(\Omega)} \|\varphi\|_{L^2(\Omega)}^3
	+ C_3 \|\varphi\|_{L^2(\Omega)}^4
	\qquad \mbox{for all } \varphi\in W^{1,2}(\Omega).
  \ee
  Finally abbreviating $C_4:=\|f'\|_{L^\infty((0,M))}$, $B:=8 \cdot \Big\{ \frac{C_1 L}{4} + C_1 m\Big\}^3$
  and $C_5:=\sqrt{\frac{m^3}{128B}}$, given $\delta>0$ we can find some small $\delta_0>0$ such that
  \be{44.4}
	\frac{|\Omega|}{C_5} \cdot \delta_0 < \frac{\delta}{2}
  \ee
  and 
  \be{44.5}
	\frac{C_2 C_4 |\Omega|}{C_5} \cdot \bigg\{ \frac{C_3 m^\frac{3}{2} L^\frac{1}{2}}{2} + C_3 m^2 \bigg\}^\frac{1}{2}
	\cdot \delta_0^\frac{1}{2} < \frac{\delta}{2}.
  \ee
  Then Lemma \ref{lem43} says that there exists $\tau>0$ such that
  \bas
	\sup_{(n,c,F)\in \set} \inf_{t\in [\tz,\tz+\tau]} \int_t^{t+1} \io \Big\{ F(n)f(c)+|\nabla c|^2 \Big\} < \delta_0.
  \eas
  Thus, if we now pick any $(n,c,F)\in \set$, then we can pick $t_0\in [\tz,\tz+\tau]$ such that
  \be{44.7}
	\int_{t_0}^{t_0+1} \io F(n)f(c) < \delta_0
	\qquad \mbox{and} \qquad
	\int_{t_0}^{t_0+1} \io |\nabla c|^2 < \delta_0.
  \ee
  With this number $t_0$ fixed henceforth, we observe that
  \be{44.777}
	\int_{t_0}^{t_0+1} \|\nabla n^\frac{1}{2}(\cdot,t)\|_{L^2(\Omega)}^2 dt
	= \frac{1}{4} \int_{t_0}^{t_0+1} \io \frac{|\nabla n|^2}{n} \le \frac{L}{4}
  \ee
  by definition of $\set$, so that using (\ref{44.01}) we obtain
  \bas
	\int_{t_0}^{t_0+1} \|n(\cdot,t)\|_{L^3(\Omega)} dt
	&=& \int_{t_0}^{t_0+1} \|n^\frac{1}{2}(\cdot,t)\|_{L^6(\Omega)}^2 dt \\
	&\le& C_1 \int_{t_0}^{t_0+1} \|\nabla n^\frac{1}{2}(\cdot,t)\|_{L^2(\Omega)}^2 dt 
	+ C_1 \int_{t_0}^{t_0+1} \|n^\frac{1}{2}(\cdot,t)\|_{L^2(\Omega)}^2 dt \\
	&\le& \frac{C_1 L}{4} + C_1 m \\
	&=& \frac{B^\frac{1}{3}}{2},
  \eas
  because
  \be{44.77}
	\|n^\frac{1}{2}(\cdot,t)\|_{L^2(\Omega)}^2 = \io n(\cdot,t)=m
	\qquad \mbox{for a.e.~} t>\tz
  \ee
  again due to the definition of $\set$.
  The measurable set
  \bas
	S:=\Big\{ t\in (t_0,t_0+1) \ \Big| \ \|n(\cdot,t)\|_{L^3(\Omega)} \le B^\frac{1}{3}\Big\}
  \eas
  therefore satisfies $|S| \ge \frac{1}{2}$ by the Chebyshev inequality, and in view of our definition of $C_5$ we infer from 
  Lemma \ref{lem433} that
  \be{44.8}
	\io F(n(x,t)) dx
	\ge C_5
	\qquad \mbox{for all } t\in S.
  \ee
  For the proof of (\ref{44.1}), we now decompose the first integral in (\ref{44.7}) according to
  \bas
	\int_{t_0}^{t_0+1} \io F(n)f(c)
	&=& \int_{t_0}^{t_0+1} \io F(n(x,t)) \cdot \Big\{ f(c(x,t)) - \overline{f(c(\cdot,t))} \Big\} dxdt \\
	& & + \int_{t_0}^{t_0+1} \io F(n(x,t)) \cdot \overline{f(c(\cdot,t))} dxdt,
  \eas
  where by (\ref{44.8}),
  \bas
	\int_{t_0}^{t_0+1} \io F(n(x,t)) \cdot \overline{f(c(\cdot,t))} dxdt
	&=& \int_{t_0}^{t_0+1}  \overline{f(c(\cdot,t))} \cdot \bigg\{ \io F(n(x,t))dx \bigg\} dt \\
	&\ge& \int_S \overline{f(c(\cdot,t))} \cdot \bigg\{ \io F(n(x,t))dx \bigg\} dt \\
	&\ge& C_5 \int_S \overline{f(c(\cdot,t))} dt \\
	&=& \frac{C_5}{|\Omega|} \int_S \io f(c(x,t)) dxdt.
  \eas
  Hence, for the integral in question we obtain the inequality
  \bea{44.9}
	\int_S \io f(c(x,t)) dxdt
	&\le& \frac{|\Omega|}{C_5} \int_{t_0}^{t_0+1} \io F(n)f(c) \nn\\
	& & - \frac{|\Omega|}{C_5} \int_{t_0}^{t_0+1} \io F(n(x,t)) \cdot \Big\{f(c(x,t))-\overline{f(c(\cdot,t))}\Big\} dxdt,
  \eea
  in which thanks to (\ref{44.7}) and (\ref{44.4}),
  \be{44.10}
	\frac{|\Omega|}{C_5} \int_{t_0}^{t_0+1} \io F(n)f(c)
	<\frac{|\Omega|}{C_5} \cdot \delta_0 < \frac{\delta}{2}.
  \ee
  Moreover, invoking the H\"older inequality we can estimate
  \bea{44.11}
	& & \hspace*{-20mm}
	- \frac{|\Omega|}{C_5} \int_{t_0}^{t_0+1} \io F(n(x,t)) \cdot \Big\{f(c(x,t))-\overline{f(c(\cdot,t))}\Big\} dxdt 
	\nn\\
	&\le& \frac{|\Omega|}{C_5} \int_{t_0}^{t_0+1} \Big\|F(n(\cdot,t))\Big\|_{L^\frac{6}{5}(\Omega)} \cdot
	\Big\| f(c(\cdot,t))- \overline{f(c(\cdot,t))} \Big\|_{L^6(\Omega)} dt \nn\\
	&\le& \frac{|\Omega|}{C_5} \cdot 
	\bigg( \int_{t_0}^{t_0+1} \Big\|F(n(\cdot,t))\Big\|_{L^\frac{6}{5}(\Omega)}^2 dt \bigg)^\frac{1}{2} \cdot
	\bigg( \int_{t_0}^{t_0+1} \Big\| f(c(\cdot,t))- \overline{f(c(\cdot,t))} \Big\|_{L^6(\Omega)}^2 dt \bigg)^\frac{1}{2},
  \eea
  where as a consequence of (\ref{44.2}) and our choice of $C_4$ we have
  \bea{44.12}
	\Big\|f(c(\cdot,t))-\overline{f(c(\cdot,t))} \Big\|_{L^6(\Omega)}^2
	&\le& C_2^2 \Big\|\nabla f(c(\cdot,t))\Big\|_{L^2(\Omega)}^2 \nn\\
	&=& C_2^2 \io f'^2(c(\cdot,t)) |\nabla c(\cdot,t)|^2 \nn\\
	&\le& C_2^2 C_4^2 \|\nabla c(\cdot,t)\|_{L^2(\Omega)}^2
	\qquad \mbox{for a.e.~} t\in (t_0,t_0+1),
  \eea
  because $c\le M$ a.e.~in $\Omega\times (\tz,\infty)$ by definition of $\set$.\\
  As for the factor in (\ref{44.11}) containing $n$, we use (\ref{44.3}) and the fact that $F(n)\le n$ by (\ref{F1})
  to see that
  \bea{44.13}
	\Big\|F(n(\cdot,t))\Big\|_{L^\frac{6}{5}(\Omega)}^2
	&\le& \|n(\cdot,t)\|_{L^\frac{6}{5}(\Omega)}^2 \nn\\
	&=& \|n^\frac{1}{2}(\cdot,t)\|_{L^\frac{12}{5}(\Omega)}^4 \nn\\
	&\le& C_3 \|\nabla n^\frac{1}{2}(\cdot,t)\|_{L^2(\Omega)} \|n^\frac{1}{2}(\cdot,t)\|_{L^2(\Omega)}^3
	+ C_3 \|n^\frac{1}{2}(\cdot,t)\|_{L^2(\Omega)}^4 \nn\\
	&=& C_3 m^\frac{3}{2} \|\nabla n^\frac{1}{2}(\cdot,t)\|_{L^2(\Omega)} + C_3 m^2
	\qquad \mbox{for a.e.~} t\in (t_0,t_0+1),
  \eea
  again due to (\ref{44.77}). \\
  In summary, from (\ref{44.11}), (\ref{44.12}) and (\ref{44.13}) we obtain upon employing 
  the Cauchy-Schwarz inequality that
  \bas
	& & \hspace*{-10mm}
	- \frac{|\Omega|}{C_5} \int_{t_0}^{t_0+1} \io F(n(x,t)) \cdot \Big\{f(c(x,t))-\overline{f(c(\cdot,t))}\Big\} dxdt 
	\nn\\
	&\le& \frac{|\Omega|}{C_5} \cdot 
	\Bigg\{ C_3 m^\frac{3}{2} \int_{t_0}^{t_0+1} \|\nabla n^\frac{1}{2}(\cdot,t)\|_{L^2(\Omega)} dt + C_3 m^2 
	\Bigg\}^\frac{1}{2}  \cdot
	\Bigg\{ C_2^2 C_4^2 \int_{t_0}^{t_0+1} \|\nabla c(\cdot,t)\|_{L^2(\Omega)}^2 dt \Bigg\}^\frac{1}{2} \\[2mm]
	&\le& \frac{|\Omega|}{C_5} \cdot 
	\Bigg\{ C_3 m^\frac{3}{2} \cdot 
	\bigg(\int_{t_0}^{t_0+1} \|\nabla n^\frac{1}{2}(\cdot,t)\|_{L^2(\Omega)}^2 dt \bigg)^\frac{1}{2}
	+ C_3 m^2 \Bigg\}^\frac{1}{2} \cdot
	\Bigg\{ C_2^2 C_4^2 \int_{t_0}^{t_0+1} \|\nabla c(\cdot,t)\|_{L^2(\Omega)}^2 dt \Bigg\}^\frac{1}{2},
  \eas
  so that (\ref{44.777}), (\ref{44.7}) and then (\ref{44.5}) become applicable so as to warrant that
  \bas
	& & \hspace*{-40mm}
	- \frac{|\Omega|}{C_5} \int_{t_0}^{t_0+1} \io F(n(x,t)) \cdot \Big\{f(c(x,t))-\overline{f(c(\cdot,t))}\Big\} dxdt 
	\nn\\
	&\le& \frac{|\Omega|}{C_5} \cdot \Big\{ C_3 m^\frac{3}{2} \cdot \frac{L^\frac{1}{2}}{2} + C_3 m^2 
	\Big\}^\frac{1}{2}
	\cdot \Big\{ C_2^2 C_4^2 \delta_0\Big\}^\frac{1}{2} 	\nn\\[1mm]
	&<& \frac{\delta}{2}.
  \eas
  Combined with (\ref{44.10}) and (\ref{44.9}), this shows (\ref{44.1}) and thereby completes the proof.
\qed
Once more relying on the regularity features in (\ref{set2}), we can show that the above convergence of $f(c)$
does not only take place in $L^1(\Omega)$ but actually even in $L^\infty(\Omega)$.
\begin{lem}\label{lem45}
  Let $m>0, M>0, L>0$ and $\tz\ge 0$. Then
  \bas
	\sup_{(n,c,F)\in\set} 
	\inf_{\begin{array}{c} \scriptstyle S\subset (T_0,T_0+\tau) \\[-1mm]	
	\scriptstyle S \mbox{ \rm \small is measurable with } |S|\ge \frac{1}{2} \\[-1mm]
	\scriptstyle \rm{and \ diam} \, S \le 1
	\end{array} }
	\int_S \|f(c(\cdot,t))\|_{L^\infty(\Omega)} \to 0
	\qquad \mbox{as } \tau\to\infty.
  \eas
\end{lem}
\proof
  Fixing $m>0, M>0, L>0$ and $\tz\ge 0$, we need to make sure that for each $\delta>0$ there exists $\tau>0$ 
  such that given any $(n,c,F)\in\set$ we can find a measurable $S\subset (\tz,\tz+\tau)$ such that
  $|S|\ge \frac{1}{2}$ and $\rm{diam} \, S \le 1$ as well as
  \be{45.1}
	\int_S \|f(c(\cdot,t))\|_{L^\infty(\Omega)} dt < \delta.
  \ee
  For this purpose, we first use that $W^{1,4}(\Omega) \hra L^\infty(\Omega)$ to interpolate by means of the 
  Gagliardo-Nirenberg inequality to find $C_1>0$ fulfilling
  \be{45.2}
	\|\varphi\|_{L^\infty(\Omega)} 
	\le C_1 \|\nabla\varphi\|_{L^4(\Omega)}^\frac{12}{13} \|\varphi\|_{L^1(\Omega)}^\frac{1}{13} 
	+ C_1\|\varphi\|_{L^1(\Omega)}
	\qquad \mbox{for all } \varphi\in W^{1,4}(\Omega),
  \ee
  and abbreviate $C_2:=\|f'\|_{L^\infty((0,M))}$.
  Then for arbitrary $\delta>0$ we can pick $\delta_0>0$ small enough satisfying
  \be{45.4}
	C_1\delta_0<\frac{\delta}{2}
  \ee
  and
  \be{45.5}
	C_1 C_2^\frac{12}{13} L^\frac{3}{13} \delta_0^\frac{1}{13}
	<\frac{\delta}{2},
  \ee
  and apply Lemma \ref{lem44} to obtain $\tau>0$ such that
  \bas
	\sup_{(n,c,F)\in\set} 
	\inf_{\begin{array}{c} \scriptstyle S\subset (T_0,T_0+\tau) \\[-1mm]	
	\scriptstyle S \mbox{ \rm \small is measurable with } |S|\ge \frac{1}{2} \\[-1mm]
	\scriptstyle \rm{and \ diam} \, S \le 1
	\end{array} }
	\int_S \io f(c(x,t))dxdt <\delta_0.
  \eas
  This means that if we fix $(n,c,F)\in \set$, then we can find a measurable $S\subset (\tz,\tz+\tau)$ such that
  $|S|\ge \frac{1}{2}, \rm{diam} \, S \le 1$ and
  \be{45.6}
	\int_S \io f(c(x,t))dxdt < \delta_0.
  \ee
  This entails that if we estimate the integral under consideration by using (\ref{45.2}) according to
  \be{45.7}
	\int_S \|f(c(\cdot,t))\|_{L^\infty(\Omega)} dt
	\le C_1 \int_{t_0}^{t_0+1} \|\nabla f(c(\cdot,t))\|_{L^4(\Omega)}^\frac{12}{13} 
	\|f(c(\cdot,t))\|_{L^1(\Omega)}^\frac{1}{13} dt
	+ C_1\int_{t_0}^{t_0+1} \|f(c(\cdot,t))\|_{L^1(\Omega)} dt,
  \ee
  then the rightmost term herein satisfies
  \be{45.8}
	C_1\int_{t_0}^{t_0+1} \|f(c(\cdot,t))\|_{L^1(\Omega)} dt
	< C_1\delta_0 < \frac{\delta}{2}
  \ee
  in view of (\ref{45.4}).
  As the inclusion $(n,c,F)\in\set$ ensures that $c\le M$ a.e.~in $\Omega\times (\tz,\infty)$ and hence
  \bas
	\|\nabla f(c(\cdot,t))\|_{L^4(\Omega)}^4
	= \io f'^4(c(\cdot,t)) |\nabla c(\cdot,t)|^4
	\le C_2^4 \io |\nabla c(\cdot,t)|^4
	\qquad \mbox{for a.e.~} t>\tz
  \eas
  by definition of $C_2$, two applications of the H\"older inequality to the first integral on the right of (\ref{45.7})
  yield
  \bas
	& & \hspace*{-20mm}
	C_1 \int_{t_0}^{t_0+1} \|\nabla f(c(\cdot,t))\|_{L^4(\Omega)}^\frac{12}{13} 
	\|f(c(\cdot,t))\|_{L^1(\Omega)}^\frac{1}{13} dt	\\
	&\le& C_1 \bigg\{ \int_{t_0}^{t_0+1} \|\nabla f(c(\cdot,t))\|_{L^4(\Omega)}^4 dt \bigg\}^\frac{3}{13}
	\cdot \bigg\{\int_{t_0}^{t_0+1} \|f(c(\cdot,t))\|_{L^1(\Omega)}^\frac{1}{10} dt \bigg\}^\frac{10}{13} \\
	&\le& C_1 \bigg\{ \int_{t_0}^{t_0+1} \|\nabla f(c(\cdot,t))\|_{L^4(\Omega)}^4 dt \bigg\}^\frac{3}{13}
	\cdot \bigg\{\int_{t_0}^{t_0+1} \|f(c(\cdot,t))\|_{L^1(\Omega)}dt \bigg\}^\frac{1}{13} \\
	&\le& C_1 \cdot \bigg\{ C_2^4 \int_{t_0}^{t_0+1} \io |\nabla c|^4 \bigg\}^\frac{3}{13}
	\cdot \bigg\{\int_{t_0}^{t_0+1} \|f(c(\cdot,t))\|_{L^1(\Omega)}dt \bigg\}^\frac{1}{13}.
  \eas
  Since $\int_{t_0}^{t_0+1} \io |\nabla c|^4 \le L$ by definition of $\set$, (\ref{45.6}) and (\ref{45.5}) guarantee that
  \bas
	C_1 \int_{t_0}^{t_0+1} \|\nabla f(c(\cdot,t))\|_{L^4(\Omega)}^\frac{12}{13} 
	\|f(c(\cdot,t))\|_{L^1(\Omega)}^\frac{1}{13} dt
	\le C_1 \cdot \Big\{ C_2^4 L \Big\}^\frac{3}{13} \cdot \delta_0^\frac{1}{13}
	< \frac{\delta}{2},
  \eas
  which in conjunction with (\ref{45.8}) and (\ref{45.7}) establishes (\ref{45.1}).
\qed
Now thanks to the assumed positivity of $f$ on $(0,\infty)$, and in view of the downward monotonicity
of $t\mapsto \|c(\cdot,t)\|_{L^\infty(\Omega)}$ asserted by Corollary \ref{cor32}, the latter implies 
the following doubly uniform decay property of (\ref{d3.1}) which constitutes the main result of this section.
\begin{lem}\label{lem46}
  Let $m>0, M>0, L>0$ and $\tz\ge 0$. Then
  \be{46.1}
	\sup_{(n,c,F)\in\set} \|c\|_{L^\infty(\Omega\times (t,\infty))} \to 0
	\qquad \mbox{as } t\to\infty.
  \ee
\end{lem}
\proof
  We fix $m>0, M>0, L>0$ and $\tz\ge 0$ and note that proving (\ref{46.1}) amounts to showing that for each $\delta>0$
  we can find $t_0>\tz$ such that whenever $(n,c,F)\in\set$, we have
  \be{46.2}
	\|c(\cdot,t)\|_{L^\infty(\Omega)} \le \delta
	\qquad \mbox{for a.e.~} t>t_0.
  \ee
  To verify this, we may assume that $\delta<M$ and then observe that since $f$ is continuous and positive on $(0,\infty)$,
  the number
  \bas
	\delta_0:=\min \Big\{ f(s) \ \Big| \ s\in [\delta,M] \Big\}
  \eas
  is well-defined and positive.
  An application of Lemma \ref{lem45} thus yields $\tau>0$ with the property that
  \be{46.3}
	\sup_{(n,c,F)\in\set} 
	\inf_{\begin{array}{c} \scriptstyle S\subset (T_0,T_0+\tau) \\[-1mm]	
	\scriptstyle S \mbox{ \rm \small is measurable with } |S|\ge \frac{1}{2} \\[-1mm]
	\scriptstyle \rm{and \ diam} \, S \le 1
	\end{array} }
	\int_S \|f(c(\cdot,t))\|_{L^\infty(\Omega)}
	<\frac{\delta_0}{4}.
  \ee
  In order to show that the desired conclusion holds for $t_0:=\tz+\tau$, we now fix $(n,c,F)\in\set$
  and then obtain from (\ref{46.3}) that there exists a measurable set $S\subset (\tz,\tz+\tau)$ such that
  $|S|\ge\frac{1}{2}, \rm{diam} \, S\le 1$ and
  \bas
	\int_S \|f(c(\cdot,t))\|_{L^\infty(\Omega)} <\frac{\delta_0}{4}.
  \eas
  Since
  \bas
	\delta_0 \cdot \bigg| \Big\{ t\in S \ \Big| \ \|f(c(\cdot,t))\|_{L^\infty(\Omega)} \ge \delta_0\Big\} \bigg|
	\le \int_S \|f(c(\cdot,t))\|_{L^\infty(\Omega)} dt
  \eas
  by the Chebyshev inequality, this guarantees that
  \bas
	S_0:=\Big\{ t\in S \ \Big| \ \|f(c(\cdot,t))\|_{L^\infty(\Omega)} < \delta_0 \Big\}
  \eas
  satisfies $|S_0| \ge |S|-\frac{1}{4} \ge \frac{1}{4}$.
  Moreover, the definition of $\delta_0$ ensures that for each $t\in S_0$ we necessarily have 
  $c(\cdot,t)\le \delta$ a.e.~in $\Omega$, that is,
  \be{46.4}
	\|c(\cdot,t_1)\|_{L^\infty(\Omega)} \le \delta
	\qquad \mbox{for all } t_1\in S_0.
  \ee
  We now invoke Corollary \ref{cor32} to find a null set $N\subset (\tz,\infty)$ such that
  $\|c(\cdot,t)\|_{L^\infty(\Omega)} \le \|c(\cdot,t_1)\|_{L^\infty(\Omega)}$ for all $t_1\in (\tz,\infty)\setminus N$
  and each $t\in (t_1,\infty)\setminus N$.
  Since $|S_0\setminus N| \ge \frac{1}{4}$, we may thus pick an arbitrary $t_1\in S_0\setminus N$ and apply (\ref{46.4})
  to infer that
  \bas
	\|c(\cdot,t)\|_{L^\infty(\Omega)} \le \|c(\cdot,t_1)\|_{L^\infty(\Omega)} \le \delta
	\qquad \mbox{for all } t\in (t_1,\infty)\setminus N,
  \eas
  which implies (\ref{46.2}) due to the fact that the inclusion $t_1\in S_0\subset S \subset (\tz,\tz+\tau)$
  along with our choice of $t_0$ warrants that $t_1<t_0$.
\qed
\mysection{Dissipation in (\ref{d3.1}) implied by uniform smallness of $c$}\label{sect5}
Our next goal is to make sure that as soon as $c$ becomes suitably small, solutions to (\ref{d3.1})
enjoy further regularity properties. 
This will arise as a consequence of the following lemma which provides a weak formulation of a corresponding
differential inequality which can formally be obtained on computing 
\bas
	\frac{d}{dt} \io \psi(n)\rho(c)
\eas
for certain convex functions $\psi$ and $\rho$, 
and which at this heuristic level can be seen to yield an entropy-type
inequality under the structural assumption (\ref{50.200}) below which is satisfied e.g.~for the couple
\be{psi_rho}
	\psi(s):=s^p, \quad s\ge 0, 
	\qquad \mbox{and} \qquad
	\rho(\sigma):=\frac{1}{(2\eta-\sigma)^\theta}, \quad \sigma\in [0,2\eta),
\ee
for any fixed $p>1$ and $\theta\in (0,\frac{p-1}{p+1})$, provided that $\eta=\eta(p,\theta)>0$ is appropriately small
(cf.~also \cite[Lemma 5.1]{win_ARMA} for a similar reasoning for classical solutions in a related problem).
In view of our weak assumptions on regularity of solutions, our rigorous justification of the respective integrated
version, in its overall strategy inspired by the testing procedure presented in \cite{alt_luckhaus}, 
will require certain additional restrictions on the growth of $\psi$ with respect to $n$, which 
at this stage substantially reduce the range of admissible $p$ in (\ref{psi_rho}), but which
in a later step can be removed upon an approximation argument so as to finally allow for any choice of $p$
in Lemma \ref{lem61}.
\begin{lem}\label{lem50}
  Assume (\ref{F1}) and (\ref{F2}). 
  Let $\eta>0$, and suppose that $\psi\in C^2([0,\infty))$ and $\rho\in C^2([0,2\eta))$ are such that
  \be{50.01}
	\psi(s)>0, \quad
	\psi'(s)\ge 0 
	\quad \mbox{and} \quad
	\psi''(s)>0
	\qquad \mbox{for all } s\ge 0
  \ee
  with
  \be{50.1}
	\limsup_{s\to\infty} s^\frac{1}{5} \psi''(s)<\infty,
  \ee
  that
  \be{50.2}
	\rho(\sigma)>0, \quad
	\rho'(\sigma)\ge 0 
	\quad \mbox{and} \quad
	\rho''(\sigma)>0
	\qquad \mbox{for all } \sigma\in [0,2\eta), 
  \ee
  and that
  \be{50.200}
	4\psi'^2(s) \rho'^2(\sigma)
	+ \chi_0^2 s^2 \psi''^2(s) \rho^2(\sigma)
	\le 2 \psi(s)\psi''(s) \rho(\sigma)\rho''(\sigma)
	\qquad \mbox{for all $s\ge 0$ and $\sigma \in [0,\eta]$,}
  \ee
  where $\chi_0:=\|\chi\|_{L^\infty((0,\eta))}$.\\
  Then if $\tz \ge 0$ and 
  $\tu\in L^\infty_{loc}((\tz,\infty);L^2_\sigma(\Omega)) \cap L^2_{loc}((\tz,\infty);W_0^{1,2}(\Omega))$,
  given any strong solution $(n,c)$ of (\ref{d3.1}) in $\Omega\times (\tz,\infty)$ with the additional property that
  \be{50.3}
	\|c\|_{L^\infty(\Omega\times (\tz,\infty))} \le \eta,
  \ee
  one can find a null set $N\subset (\tz,\infty)$ such that
  \bea{50.4}
	& & \hspace*{-16mm}
	\io \psi(n(\cdot,t)) \rho(c(\cdot,t))
	+ \frac{1}{2} \int_{t_0}^t \io \psi''(n) \rho(c) |\nabla n|^2
	\le \io \psi(\cdot,t_0)) \rho(\cdot,t_0)) 
	\quad \mbox{for each $t_0\in (\tz,\infty)\setminus N$} \nn\\[2mm]
	& & \hspace*{102mm}
	\mbox{and all $t\in (t_0,\infty)\setminus N$.}
  \eea
\end{lem}
\proof
  In order to collect some regularity properties needed in the course of our testing procedure, we first observe that as a
  consequence of (\ref{50.1}), we can find $C_1>0$ such that
  \be{50.5}
	\psi''(s) \le C_1(s+1)^{-\frac{1}{5}}
	\qquad \mbox{for all } s\ge 0,
  \ee
  whence there exist $C_2>0$ and $C_3>0$ such that
  \be{50.6}
	\psi'(s) \le C_2 (s+1)^\frac{4}{5}
	\quad \mbox{and} \quad
	\psi(s) \le C_3(s+1)^\frac{9}{5}
	\qquad \mbox{for all } s\ge 0.
  \ee
  Since our hypotheses imply that $\rho\in C^2([0,\eta])$, we can moreover fix positive constants $C_4, C_5$ and $C_6$
  such that
  \be{50.7}
	\rho(\sigma)\le C_4, 
	\quad
	\rho'(\sigma)\le C_5
	\quad \mbox{and} \quad
	\rho''(\sigma) \le C_6
	\qquad \mbox{for all } \sigma\in [0,\eta].
  \ee
  We claim that these inequalities ensure that if $(n,c)$ has the assumed strong solution property and additionally satisfies
  (\ref{50.3}), then
  \be{50.81}
	\psi'(n)\rho(c) \in L^5_{loc}(\bar\Omega\times (\tz,\infty)) \cap L^2_{loc}((\tz,\infty);W^{1,2}(\Omega))
  \ee
  and
  \bea{50.82}
	& & \hspace*{-0mm}
	\Big(\psi(n_{-h})\rho'(c)\Big)_{h\in (0,1)} \subset 
	L^\frac{20}{9}_{loc}(\bar\Omega\times (\tz,\infty)) 
	\quad \mbox{with} \quad \nn\\
	& & \psi(n_{-h})\rho'(c) \to \psi(n)\rho'(c)
	\quad \mbox{in } 
	L^\frac{20}{9}_{loc}(\bar\Omega\times (\tz,\infty)) 
	\quad \mbox{as } h\searrow 0
  \eea
  as well as
  \bea{50.8222}
	& & \bigg(\nabla \Big(\psi(n_{-h})\rho'(c)\Big)\bigg)_{h\in (0,1)} \subset 
	L^\frac{4}{3}_{loc}(\bar\Omega\times (\tz,\infty)) 
	\quad \mbox{with} \quad \nn\\
	& & \nabla \Big(\psi(n_{-h})\rho'(c)\Big)
	\to \nabla \Big(\psi(n)\rho'(c)\Big)
	\quad \mbox{in }
	L^\frac{4}{3}_{loc}(\bar\Omega\times (\tz,\infty)) 
	\quad \mbox{as } h\searrow 0,
  \eea
  where for $h\in (0,1)$ we have set
  \bas
	n_{-h}(x,t):=\left\{ \begin{array}{ll}
	n(x,t-h), \qquad & (x,t)\in \Omega\times (\tz+h,\infty), \\[1mm]
	0, \qquad & (x,t)\in\Omega\times (\tz,\tz+h].
	\end{array} \right.
  \eas
  To this end, we first note that 
  since
  by Young's inequality we obtain as a particular consequence of (\ref{50.5})-(\ref{50.7}) and (\ref{50.3})
  that
  \bas
	|\psi'(n)\rho(c)|^5
	&\le& C_2^5 C_4^5 (n+1)^4
	\qquad \mbox{in } \Omega\times (\tz,\infty)
  \eas
  and
  \bas
	\Big|\nabla \Big(\psi'(n)\rho(c)\Big)\Big|^2
	&\le& \Big\{\psi''(n) \rho(c)|\nabla n| + \psi'(n) \rho'(c) |\nabla c| \Big\}^2\\
	&\le& \Big\{ C_1 C_4 |\nabla n| + C_2 C_5(n+1)^\frac{4}{5} |\nabla c| \Big\}^2 \\
	&\le& 2C_1^2 C_4^2 |\nabla n|^2 + 2C_2^2 C_5^2(n+1)^\frac{16}{5} + 2C_2^2 C_5^2 |\nabla c|^4
	\qquad \mbox{in } \Omega\times (\tz,\infty),
  \eas
  (\ref{50.81}) results upon recalling that $(n+1)^4, |\nabla n|^2$ and $|\nabla c|^4$ 
  and hence clearly also $(n+1)^\frac{16}{5}$ belong to
  $L^1_{loc}(\bar\Omega\times (\tz,\infty))$ according to Definition \ref{defi3}.
  Moreover, since (\ref{50.6}) warrants that 
  \bas
	\psi^\frac{20}{9}(n) \le C_3^\frac{20}{9} (n+1)^4
	\qquad \mbox{in } \Omega\times (\tz,\infty),
  \eas
  it follows from the fact that $n\in L^4_{loc}(\bar\Omega\times (\tz,\infty))$ that
  \be{50.822}
	\psi(n_{-h}) \to \psi(n)
	\quad \mbox{in } L^\frac{20}{9}_{loc}(\bar\Omega\times (\tz,\infty))
	\quad \mbox{as } h\searrow 0,
  \ee
  and that hence in particular (\ref{50.82}) holds
  due to the fact that $\rho'(c)$ is bounded by (\ref{50.7}).
  Apart from this, (\ref{50.822}) also implies that in
  \be{50.833}
	\nabla \Big(\psi(n_{-h})\rho'(c)\Big)
	= \psi'(n_{-h}) \rho'(c)\nabla n_{-h}
	+ \psi(n_{-h})\rho''(c)\nabla c,
  \ee
  we have
  \be{50.84}
	\psi(n_{-h})\rho''(c)\nabla c \to \psi(n)\rho''(c)\nabla c
	\quad \mbox{in } L^\frac{4}{3}_{loc}(\bar\Omega\times (\tz,\infty))
	\quad \mbox{as } h\searrow 0,
  \ee
  because $\rho''(c)\nabla c \in L^4_{loc}(\bar\Omega\times (0,\infty))$ according to (\ref{50.7}) and the requirement
  $\nabla c \in L^4_{loc}(\bar\Omega\times (\tz,\infty))$ in Definition \ref{defi3}, 
  and because $\frac{9}{20}+\frac{1}{4}=\frac{7}{10}<\frac{3}{4}$.
  Since (\ref{50.6}) combined with Young's inequality shows that
  \bas
	|\psi'(n)\nabla n|^\frac{4}{3}
	\le C_2^\frac{4}{3} (n+1)^\frac{16}{15} |\nabla n|^\frac{4}{3}
	\le C_2^\frac{4}{3} (n+1)^\frac{16}{5} + C_2^\frac{4}{3} |\nabla n|^2
	\qquad \mbox{in } \Omega\times (\tz,\infty),
  \eas
  it follows again from the inclusions $n\in L^4_{loc}(\bar\Omega\times (\tz,\infty))$ and 
  $\nabla n\in L^2_{loc}(\bar\Omega\times (\tz,\infty))$
  that $\psi'(n_{-h})\nabla n_{-h} \to \psi'(n)\nabla n$
  in $L^\frac{4}{3}_{loc}(\bar\Omega\times (\tz,\infty))$ as $h\searrow 0$.
  Once more using the boundedness of $\rho'(c)$, we thus obtain that
  \bas
	\psi'(n_{-h})\rho'(c)\nabla n_{-h} \to \psi'(n)\rho'(c)\nabla n
	\quad \mbox{in } L^\frac{4}{3}_{loc}(\bar\Omega\times (\tz,\infty))
	\quad \mbox{as } h\searrow 0,
  \eas
  which together with (\ref{50.833}) and (\ref{50.84}) proves (\ref{50.8222}).\abs
  Since (\ref{50.6}) and (\ref{50.7}) furthermore ensure that $\psi(n)\rho(c) \le C_3 C_4 (n+1)^\frac{9}{5}$
  in $\Omega\times (\tz,\infty)$, it is evident that $\psi(n)\rho(c)$ belongs to 
  $L^1_{loc}(\bar\Omega\times (\tz,\infty))$, so that we can pick a null set $N\subset (\tz,\infty)$ such that
  $(\tz,\infty)\setminus N$ exclusively contains Lebesgue points of $(\tz,\infty)\ni t \mapsto 
  \io \psi(n(\cdot,t))\rho(c(\cdot,t))$. 
  We now fix $t_0\in (\tz,\infty)\setminus N$ and $t_1\in (t_0,\infty)\setminus N$, let	
  $\zeta_\delta$ be as given by (\ref{zeta2}), and define
  \bas
	\phi(x,t):=\zeta_\delta(t) \cdot S_h [\psi'(n)\rho(c)] (x,t),
	\qquad x\in\Omega, \ t>\tz,
  \eas
  for $\delta\in (0,\delta_0)$ and $h\in (0,\delta_0-\delta)$ with $\delta_0:=\min\{1,\frac{t_0-\tz}{2}\}$, where
  the averaging operator $S_h$ is as introduced in (\ref{S}). 
  Then $\phi$ has compact support in $\bar\Omega\times (\tz,t_1+\delta]$, and since evidently 
  $\nabla S_h[\psi'(n)\rho(c)]=S_h[\nabla(\psi'(n)\rho(c))]$, it follows from (\ref{50.81}) that
  $\phi\in L^5(\bar\Omega\times (\tz,\infty)) \cap L^2((\tz,\infty);W^{1,2}(\Omega))$. 
  Computing

  \bas
	\phi_t(x,t)=\zeta_\delta'(t) \cdot S_h[\psi'(n)\rho(c)](x,t)
	+ \zeta_\delta(t)\cdot \frac{\psi'(n(x,t+h))\rho(c(x,t+h))-\psi'(n(x,t))\rho(c(x,t))}{h}
  \eas
  for a.e.~$x\in\Omega$ and $t>\tz$,

  from (\ref{50.81}) we furthermore see that $\phi_t\in L^5(\Omega\times (\tz,\infty)) \subset 
  L^\frac{4}{3}(\Omega\times (\tz,\infty))$, so that Lemma \ref{lem51} ii) guarantees that we may
  use $\phi$ in (\ref{d3.3}) to gain the identity
  \bea{50.10}
	& & \hspace*{-20mm}
	I_1(\delta,h)+I_2(\delta,h)+I_3(\delta,h) \nn\\[2mm]
	&:=& \frac{1}{\delta} \int_{t_1}^{t_1+\delta} \io n(x,t) S_h[\psi'(n)\rho(c)](x,t) dxdt
	- \frac{1}{\delta} \int_{t_0-\delta}^{t_0} \io n(x,t) S_h[\psi'(n)\rho(c)](x,t) dxdt \nn\\
	& & - \frac{1}{h} \int_{\tz}^{t_1+1} \io \zeta_\delta(t) n(x,t) 
	\Big\{ \psi'(n(x,t+h))\rho(c(x,t+h)) - \psi'(n(x,t))\rho(c(x,t)) \Big\} dxdt \nn\\
	&=& - \int_{\tz}^{t_1+1} \io \zeta_\delta(t) \nabla n(x,t) 
	\cdot S_h\Big[\nabla\Big(\psi'(n)\rho(c)\Big)\Big](x,t) dxdt \nn\\
	& & + \int_{\tz}^{t_1+1} \io \zeta_\delta(t) n(x,t) F'(n(x,t)) \chi(c(x,t)) \nabla c(x,t) 
	\cdot S_h\Big[\nabla\Big(\psi'(n)\rho(c)\Big)\Big](x,t) dxdt \nn\\
	& & - \int_{\tz}^{t_1+1} \io \zeta_\delta(t) \Big( \tu(x,t) \cdot \nabla n(x,t) \Big)
	S_h[\psi'(n)\rho(c)](x,t) dxdt \nn\\[2mm]
	&=:& I_4(\delta,h)+I_5(\delta,h)+I_6(\delta,h)
  \eea
  for all $\delta\in (0,\delta_0)$ and $h\in (0,\delta_0-\delta)$.
  Here thanks to (\ref{50.81}) and Lemma \ref{lem60} we have
  \bas	
	S_h\Big[\nabla\Big(\psi'(n)\rho(c)\Big)\Big] \to \nabla \Big(\psi'(n)\rho(c)\Big)
	\quad \mbox{in } L^2_{loc}(\bar\Omega\times (\tz,\infty))
	\qquad \mbox{as } h\searrow 0,
  \eas
  so that since both $\nabla n$ and $nF'(n)\chi(c)\nabla c$ belong to 
  $L^2_{loc}(\bar\Omega\times (\tz,\infty))$
  by Definition \ref{defi3} and Lemma \ref{lem51}, we obtain
  \be{50.12}
	I_4(\delta,h)\to -\int_{\tz}^{t_1+1} \io \zeta_\delta(t) \nabla n \cdot \nabla \Big(\psi'(n)\rho(c)\Big)
	\qquad \mbox{as } h\searrow 0
  \ee
  and
  \be{50.13}
	I_5(\delta,h) \to \int_{\tz}^{t_1+1} \io \zeta_\delta(t) nF'(n)\chi(c) \nabla c\cdot \nabla \Big(\psi'(n)\rho(c)\Big)
	\qquad \mbox{as } h\searrow 0.
  \ee
  We next note that (\ref{50.81}) in light of Lemma \ref{lem60} also warrants that 
  \bas	
	S_h[\psi'(n)\rho(c)] \to \psi'(n)\rho(c)
	\qquad \mbox{in } L^5_{loc}(\bar\Omega\times (\tz,\infty))
	\qquad \mbox{as } h\searrow 0,
  \eas
  which entails that
  \be{50.14}
	I_6(\delta,h) \to -\int_{\tz}^{t_1+1} \io \zeta_\delta(t) \psi'(n)\rho(c) (\tu\cdot\nabla n)
	\qquad \mbox{as } h\searrow 0,
  \ee
  because $\tu \cdot \nabla n \in L^\frac{5}{4}_{loc}(\bar\Omega\times (\tz,\infty))$ according to Lemma \ref{lem51}, and that
  \bea{50.15}
	I_1(\delta,h)+I_2(\delta,h)
	&\to& \frac{1}{\delta} \int_{t_1}^{t_1+\delta} \io n\psi'(n)\rho(c)
	- \frac{1}{\delta} \int_{t_0-\delta}^{t_0} \io n\psi'(n)\rho(c) \nn\\[2mm]
	&=:& I_{12}^\infty(\delta)
	\qquad \mbox{as } h\searrow 0,
  \eea
  for clearly also $n$ lies in $L^\frac{5}{4}_{loc}(\bar\Omega\times (\tz,\infty))$.\\
  As for the remaining term $I_3(\delta,h)$ in (\ref{50.10}), we follow a well-known argument (\cite{alt_luckhaus})
  in using the convexity of $\psi$ to firstly obtain the pointwise estimate
  \bas
	\psi(n(x,t)) - \psi(n(x,t-h)) \le \psi'(n(x,t)) \cdot [n(x,t)-n(x,t-h)]
	\qquad \mbox{for a.e.~$x\in\Omega$ and } t\in (\tz,t_1+1),
  \eas
  which on integration implies that
  \bea{50.16}
	J(\delta,h) 
	&:=& \frac{1}{h} \int_{\tz}^{t_1+1} \io \zeta_\delta(t) \cdot \Big\{\psi(n(x,t))-\psi(n(x,t-h)) \Big\}
	\cdot \rho(c(x,t)) dxdt \nn\\
	&\le& \frac{1}{h} \int_{\tz}^{t_1+1} \io \zeta_\delta(t) n(x,t)\psi'(n(x,t)) \rho(c(x,t)) dxdt \nn\\
	& & - \frac{1}{h} \int_{\tz}^{t_1+1} \io \zeta_\delta(t) n(x,t-h) \psi'(n(x,t)) \rho(c(x,t)) dxdt \nn\\
	&=& \Bigg\{ \frac{1}{h} \int_{\tz}^{t_1+1} \io \zeta_\delta(t) n(x,t)\psi'(n(x,t)) \rho(c(x,t)) dxdt \nn\\
	& & - \frac{1}{h} \int_{\tz}^{t_1+1} \io \zeta_\delta(t-h) n(x,t-h) \psi'(n(x,t)) \rho(c(x,t)) dxdt \Bigg\} \nn\\
	& & + \Bigg\{ \frac{1}{h} \int_{\tz}^{t_1+1} \io \zeta_\delta(t-h)n(x,t-h) \psi'(n(x,t)) \rho(c(x,t)) dxdt \nn\\
	& & - \frac{1}{h} \int_{\tz}^{t_1+1} \io \zeta_\delta(t) n(x,t-h) \psi'(n(x,t)) \rho(c(x,t)) dxdt \Bigg\} \nn\\[2mm]
	&=:& J_1(\delta,h)+J_2(\delta,h)
  \eea
  for all $\delta\in (0,\delta_0)$ and $h\in (0,\delta_0-\delta)$.
  Here the substitution $s=t-h$ reveals that
  \be{50.17}
	J_1(\delta,h)=I_3(\delta,h),
  \ee
  whereas with $I_{12}^\infty(\delta)$ as in (\ref{50.15}) we have
  \bea{50.18}
	J_2(\delta,h)
	&=& - \int_{\tz}^{t_1+1} \io \frac{\zeta_\delta(t-h)-\zeta_\delta(t)}{-h} \cdot n(x,t-h) \psi'(n(x,t)) \rho(c(x,t))
	dxdt \nn\\
	&\to& - \int_{\tz}^{t_1+1} \io \zeta_\delta'(t) n(x,t) \psi'(n(x,t)) \rho(c(x,t)) dxdt \nn\\[2mm]
	&=& I_{12}^\infty(\delta)
	\qquad \mbox{as } h\searrow 0
  \eea
  because 
  clearly $\frac{\zeta_\delta(\cdot-h)-\zeta_\delta}{-h} \wsto \zeta_\delta'$ in $L^\infty((\tz,t_1+1))$
  and 
  $n_{-h} \to n$ in $L^4_{loc}(\bar\Omega\times (\tz,\infty))$ as $h\searrow 0$, and once more because of (\ref{50.81}).\\
  On the other hand, again by substitution we can rewrite the expression on the left-hand side of (\ref{50.16}) 
  according to
  \bea{50.19}
	J(\delta,h)
	&=& \frac{1}{h} \int_{\tz}^{t_1+1} \io \zeta_\delta(t) \psi(n(x,t)) \rho(c(x,t)) dxdt \nn\\
	& & - \frac{1}{h} \int_{\tz}^{t_1+1} \io \zeta_\delta(t+h) \psi(n(x,t)) \rho(c(x,t+h)) dxdt \nn\\
	&=& - \int_{\tz}^{t_1+1} \io \frac{\zeta_\delta(t+h)-\zeta_\delta(t)}{h} \cdot \psi(n(x,t)) \rho(c(x,t+h)) dxdt \nn\\
	& & - \Bigg\{ \frac{1}{h} \int_{\tz}^{t_1+1} \io \zeta_\delta(t) \psi(n(x,t)) \rho(c(x,t+h)) dxdt \nn\\
	& & - \frac{1}{h} \int_{\tz}^{t_1+1} \io \zeta_\delta(t) \psi(n(x,t)) \rho(c(x,t)) dxdt \Bigg\} \nn\\[2mm]
	&=:& J_3(\delta,h)+J_4(\delta,h)
  \eea
  for $\delta\in (0,\delta_0)$ and $h\in (0,\delta_0-\delta)$, where arguing as above we see that
  \bas
	J_3(\delta,h) 
	&\to& - \int_{\tz}^{t_1+1} \io \zeta_\delta'(t) \psi(n(x,t))\rho(c(x,t)) dxdt \\
	&=& - \frac{1}{\delta} \int_{t_0-\delta}^{t_0} \io \psi(n)\rho(c)
	+ \frac{1}{\delta} \int_{t_1}^{t_1+\delta} \io \psi(n)\rho(c)
	\qquad \mbox{as } h\searrow 0.
  \eas
  In summary, from (\ref{50.10}) and (\ref{50.12})-(\ref{50.19}) we thus infer that for all $\delta\in (0,\delta_0)$,
  \bea{50.20}
	& & \hspace*{-40mm}
	- \int_{\tz}^{t_1+1} \io \zeta_\delta(t) \nabla n \cdot \nabla \Big(\psi'(n)\rho(c)\Big) 
	+ \int_{\tz}^{t_1+1} \io \zeta_\delta(t) nF'(n) \chi(c) \nabla c \cdot \nabla \Big(\psi'(n)\rho(c)\Big) \nn\\
	& & - \int_{\tz}^{t_1+1} \io \zeta_\delta(t) \psi'(n)\rho(c) (\tu \cdot \nabla n) \nn\\
	&=& \liminf_{h\searrow 0} \Big\{I_1(\delta,h)+I_2(\delta,h)+I_3(\delta,h)\Big\} \nn\\
	&=& I_{12}^\infty(\delta) + \liminf_{h\searrow 0} J_1(\delta,h) \nn\\
	&\ge& I_{12}^\infty(\delta) + \liminf_{h\searrow 0} \Big\{J(\delta,h)-J_2(\delta,h)\Big\} \nn\\
	&=& I_{12}^\infty(\delta) + \liminf_{h\searrow 0} J(\delta,h) - I_{12}^\infty(\delta) \nn\\
	&=& \liminf_{h\searrow 0} \Big\{ J_3(\delta,h)+J_4(\delta,h)\Big\} \nn\\
	&=& - \frac{1}{\delta} \int_{t_0-\delta}^{t_0} \io \psi(n)\rho(c)
	+ \frac{1}{\delta} \int_{t_1}^{t_1+\delta} \io \psi(n)\rho(c)
	+ \liminf_{h\searrow 0} J_4(\delta,h)
  \eea
  Now an estimate for $J_4(\delta,h)$ can be obtained by pursuing a variant of the above strategy:
  First, by convexity of $\rho$ we see that
  \bas
	\rho(c(x,t+h))-\rho(c(x,t))
	\le \rho'(c(x,t+h)) \cdot [c(x,t+h)-c(x,t)]
	\qquad \mbox{for a.e.~$x\in\Omega$ and } t\in (\tz,t_1+1),
  \eas
  and that hence
  \bea{50.210}
	J_4(\delta,h)
	&=& - \frac{1}{h} \int_{\tz}^{t_1+1} \io \zeta_\delta(t) \psi(n(x,t))
	\cdot \Big\{ \rho(c(x,t+h))-\rho(c(x,t))\Big\} dxdt \nn\\
	&\ge& - \frac{1}{h} \int_{\tz}^{t_1+1} \io \zeta_\delta(t) \psi(n(x,t) \rho'(c(x,t+h)) \cdot
	\Big\{ c(x,t+h)-c(x,t) \Big\} dxdt \nn\\
	&=& - \frac{1}{h} \int_{\tz}^{t_1+1} \io \zeta_\delta(t) c(x,t+h) \psi(n(x,t) \rho'(c(x,t+h)) dxdt \nn\\
	& & + \frac{1}{h} \int_{\tz}^{t_1+1} \io \zeta_\delta(t-h)c(x,t) \psi(n(x,t)) \rho'(c(x,t+h)) dxdt \nn\\
	& & + \frac{1}{h} \int_{\tz}^{t_1+1} \io \zeta_\delta(t) c(x,t) \psi(n(x,t)) \rho'(c(x,t+h)) dxdt \nn\\
	& & - \frac{1}{h} \int_{\tz}^{t_1+1} \io \zeta_\delta(t-h) c(x,t) \psi(n(x,t)) \rho'(c(x,t+h)) dxdt \nn\\
	&=& \frac{1}{h} \int_{\tz}^{t_1+1} \io \zeta_\delta(t-h) c(x,t) \cdot 
	\Big\{ \psi(n(x,t))\rho'(c(x,t+h)) - \psi(n(x,t-h)) \rho'(c(x,t))\Big\} dxdt \nn\\
	& & + \int_{\tz}^{t_1+1} \io \frac{\zeta_\delta(t-h)-\zeta_\delta(t)}{-h} \cdot
	c(x,t) \psi(n(x,t)) \rho'(c(x,t+h)) dxdt \nn\\[2mm]
	&=:& J_{41}(\delta,h) + J_{42}(\delta,h)
	\qquad \mbox{for all $\delta\in (0,\delta_0)$ and } h\in (0,\delta_0-\delta),
  \eea
  where we have substituted $t$ by $t-h$ in one of the integrals making up $J_{41}(\delta,h)$.
  Here, arguing as above we infer that
  \bea{50.211}
	J_{42}(\delta,h)
	&\to& \int_{\tz}^{t_1+1} \io \zeta_\delta'(t) c(x,t) \psi(n(x,t)) \rho'(c(x,t)) dxdt \nn\\
	&=& - \frac{1}{\delta} \int_{t_1}^{t_1+\delta} \io c\psi(n)\rho'(c)
	+ \frac{1}{\delta} \int_{t_0-\delta}^{t_0} c\psi(n)\rho'(c) \nn\\[2mm]
	&=:& \tilde L_{21}(\delta)
	\qquad \mbox{as } h\searrow 0,
  \eea
  and in order to gain appropriate information on $J_{41}(\delta,h)$ through the PDE satisfied by $c$, we note that since
  $\nabla c\in L^4_{loc}(\bar\Omega\times (\tz,\infty))$ by assumption and
  $(\psi(n_{-h})\rho'(c))_{h\in (0,1)} \subset L^\frac{20}{9}_{loc}(\bar\Omega\times (\tz,\infty))$
  as well as
  $(\nabla(\psi(n_{-h})\rho'(c)))_{h\in (0,1)} \subset L^\frac{4}{3}_{loc}(\bar\Omega\times (\tz,\infty))$
  by (\ref{50.82}) and (\ref{50.8222}), for all $\delta\in (0,\delta_0)$ and
  $h\in (0,\delta_0-\delta)$ we may invoke Lemma \ref{lem51} iii) to use
  \bas
	\tilde \phi(x,t):=\zeta_\delta(t-h) \cdot S_h [\psi(n_{-h})\rho'(c)](x,t),
	\qquad x\in \Omega, \ t>\tz,
  \eas
  as a test function in (\ref{d3.4}).
  We thereby obtain the identity
  \bea{50.21}
	& & \hspace*{-10mm}
	\tilde I_1(\delta,h)+\tilde I_2(\delta,h)+\tilde I_3(\delta,h) \nn\\[2mm]
	&:=& \frac{1}{\delta} \int_{t_1+h}^{t_1+\delta+h} \io c(x,t) S_h[\psi(n_{-h}) \rho'(c)](x,t) dxdt \nn\\
	& & - \frac{1}{\delta} \int_{t_0-\delta+h}^{t_0+h} \io c(x,t) S_h[\psi(n_{-h}) \rho'(c)](x,t) dxdt \nn\\
	& & - \frac{1}{h} \int_{\tz}^{t_1+1} \io \zeta_\delta(t-h) c(x,t) \cdot 
	\Big\{ \psi(n_{-h}(x,t+h))\rho'(c(x,t+h)) - \psi(n_{-h}(x,t)) \rho'(c(x,t))\Big\} dxdt \nn\\
	&=& \int_{\tz}^{t_1+1} \io \zeta_\delta(t-h)\nabla c(x,t) \cdot 
	S_h\Big[ \nabla \Big(\psi(n_{-h})\rho'(c)\Big) \Big] (x,t) dxdt \nn\\
	& & - \int_{\tz}^{t_1+1} \io \zeta_\delta(t-h) F(n(x,t)) f(c(x,t)) S_h[\psi(n_{-h})\rho'(c)](x,t) dxdt \nn\\
	& & - \int_{\tz}^{t_1+1} \io \zeta_\delta(t-h) (\tu(x,t)\cdot \nabla c(x,t)) S_h[\psi(n_{-h})\rho'(c)](x,t) dxdt 
	\nn\\[2mm]
	&=:& \tilde I_4(\delta,h) + \tilde I_5(\delta,h) + \tilde I_6(\delta,h)
	\qquad \mbox{for all $\delta\in (0,\delta_0)$ and } h\in (0,\delta_0-\delta),
  \eea
  where evidently
  \be{50.212}
	\tilde I_3(\delta,h)=-J_{41}(\delta,h)
	\qquad \mbox{for all $\delta\in (0,\delta_0)$ and } h\in (0,\delta_0-\delta).
  \ee
  We now use that (\ref{50.82}) and (\ref{50.8222}) along with Lemma \ref{lem60} guarantee that
  \be{50.23}
	S_h[\psi(n_{-h})\rho'(c)] \to \psi(n)\rho'(c)
	\qquad \mbox{in } L^\frac{20}{9}_{loc}(\bar\Omega\times (\tz,\infty))
  \ee
  and 
  \be{50.22}
	S_h\Big[ \nabla \Big(\psi(n_{-h})\rho'(c)\Big) \Big] 
	\to \nabla \Big(\psi(n) \rho'(c)\Big)
	\qquad \mbox{in } L^\frac{4}{3}_{loc}(\bar\Omega\times (\tz,\infty))
  \ee
  as $h\searrow 0$, and that hence
  \be{50.24}
	\tilde I_4(\delta,h) \to - \int_{\tz}^{t_1+1} \io \zeta_\delta(t) \nabla c \cdot \nabla \Big(\psi(n)\rho'(c)\Big)
  \ee
  and
  \be{50.25}
	\tilde I_5(\delta,h) \to - \int_{\tz}^{t_1+1} \io \zeta_\delta(t) F(n)f(c) \psi(n)\rho'(c)
  \ee
  as well as
  \be{50.26}
	\tilde I_6(\delta,h) \to - \int_{\tz}^{t_1+1} \io \zeta_\delta(t) \psi(n)\rho'(c) (\tu\cdot\nabla c)
  \ee
  as $h\searrow 0$, 
  because $\nabla c\in L^4_{loc}(\bar\Omega\times (\tz,\infty))$,
  $F(n)f(c) \in L^4_{loc}(\bar\Omega\times (\tz,\infty))\subset L^\frac{20}{11}_{loc}(\bar\Omega\times (\tz,\infty))$ and also
  $\tu\cdot\nabla c\in L^\frac{20}{11}_{loc}(\bar\Omega\times (\tz,\infty))$ according to Definition \ref{defi3}
  and Lemma \ref{lem51}.\\
  Next, on the left of (\ref{50.21})
  we may combine (\ref{50.23}) with the 
  fact that the family $(\one_{(t_1+h,t_1+\delta+h)})_{h\in (0,1)}$ of indicator functions satisfies
  $\one_{(t_1+h,t_1+\delta+h)} \wsto \one_{(t_1,t_1+\delta)}$ in $L^\infty(\R)$ as $h\searrow 0$
  to conclude, once more relying on the boundedness of $c$, that
  \bea{50.27}
	\tilde I_1(\delta,h)+\tilde I_2(\delta,h)
	&\to& \frac{1}{\delta} \int_{t_1}^{t_1+\delta} \io c \psi(n)\rho'(c)
	- \frac{1}{\delta} \int_{t_0-\delta}^{t_0} \io c\psi(n)\rho'(c) \nn\\[2mm]
	&=:& - \tilde I_{12}^\infty(\delta)
	\qquad \mbox{as } h\searrow 0.
  \eea
  Collecting (\ref{50.212}) and (\ref{50.24})-(\ref{50.27}), from (\ref{50.210}), (\ref{50.211}) and (\ref{50.21})
  we all in all infer that
  \bas
	\liminf_{h\searrow 0} J_4(\delta,h)
	&\ge& \liminf_{h\searrow 0} \Big\{ J_{41}(\delta,h)+J_{42}(\delta,h)\Big\} \\
	&=& \liminf_{h\searrow 0} J_{41}(\delta,h) + \tilde I_{12}^\infty(\delta) \\
	&=& \liminf_{h\searrow 0} \Big\{-\tilde I_3(\delta,h)\Big\} + \tilde I_{12}^\infty(\delta) \\
	&=& \lim_{h\searrow 0} \Big\{\tilde I_1(\delta,h)+\tilde I_2(\delta,h)\Big\} 
	+ \lim_{h\searrow 0} \Big\{ -\tilde I_4(\delta,h)-\tilde I_5(\delta,h)-\tilde I_6(\delta,h)\Big\} 
	+ \tilde I_{12}^\infty(\delta) \\
	&=& \int_{\tz}^{t_1+1} \io \zeta_\delta(t) \nabla c \cdot \nabla \Big(\psi(n)\rho'(c)\Big) 
	+ \int_{\tz}^{t_1+1} \io \zeta_\delta(t) F(n)f(c) \psi(n) \rho'(c) \\
	& & + \int_{\tz}^{t_1+1} \io \zeta_\delta(t) \psi(n)\rho'(c) (\tu\cdot\nabla c)
	\qquad \mbox{for all } \delta\in (0,\delta_0).
  \eas
  Consequently, (\ref{50.20}) implies that
  \bea{50.30}
	& & \hspace*{-20mm}
	\frac{1}{\delta} \int_{t_1}^{t_1+\delta} \io \psi(n) \rho(c)
	- \frac{1}{\delta} \int_{t_0-\delta}^{t_0} \io \psi(n) \rho(c) \nn\\
	&\le& - \int_{\tz}^{t_1+1} \io \zeta_\delta(t) \nabla n \cdot \nabla\Big(\psi'(n)\rho(c)\Big) 
	+ \int_{\tz}^{t_1+1} \io \zeta_\delta(t) nF'(n) \chi(c) \nabla c \cdot \nabla \Big(\psi'(n)\rho(c)\Big) \nn\\
	& & - \int_{\tz}^{t_1+1} \io \zeta_\delta(t) \psi'(n)\rho(c) (\tu\cdot\nabla n) 
	- \int_{\tz}^{t_1+1} \io \zeta_\delta(t) \nabla c \cdot \nabla \Big( \psi(n)\rho'(c)\Big) \nn\\
	& & - \int_{\tz}^{t_1+1} \io \zeta_\delta(t) F(n)f(c)  \psi(n) \rho'(c)
	- \int_{\tz}^{t_1+1} \io \zeta_\delta(t) \psi(n)\rho'(c) (\tu\cdot\nabla c)
  \eea
  for all $\delta\in (0,\delta_0)$. Here since $\tu$ is solenoidal, two integrations by parts show that
  \bas
	& & \hspace*{-30mm}
	- \int_{\tz}^{t_1+1} \io \zeta_\delta(t) \psi'(n)\rho(c) (\tu\cdot\nabla n) 
	- \int_{\tz}^{t_1+1} \io \zeta_\delta(t) \psi(n)\rho'(c) (\tu\cdot\nabla c) \\
	&=&  - \int_{\tz}^{t_1+1} \io \zeta_\delta(t) \rho(c) \tu \cdot \nabla \psi(n)
	- \int_{\tz}^{t_1+1} \io \zeta_\delta(t) \psi(n) \tu \cdot \nabla \rho(c) \\
	&=& - \int_{\tz}^{t_1+1} \io \zeta_\delta(t) \tu \cdot \nabla \Big(\psi(n)\rho(c)\Big) \\[2mm]
	&=& 0
	\qquad \mbox{for all } \delta\in (0,\delta_0),
  \eas
  whereas clearly
  \bas
	- \int_{\tz}^{t_1+1} \io \zeta_\delta(t) F(n)f(c) \psi(n) \rho'(c) \le 0
	\qquad \mbox{for all } \delta\in (0,\delta_0).
  \eas
  Therefore, from (\ref{50.30}) we infer that
  \bea{50.32}
	& & \hspace*{-15mm}
	\frac{1}{\delta} \int_{t_1}^{t_1+\delta} \io \psi(n) \rho(c)
	- \frac{1}{\delta} \int_{t_0-\delta}^{t_0} \io \psi(n) \rho(c) \nn\\
	&\le& - \int_{\tz}^{t_1+1} \io \zeta_\delta(t) \nabla n \cdot \nabla\Big(\psi'(n)\rho(c)\Big) \nn\\
	& & + \int_{\tz}^{t_1+1} \io \zeta_\delta(t) nF'(n) \chi(c) \nabla c \cdot \nabla \Big(\psi'(n)\rho(c)\Big) \nn\\
	& & - \int_{\tz}^{t_1+1} \io \zeta_\delta(t) \nabla c \cdot \nabla \Big( \psi(n)\rho'(c)\Big) \nn\\
	&=& - \int_{\tz}^{t_1+1} \io \zeta_\delta(t) \psi''(n) \rho(c)|\nabla n|^2
	- \int_{\tz}^{t_1+1} \io \zeta_\delta(t) \psi'(n) \rho'(c) \nabla n\cdot \nabla c \nn\\
	& & + \int_{\tz}^{t_1+1} \io \zeta_\delta(t) nF'(n) \psi''(n) \chi(c)\rho(c) \nabla n\cdot\nabla c
	+ \int_{\tz}^{t_1+1} \io \zeta_\delta(t) nF'(n) \psi'(n) \chi(c) \rho'(c) |\nabla c|^2 \nn\\
	& & - \int_{\tz}^{t_1+1} \io \zeta_\delta(t) \psi'(n) \rho'(c) \nabla n \cdot \nabla c
	- \int_{\tz}^{t_1+1} \io \zeta_\delta(t) \psi(n)\rho''(c) |\nabla c|^2 \nn\\
	&=& - \int_{\tz}^{t_1+1} \io \zeta_\delta(t) \psi''(n) \rho(c)|\nabla n|^2 \nn\\
	& & - \int_{\tz}^{t_1+1} \io \zeta_\delta(t) \cdot \Big\{
	2\psi'(n)\rho'(c) - nF'(n) \psi''(n) \chi(c)\rho(c)\Big\} \nabla n\cdot \nabla c \nn\\
	& & - \int_{\tz}^{t_1+1} \io \zeta_\delta(t) \cdot \Big\{
	\psi(n) \rho''(c) - nF'(n) \psi'(n) \chi(c)\rho'(c)\Big\} |\nabla c|^2
	\qquad \mbox{for all } \delta\in (0,\delta_0).
  \eea
  Here we can estimate the second integral on the right by means of Young's inequality according to
  \bea{50.333}
	& & \hspace*{-20mm}
	- \int_{\tz}^{t_1+1} \io \zeta_\delta(t) \cdot \Big\{
	2\psi'(n)\rho'(c) - nF'(n) \psi''(n) \chi(c)\rho(c)\Big\} \nabla n\cdot \nabla c \nn\\
	&\le& \frac{1}{2} \int_{\tz}^{t_1+1} \io \zeta_\delta(t) \psi''(n) |\nabla n|^2 \rho(c) \nn\\
	& & + \int_{\tz}^{t_1+1} \io \zeta_\delta(t) \cdot 
	\frac{\Big\{2\psi'(n)\rho'(c)-nF'(n)\psi''(n)\chi(c)\rho(c)\Big\}^2}{2\psi''(n)\rho(c)} \cdot |\nabla c|^2 
  \eea
  for all $\delta\in (0,\delta_0)$.
  Now since $F', \psi',\psi'', \chi,\rho$ and $\rho'$ are all nonnegative, and since $F'^2 \le 1$ on $[0,\infty)$
  by (\ref{F1}) and $\chi^2(c) \le \chi_0^2$ a.e.~in $\Omega\times (\tz,\infty)$ by (\ref{50.3}), the hypothesis
  (\ref{50.200}) warrants that
  \bas
	& & \hspace*{-12mm}
	\frac{\Big\{2\psi'(n)\rho'(c)-nF'(n)\psi''(n)\chi(c)\rho(c)\Big\}^2}{2\psi''(n)\rho(c)} 
	- \Big\{\psi(n)\rho''(c) - nF'(n) \psi'(n)\chi(c)\rho'(c) \Big\} \\
	&=& \frac{1}{2\psi''(n)\rho(c)} \cdot 
	\Bigg\{ \Big\{ 4\psi'^2(n)\rho'^2(c) - 4nF'(n)\psi'(n)\psi''(n)\chi(c)\rho(c)\rho'(c)
	+ n^2 F'^2(n) \psi''^2(n) \chi^2(c) \rho^2(c) \Big\} \\
	& & \hspace*{23mm}
	- \Big\{2\psi(n)\psi''(n) \rho(c)\rho''(c) + 2nF'(n) \psi'(n)\psi''(n) \chi(c) \rho(c)\rho'(c)\Big\} \Bigg\} \\[2mm]
	&=& \frac{1}{2\psi''(n)\rho(c)} \cdot 
	\bigg\{ 4\psi'^2(n)\rho'^2(c) + n^2 F'^2(n) \psi''^2(n) \chi^2(c) \rho^2(c) \\
	& & \hspace*{24mm}
	- 2\psi(n)\psi''(n)\rho(c)\rho''(c)
	- 2nF'(n) \psi'(n) \psi''(n) \chi(c)\rho(c)\rho'(c) \bigg\} \\[2mm]
	&\le& \frac{1}{2\psi''(n)\rho(c)} \cdot 
	\bigg\{ 4\psi'^2(n)\rho'^2(c) 
	+ n^2 \psi''^2(n) \chi_0^2 \rho^2(c) 
	- 2\psi(n)\psi''(n)\rho(c)\rho''(c) \bigg\} \\[2mm]
	&\le& 0
	\qquad \mbox{a.e.~in } \Omega\times (\tz,\infty).
  \eas
  Inserting (\ref{50.333}) into (\ref{50.32}), we thus infer that
  \bea{50.34}
	& & \hspace*{-20mm}
	\frac{1}{\delta} \int_{t_1}^{t_1+\delta} \io \psi(n)\rho(c)
	+ \frac{1}{2} \int_{\tz}^{t_1+1} \io \zeta_\delta(t) \psi''(n) \rho(c) |\nabla n|^2 \nn\\
	&\le& \frac{1}{\delta} \int_{t_0-\delta}^{t_0} \io \psi(n)\rho(c)
	\qquad \mbox{for all } \delta\in (0,\delta_0),
  \eea
  where the Lebesgue point properties of $t_0$ and $t_1$ ensure that 
  \bas
	\frac{1}{\delta} \int_{t_0-\delta}^{t_0} \io \psi(n)\rho(c) 
	\to \io \psi(n(\cdot,t_0))\rho(c(\cdot,t_0))
	\quad \mbox{and} \quad
	\frac{1}{\delta} \int_{t_1}^{t_1+\delta} \io \psi(n)\rho(c) 
	\to \io \psi(n(\cdot,t_1))\rho(c(\cdot,t_1))
  \eas
  as $h\searrow 0$. Again since $\psi''$ is nonnegative, using that also $\zeta_\delta \ge 0$ and that $\zeta_\delta\equiv 1$
  in $(t_0,t_1)$ we thereby obtain from (\ref{50.34}) that indeed (\ref{50.4}) is valid.
\qed
\mysection{Estimates implied by Lemma \ref{lem50}}\label{sect6}
Our application of Lemma \ref{lem50} will be prepared by the following statement which partly explains 
the particular approximation of $[0,\infty)\mapsto s^p$ to be pursued in Lemma \ref{lem61}.
It may be worthwhile mentioning here that for given $p\ge 2$
it seems impossible to adjust $q<p$ in such a way that for small $\delta>0$,
the alternative and 
apparently more straightforward choice $\psi_\delta(s):=\frac{s^p}{1+\delta s^q}$ is admissible in Lemma \ref{lem50}.
We therefore employ a certain integrated variant thereof, with its precise form and some of its properties
described as follows.
\begin{lem}\label{lem62}
  For $p\ge 2, q\in (0,p-1)$ and $\delta\in (0,1)$, we let
  \be{62.1}
	\psi_\delta(s):=p\int_0^s \frac{\sigma^{p-1}}{1+\delta \sigma^q} d\sigma,
	\qquad s\ge 0.
  \ee
  Then $\psi_\delta>0, \psi_\delta'>0$ and $\psi_\delta''>0$ on $(0,\infty)$, and for all $s\ge 0$ we have
  \be{62.3}
	\psi_\delta(s)\nearrow s^p
	\quad \mbox{and} \quad
	\psi_\delta''(s) \to p(p-1) s^{p-2}
	\qquad \mbox{as } \delta\searrow 0.
  \ee
  Moreover,
  \be{62.4}
	\frac{\psi_\delta'^2(s)}{\psi_\delta(s) \psi_\delta''(s)} \le \frac{p}{p-q-1}
	\qquad \mbox{for all } s>0
  \ee
  and 
  \be{62.5}
	\frac{s^2 \psi_\delta''(s)}{\psi_\delta(s)} \le p(p-1)
	\qquad \mbox{for all } s>0,
  \ee
  and for each fixed $\delta\in (0,1)$ we have
  \be{62.6}
	s^{-p+q+2} \psi_\delta''(s) \to \frac{p(p-q-1)}{\delta}
	\qquad \mbox{as } s\to\infty.
  \ee
\end{lem}
\proof
  We first note that
  \bas
	(1+\delta s^q) \psi_\delta(s)
	= p\int_0^s \frac{1+\delta s^q}{1+\delta \sigma^q} \cdot \sigma^{p-1} d\sigma
	\ge p\int_0^s \sigma^{p-1} d\sigma 
	= s^p
	\qquad \mbox{for all } s\ge 0,
  \eas
  so that
  \be{62.7}
	\psi_\delta(s) \ge \frac{s^p}{1+\delta s^q}
	\qquad \mbox{for all } s\ge 0.
  \ee
  Next, computing
  \be{62.8}
	\psi_\delta'(s)
	= \frac{ps^{p-1}}{1+\delta s^q}
	\quad \mbox{and} \quad
	\psi_\delta''(s) = p \cdot \frac{(p-1)s^{p-2} + (p-q-1) \delta s^{p+q-2}}{(1+\delta s^q)^2}
  \ee
  for $s\ge 0$, using that $q<p-1$ we see that $\psi_\delta, \psi_\delta'$ and $\psi_\delta''$ indeed 
  are all positive on $(0,\infty)$.\\
  Furthermore, combining (\ref{62.8}) with (\ref{62.7}) shows that
  \bas
	\frac{\psi_\delta'^2(s)}{\psi_\delta(s)\psi_\delta''(s)}
	&\le& \frac{\psi_\delta'^2(s)}{\psi_\delta''(s)} \cdot \frac{1+\delta s^q}{s^p} \\
	&=& \frac{p^2 s^{2p-2}}{(1+\delta s^q)^2}
	\cdot \frac{1}{p} \frac{(1+\delta s^q)^2}{(p-1) s^{p-2} + (p-q-1) \delta s^{p+q-2}}	
	\cdot \frac{1+\delta s^q}{s^p} \\
	&=& p\cdot \frac{1+\delta s^q}{p-1 + (p-q-1)\delta s^q} \\
	&\le& p \cdot \frac{1+\delta s^q}{p-q-1 + (p-q-1) \delta s^q}
	= \frac{p}{p-q-1}
	\qquad \mbox{for all } s>0,
  \eas
  because $q>0$, and that similarly
  \bas
	\frac{s^2 \psi_\delta''(s)}{\psi_\delta(s)}
	&\le& s^2 \psi_\delta''(s) \cdot \frac{1+\delta s^q}{s^p} \\
	&=& s^2 \cdot p\frac{(p-1) s^{p-2} + (p-q-1) \delta s^{p+q-2}}{(1+\delta s^q)^2}
	\cdot \frac{1+\delta s^q}{s^p} \\
	&=& p \cdot \frac{p-1+(p-q-1) \delta s^q}{1+\delta s^q} \\
	&\le& p\cdot \frac{p-1+(p-1)\delta s^q}{1+\delta s^q}
	= p(p-1)
	\qquad \mbox{for all } s>0.
  \eas
  Finally, (\ref{62.6}) and the second statement in (\ref{62.3}) are evident from (\ref{62.8}), whereas the first
  claim in (\ref{62.3}) results from Beppo Levi's theorem.
\qed
Along with the originally intended choice (\ref{psi_rho}) of $\rho$, making use of these functions in Lemma \ref{lem50}
allows us to deduce an entropy-type inequality under a smallness assumption on $c$.
For simplicity in presentation, we confine ourselves to proving the inequality (\ref{61.2}) which
actually slightly differs from an inequality indicating genuine decrease
of the functional $\io n^p$ in that it involves a factor $2$ on its right, but the boundedness and dissipation
properties thereby implied will be sufficient for our purpose.
\begin{lem}\label{lem61}
  For all $p\ge 2$ there exists $\eta>0$ with the following property: 
  Whenever $\tz\ge 0, F$ satisfies (\ref{F1}) and
  $\tu\in L^2_{loc}((\tz,\infty);W_0^{1,2}(\Omega)) \cap L^\infty_{loc}((\tz,\infty);L^2_\sigma(\Omega))$,
  for any strong solution $(n,c)$ of (\ref{d3.1}) in $\Omega\times (\tz,\infty)$ fulfilling
  \be{61.1}
	\|c\|_{L^\infty(\Omega\times (\tz,\infty))} \le \eta
  \ee
  there exists a null set $N(p)\subset (\tz,\infty)$ such that
  \bea{61.2}
	\io n^p(\cdot,t) + \frac{p(p-1)}{2} \int_{t_0}^t \io n^{p-2} |\nabla n|^2
	&\le& 2\io n^p(\cdot,t_0)
	\qquad \mbox{for all $t_0\in (\tz,\infty)\setminus N(p)$} \nn\\[1mm]
	& & \hspace*{28mm} 
	\mbox{and each } t\in (t_0,\infty)\setminus N(p).
  \eea
\end{lem}
\proof
  Given $p\ge 2$, we first choose $\theta\in (0,1)$ small enough such that
  \be{61.3}
	\frac{5p\theta}{\theta+1} \le 1
  \ee
  and then fix a small number $\eta\in (0,1)$ satisfying
  \be{61.4}
	\frac{4p(p-1)\chi_1^2 \eta^2}{\theta(\theta+1)} \le 1,
  \ee
  where $\chi_1:=\|\chi\|_{L^\infty((0,1))}$.
  With any fixed sequence $(\delta_j)_{j\in\N} \subset (0,1)$ satisfying $\delta_j\searrow 0$ as $j\to\infty$,
  we now let $\psi_\delta$ be as defined in Lemma \ref{lem62} with $q:=p-\frac{9}{5}$, that is, we let
  \bas
	\psi_\delta(s):= p\int_0^s \frac{\sigma^{p-1}}{1+\delta \sigma^{p-\frac{9}{5}}} d\sigma
	\qquad \mbox{for } s\ge 0,
  \eas
  and moreover we set
  \be{61.44}
	\rho(\sigma):=\frac{1}{(2\eta-\sigma)^\theta}
	\qquad \mbox{for } \sigma\in [0,2\eta).
  \ee
  Then $\psi_\delta>0, \psi_\delta'>0$ and $\psi_\delta''>0$ on $(0,\infty)$ by Lemma \ref{lem62}, whereas computing
  \be{61.5}
	\rho'(\sigma)=\frac{\theta}{(2\eta-\sigma)^{\theta+1}}
	\quad \mbox{and} \quad
	\rho''(\sigma)=\frac{\theta(\theta+1)}{(2\eta-\sigma)^{\theta+2}}
	\qquad \mbox{for } \sigma\in [0,2\eta)
  \ee
  we see that also $\rho>0, \rho'>0$ and $\rho''>0$ throughout $[0,2\eta)$.\\
  By means of (\ref{62.4}) and (\ref{61.5}), we now estimate
  \bea{61.6}
	\frac{4\psi_\delta'^2(s) \rho'^2(\sigma)}{\psi_\delta(s)\psi_\delta''(s) \rho(\sigma)\rho''(\sigma)}
	&\le& 4 \cdot \frac{5p}{4} \cdot \frac{\rho'^2(\sigma)}{\rho(\sigma)\rho''(\sigma)} \nn\\
	&=& 5p \cdot \frac{\theta^2 (2\eta-\sigma)^{-2\theta-2}}
		{(2\eta-\sigma)^{-\theta} \cdot \theta(\theta+1)(2\eta-\sigma)^{-\theta-2}} \nn\\
	&=& \frac{5p\theta}{\theta+1} \nn\\[2mm]
	&\le& 1
	\qquad \mbox{for all $s>0$ and } \sigma \in [0,2\eta)
  \eea
  thanks to (\ref{61.3}), while (\ref{62.5}) in conjunction with (\ref{61.44}) and (\ref{61.5}) shows that
  \bea{61.7}
	\frac{\chi_1^2 s^2 \psi_\delta''^2(s)\rho^2(\sigma)}{\psi_\delta(s)\psi_\delta''(s) \rho(\sigma)\rho''(\sigma)}
	&=& \chi_1^2 \cdot \frac{s^2 \psi_\delta''(s)}{\psi_\delta(s)} \cdot \frac{\rho(\sigma)}{\rho''(\sigma)} \nn\\
	&\le& \chi_1^2 \cdot p(p-1) \cdot \frac{(2\eta-\sigma)^2}{\theta(\theta+1)} \nn\\
	&\le& \chi_1^2 \cdot p(p-1) \cdot \frac{4\eta^2}{\theta(\theta+1)} \nn\\[2mm]
	&\le& 1
	\qquad \mbox{for all $s>0$ and } \sigma \in [0,2\eta)
  \eea
  according to (\ref{61.4}).
  Since $\eta<1$ entails that $\chi_1 \ge \chi_0:=\|\chi\|_{L^\infty((0,\eta))}$, combining (\ref{61.6}) with (\ref{61.7})
  ensures that
  \bas
	4\psi_\delta'^2(s) \rho'^2(\sigma) + \chi_0^2 s^2 \psi_\delta''^2(s)\rho^2(\sigma)
	\le 2\psi_\delta(s) \psi_\delta''(s) \rho(\sigma)\rho''(\sigma)
	\qquad \mbox{for all $s>0$ and } \sigma \in [0,2\eta).
  \eas
  As furthermore our choice of $q$ guarantees that
  $\limsup_{s\to\infty} s^\frac{1}{5} \psi_\delta''(s)=\frac{4p}{5\delta}$ is finite for each $\delta\in (\delta_j)_{j\in\N}$,
  Lemma \ref{lem50} becomes applicable so as to assert that whenever $\tz\ge 0$ and
  $\tu, n$ and $c$ have the assumed properties, for any $\delta\in (\delta_j)_{j\in\N}$ we can find a null set 
  $N_\delta \subset (\tz,\infty)$ such that
  \bea{61.8}
	& & \hspace*{-20mm}
	\io \psi_\delta(n(\cdot,t)) \rho(c(\cdot,t))
	+ \frac{1}{2} \int_{t_0}^t \io \psi_\delta''(n) \rho(c) |\nabla n|^2
	\le \io \psi_\delta(n(\cdot,t_0)) \rho(c(\cdot,t_0)) \nn\\[2mm]
	& & \qquad \mbox{for all $t_0\in (\tz,\infty)\setminus N_\delta$ and } t\in (t_0,\infty)\setminus N_\delta.
  \eea
  Here as $0\le c \le \eta$ a.e.~in $\Omega\times (\tz,\infty)$, recalling (\ref{61.44}) we can estimate
  \bas
	\frac{1}{(2\eta)^\theta} \le \rho(c) \le \frac{1}{\eta^\theta}
	\qquad \mbox{a.e.~in } \Omega\times (\tz,\infty),
  \eas
  so that since the countable union $\bigcup_{\delta\in (\delta_j)_{j\in\N}} N_\delta$ has measure zero, it follows from
  (\ref{61.8}) that with some null set $N=N(p)\subset (\tz,\infty)$ we have
  \bas
	& & \hspace*{-16mm}
	\frac{1}{(2\eta)^\theta} \io \psi_\delta(n(\cdot,t))
	+ \frac{1}{2} \cdot \frac{1}{(2\eta)^\theta} \int_{t_0}^t \io \psi_\delta''(n) |\nabla n|^2
	\le \frac{1}{\eta^\theta} \io \psi_\delta(n(\cdot,t_0))
	\qquad \mbox{for all $t_0\in (\tz,\infty)\setminus N$} \\[0mm]
	& & \hspace*{108mm}
	\mbox{and } t\in (t_0,\infty)\setminus N
  \eas
  and any $\delta\in (\delta_j)_{j\in\N}$.
  Now since from Lemma \ref{lem62} we know that $\psi_\delta(s)\nearrow s^p$ and $\psi_\delta''(s) \to p(p-1) s^{p-2}$ as 
  $\delta=\delta_j\searrow 0$, we may invoke the monotone convergence theorem and Fatou's lemma to conclude that indeed
  \bas	
	\hspace*{-10mm}
	\io n^p(\cdot,t)
	+ \frac{1}{2} \cdot p(p-1) \int_{t_0}^t \io n^{p-2} |\nabla n|^2
	&\le& \frac{(2\eta)^\theta}{\eta^\theta} \io n^p(\cdot,t_0) \nn\\
	&\le& 2\io n^p(\cdot,t_0)
  \eas
  for all $t\in (\tz,\infty)\setminus N$ and $t\in (t_0,\infty)\setminus N$,
  because $\theta<1$ entails that $2^\theta<2$.
\qed
Indeed, the above lemma entails the following boundedness properties, uniform with respect to 
functions in ${\cal S}_{m,\eta,L,\tz}$, provided that $\eta>0$ is small.
We note that since this result will be applied to finitely many $p$ only, the dependence of $\eta$ and the number
$\tau$ therein on $p$ will actually be irrelevant in the sequel.
\begin{lem}\label{lem64}
  Let $p\ge 2$. Then there exist $\eta>0$ and $\tau>0$ with the property that for all $m>0$ and $L>0$ one can find
  $C(m,L)>0$ such that if $\tz\ge 0$, if $F$ complies with (\ref{F1}), if
  $\tu\in L^2_{loc}((\tz,\infty);W_0^{1,2}(\Omega)) \cap L^\infty_{loc}((\tz,\infty);L^2_\sigma(\Omega))$,
  and if $(n,c)$ is a strong solution of (\ref{d3.1}) in $\Omega\times (\tz,\infty)$ with
  \be{64.1}
	\|c\|_{L^\infty(\Omega\times (\tz,\infty))} \le \eta
  \ee
  as well as
  \be{64.2}
	\|n(\cdot,t)\|_{L^1(\Omega)} \le m 	
	\quad \mbox{for all } t>\tz
	\qquad \mbox{and} \qquad
	\int_{\tz}^{\tz+1} \io \frac{|\nabla n|^2}{n} \le L,
  \ee
  then
  \be{64.3}
	\io n^p(\cdot,t) \le C(m,L)
	\qquad \mbox{for a.e.~} t>\tz+\tau
  \ee
  and
  \be{64.4}
	\int_{\tz+\tau}^\infty \io n^{p-2} |\nabla n|^2 \le C(m,L).
  \ee
\end{lem}
\proof
  We first observe that since $W^{1,2}(\Omega) \hra L^6(\Omega)$, (\ref{64.2}) implies that for some $C_1>0$ we have
  \bas
	\int_{\tz}^{\tz+1} \|n(\cdot,t)\|_{L^3(\Omega)} dt
	&=& \int_{\tz}^{\tz+1} \|n^\frac{1}{2}(\cdot,t)\|_{L^6(\Omega)}^2 dt \\
	&\le& C_1 \int_{\tz}^{\tz+1} \bigg\{ \|\nabla n^\frac{1}{2}(\cdot,t)\|_{L^2(\Omega)}^2 
	+ \|n^\frac{1}{2}(\cdot,t)\|_{L^2(\Omega)}^2 \bigg\} dt \\
	&=& \frac{C_1}{4} \int_{\tz}^{\tz+1} \io \frac{|\nabla n|^2}{n} + C_1 \int_{\tz}^{\tz+1} \io n \\
	&\le& \frac{C_1}{4} L + C_1 m.
  \eas
  In view of a recursive argument, to prove the lemma it is thus sufficient to show that for each $p\ge 2$ there exists
  $\eta>0$ such that for any choice of $m>0$ and $B>0$ we can fix $C_2(m,B)>0$ such that if for some $T_1\ge 0$ and some
  $\tu\in L^2_{loc}((T_1,\infty);W_0^{1,2}(\Omega)) \cap L^\infty_{loc}((T_1,\infty);L^2_\sigma(\Omega))$ we are given
  a strong solution $(n,c)$ of (\ref{d3.1}) in $\Omega\times (T_1,\infty)$ fulfilling
  $\|c\|_{L^\infty(\Omega\times (T_1,\infty))} \le \eta$
  and 
  \be{64.44}
	\|n(\cdot,t)\|_{L^1(\Omega)} \le m 	
	\qquad \mbox{for a.e.~} t>T_1
  \ee
  as well as
  \be{64.5}
	\int_{T_1}^{T_1+1} \|n(\cdot,t)\|_{L^p(\Omega)} dt \le B,
  \ee
  then
  \be{64.6}
	\io n^p(\cdot,t) \le C_2(m,B)
	\qquad \mbox{for a.e.~} t>T_1+1
  \ee
  and
  \be{64.7}
	\int_{T_1+1}^\infty \io n^{p-2} |\nabla n|^2 \le C_2(m,B)
  \ee
  as well as
  \be{64.8}
	\int_t^{t+1} \|n(\cdot,t)\|_{L^{3p}(\Omega)} dt \le C_2(m,B)
	\qquad \mbox{for all } t>T_1+1.
  \ee
  To this end, given any such $p$ we invoke Lemma \ref{lem61} to obtain $\eta>0$ with the properties listed there.
  In particular, since (\ref{64.5}) ensures that ${\rm{essinf}}_{t_0\in (T_1,T_1+1)} \io n^p(\cdot,t_0) \le B^p$,
  applying (\ref{61.2}) to some appropriately chosen $t_0\in (T_1,T_1+1)$ shows that
  \be{64.9}
	\io n^p(\cdot,t) \le 2B^p
	\qquad \mbox{for a.e.~} t>t_0
	\qquad \mbox{and} \qquad
	\int_{t_0}^\infty \io n^{p-2} |\nabla n|^2 \le \frac{4}{p(p-1)} B^p.
  \ee
  As moreover using the H\"older inequality and again the continuity of the embedding $W^{1,2}(\Omega)\hra L^6(\Omega)$
  along with (\ref{64.44}) provides $C_3>0$ such that
  \bas
	\bigg( \int_t^{t+1} \|n(\cdot,t)\|_{L^{3p}(\Omega)} dt \bigg)^p
	&\le& \int_t^{t+1} \|n(\cdot,t)\|_{L^{3p}(\Omega)}^p dt \\
	&=& \int_t^{t+1} \|n^\frac{p}{2}(\cdot,t)\|_{L^6(\Omega)}^2 dt \\
	&\le& C_3 \int_t^{t+1} \bigg\{ \|\nabla n^\frac{p}{2}(\cdot,t)\|_{L^2(\Omega)}^2
	+ \|n^\frac{p}{2}(\cdot,t)\|_{L^\frac{2}{p}(\Omega)}^2 \bigg\} dt \\
	&\le& \frac{p^2 C_3}{4} \int_t^{t+1} \io n^{p-2} |\nabla n|^2 + C_3 m^p
	\qquad \mbox{for all } t>t_0,
  \eas
  the inequalities in (\ref{64.9}) entail (\ref{64.6})-(\ref{64.8}).
\qed
\mysection{Ultimate regularity of eventual energy solutions}\label{sect7}
We now focus on the asymptotic analysis of a given particular eventual energy solution $(n,c,u)$, thus aiming at proving
Theorem \ref{theo_eventual}.
\subsection{The inclusion $(n,c,id) \in \set$}
In order to prepare an appropriate exploitation of the energy inequality (\ref{energy}),
we first assert that the dissipation rate appearing therein 
dominates the energy functional $\F$ in the following sense.
\begin{lem}\label{lem333}
  For all $m>0$, $M>0$ and $\kappa>0$ there exists $C=C(m,M,\kappa)>0$ 
  such that if $n\in L^1(\Omega)$ and $c\in L^\infty(\Omega)$
  are nonnegative with $\io n\le m$ and $\|c\|_{L^\infty(\Omega)} \le M$ as well as
  $n^\frac{1}{2} \in W^{1,2}(\Omega)$ and $c^\frac{1}{4} \in W^{1,4}(\Omega)$, and if moreover
  $u \in W_0^{1,2}(\Omega;\R^3)$, then
  \be{333.1}
	-\frac{|\Omega|}{e} \le \F[n,c,u] \le C \cdot \io \bigg\{ \frac{|\nabla n|^2}{n} + \frac{|\nabla c|^4}{c^3}
	+ |\nabla u|^2 \bigg\}
	+ C.
  \ee
\end{lem}
\proof
  We first note that our assumptions $f(0)=0$, $\chi(0)>0$ and $(\frac{f}{\chi})'>0$ on $[0,\infty)$
  imply that there exists $C_1=C_1(M)>0$ such that
  $\frac{f(s)}{\chi(s)} \ge C_1s$ for all $s\in [0,M]$, which in view of Young's inequality entails that
  \bas
	\frac{1}{2} \io \frac{\chi(c)}{f(c)} |\nabla c|^2
	\le \frac{1}{2C_1} \io \frac{|\nabla c|^2}{c}
	\le \io \frac{|\nabla c|^4}{c^3} 
	+ \frac{1}{16C_1^2} \io c
	\le \io \frac{|\nabla c|^4}{c^3} 
	+ \frac{M|\Omega|}{16C_1^2}.
  \eas
  Next, since $z\ln z \le \frac{3}{2} z^\frac{5}{3}$ for all $z\ge 0$, using the Gagliardo-Nirenberg inequality
  we find $C_2>0$ and $C_3=C_3(m)>0$ fulfilling
  \bas
	\io n\ln n
	\le \io n^\frac{5}{3}
	=\frac{3}{2} \|n^\frac{1}{2}\|_{L^\frac{10}{3}(\Omega)}^\frac{10}{3}
	\le C_2 \|\nabla n^\frac{1}{2}\|_{L^2(\Omega)}^2 \|n^\frac{1}{2}\|_{L^2(\Omega)}^\frac{4}{3}
	+ C_2 \|n^\frac{1}{2}\|_{L^2(\Omega)}^\frac{10}{3}
	\le C_3 \io \frac{|\nabla n|^2}{n} + C_3,
  \eas
  because $\|n^\frac{1}{2}\|_{L^2(\Omega)}^2=\io n\le m$.\\

  As, finally, the Poincar\'e inequality provides $C_4>0$ satisfying 
  $\io |u|^2 \le C_4 \io |\nabla u|^2$, we all in all obtain
  \bas
	\F[n,c,u] \le \max \Big\{1,C_3,C_4 \kappa \Big\} \cdot 
	\io \bigg\{ \frac{|\nabla n|^2}{n} + \frac{|\nabla c|^4}{c^3} + |\nabla u|^2 \bigg\}
	+ \frac{M|\Omega|}{16C_1^2} + C_3
  \eas
  for any such $n, c$ and $u$.
  Along with the fact that $\F[n,c,u]\ge -\frac{|\Omega|}{e}$, valid due to the inequality $z\ln z\ge -\frac{1}{e}$
  for $z>0$, this shows (\ref{333.1}).
\qed
We can thereupon make sure that all the results of the previous sections can actually be applied to such solutions, because
$(n,c,id)$ then belongs to $\set$ for adequately chosen $m,M,L$ and $\tz$.
\begin{lem}\label{lem655}
  Let $(n,c,u)$ be an eventual energy solution of (\ref{0}).
  Then there exist $m>0$, $M>0$, $L>0$ and $\tz\ge 0$ such that with $F(s):=s, s\ge 0$, the triple
  $(n,c,F)$ belongs to $\set$.
  In particular,
  \be{c_decay}
	\|c\|_{L^\infty(\Omega\times (t,\infty))} \to 0
	\qquad \mbox{as } t\to\infty.
  \ee
\end{lem}
\proof
  From Definition \ref{defi_ees} and Lemma \ref{lem_mass} it follows that $\io n(\cdot,t)=m:=\io n_0$
  for a.e.~$t>0$, and that there exist $T\ge 0$, $\kappa>0$, $K>0$ and
  $C_1>0$ such that the regularity properties in (\ref{reg_ees}) as well as (\ref{energy}) hold.
  In particular, we may therefore integrate by parts in the rightmost integrals in (\ref{w1}) and (\ref{w2}) 
  for suitably chosen $\phi$ to see that
  $(n,c)$ is a strong solution of (\ref{d3.1}) in $\Omega\times (\tz,\infty)$ with $\tz:=T+1$, $F(s):=s$ for $s\ge 0$
  and $\tu:=u$, which
  in turn allows us to apply Corollary \ref{cor32} to infer that 
  $\|c(\cdot,t)\|_{L^\infty(\Omega)} \le M:=\|c\|_{L^\infty(\Omega\times (T,T+1))}$ for a.e.~$t>\tz$.\\
  We next observe that writing $y(t):=\F[n,c,u](t)$ and
  $h(t):=\io \big\{ \frac{|\nabla n|^2}{n} + \frac{|\nabla c|^4}{c^3} + |\nabla u|^2 \big\}(\cdot,t)$ for $t>T$,
  by a completion argument we infer from (\ref{energy}) that for any nonnegative $\phi\in W^{1,\infty}((T,\infty))$
  with compact support in $(T,\infty)$ we have
  \be{655.2}
	-\int_T^\infty y(t)\phi'(t) dt 
	+ \frac{1}{K} \int_T^\infty h(t)\phi(t) 
	\le K \int_T^\infty \phi(t)dt,
  \ee 
  because both $y$ and $h$ belong to $L^1_{loc}((T,\infty))$ due to Definition \ref{defi_ees} and Lemma \ref{lem333}.\\
  In order to make sure that this implies boundedness of $y$ on $(T+1,\infty)$, we make use of the estimate provided
  by Lemma \ref{lem333} when applied to $m:=\io n_0$ and $M$ as defined above 
  to infer from (\ref{655.2}) the existence of $C_1>0$ and $C_2>0$ such that
  \bas	
	-\int_T^\infty y(t)\phi'(t) dt 
	+ C_1 \int_T^\infty y(t)\phi(t) 
	\le C_2 \int_T^\infty \phi(t)dt
  \eas
  for any such $\phi$. 
  Here we take $\phi(t):=e^{C_1 t} \zeta_\delta(t)$, $t>T$,  
  where $\zeta_\delta$ is as introduced in (\ref{zeta2}) for $\delta\in (0,t_0-T)$, with
  arbitrary Lebesgue points $t_0\in (T,T+1)$ and $t_1>T+1$ of $(T,\infty)\ni t \mapsto e^{C_1 t} y(t)$.
  Since $-\phi'(t)+C_1 \phi(t)=-e^{C_1 t} \zeta_\delta'(t)$ for a.e.~$t>T$, we thereby gain the inequality
  \bas
	\frac{1}{\delta} \int_{t_1}^{t_1+\delta} e^{C_1 t} y(t)dt
	- \frac{1}{\delta} \int_{t_0-\delta}^{t_0} e^{C_1 t} y(t)dt
	\le C_2 \int_T^\infty e^{C_1 t} \zeta_\delta(t) dt
	\qquad \mbox{for all } \delta\in (0,t_0-T),
  \eas
  which on taking $\delta\searrow 0$ shows that
  \bas
	y(t_1) 
	\le e^{-C_1 (t_1-t_0)} y(t_0) + C_2 \int_{t_0}^{t_1} e^{-C_1(t_1-t)} dt
	\le e^{-C_1 (t_1-t_0)} y(t_0) + \frac{C_2}{C_1}.
  \eas
  As the set of such Lebesgue points complements a null set in $(T,\infty)$, this implies that indeed
  \be{655.4}
	y(t) \le C_3:={\rm ess} \hspace*{-3mm} \inf_{\hspace*{-5mm} t_0\in (T,T+1)} y(t_0)+\frac{C_2}{C_1}
	\qquad \mbox{for a.e.~} t>T+1.
  \ee
  We next pick any $t>T+2$ and let $\phi(t):=\zeta_\delta(t)$, $t>T$, where again $\zeta_\delta$ is taken from 
  (\ref{zeta2}), now with $t_0:=t,t_1:=t+1$ and $\delta:=1$. From (\ref{655.2}) we thus obtain that
  \bas
	-\int_{t-1}^t y(s) ds 
	+ \int_{t+1}^{t+2} y(s)ds
	+ \frac{1}{K} \int_{t-1}^{t+2} h(s)\phi(s)ds
	\le K\int_{t-1}^{t+2} \phi(s)ds,
  \eas
  so that since $\phi\equiv 1$ in $(t,t+1)$ and $0\le \phi \le 1$ on $(T,\infty)$, using (\ref{655.4}) and the left
  inequality in (\ref{333.1}) we conclude that with $C_4:=\max\{1,M^3\}$ we have
  \bas
	\frac{1}{C_4 K} \int_t^{t+1} \io \bigg\{ \frac{|\nabla n|^2}{n} + |\nabla c|^4 \bigg\}
	&\le& \frac{1}{K} \int_t^{t+1} \io \bigg\{ \frac{|\nabla n|^2}{n} + \frac{|\nabla c|^4}{c^3} \bigg\} \\
	&\le& \frac{1}{K} \int_t^{t+1} h(s)ds \\
	&\le& C_3 + \frac{|\Omega|}{e} + 3K
	\qquad \mbox{for all } t>T+2.
  \eas
  Therefore, $(n,c,id)$ belongs to $\set$ with $L:=C_4 K\cdot(C_3 + \frac{|\Omega|}{e} + 3K)$ and $\tz:=T+2$,
  whereupon (\ref{c_decay}) becomes a consquence of Lemma \ref{lem46}.
\qed
\subsection{Preliminary statements on decay of $\nabla n$ and $u$}
Thus knowing that $c$ decays uniformly, invoking Lemma \ref{lem64} we obtain the following.
\begin{lem}\label{lem65}
  For any eventual energy solution $(n,c,u)$ of (\ref{0}), one can find $T>0$ such that
  \be{65.1}
	\int_T^\infty \io |\nabla n|^2 < \infty.
  \ee
  Moreover, for all $p\ge 2$ there exist $T(p)>0$ and $C(p)>0$ such that
  \be{65.2}
	\|n(\cdot,t)\|_{L^p(\Omega)} \le C(p)
	\qquad \mbox{for a.e.~} t>T(p),
  \ee
  and we have
  \be{65.11}
	\int_t^{t+1} \|n(\cdot,s)-\onz\|_{L^p(\Omega)} dt 
	\to 0
	\qquad \mbox{as } t\to\infty.
  \ee
\end{lem}
\proof
  In view of Lemma \ref{lem655}, we may apply Lemma \ref{lem64} firstly to $p:=2$ and $F(s):=s, \ s\ge 0$, 
  to obtain (\ref{65.1}), whereas
  (\ref{65.2}) similarly results on invoking Lemma \ref{lem64} to general $p\ge 2$.\\
  Next, given $p\ge 2$ we take $T(p)$ and 
  $C(p)$ as in (\ref{65.2}) and invoke the Poincar\'e inequality to find $C_1>0$ such that
  \bas
	\Big\|\varphi-\Mint_\Omega \varphi \Big\|_{L^2(\Omega)}^2 \le C_1 \io |\nabla\varphi|^2
	\qquad \mbox{for all } \varphi\in W^{1,2}(\Omega).
  \eas
  Then two applications of the H\"older inequality show that
  \bas
	\int_t^{t+1} \|n(\cdot,t)-\onz\|_{L^p(\Omega)} dt
	&\le& \int_t^{t+1} \|n(\cdot,t)-\onz\|_{L^{2p}(\Omega)}^\frac{p-2}{p-1} 
	\|n(\cdot,t)-\onz\|_{L^2(\Omega)}^\frac{1}{p-1} dt \\
	&\le& \Big(C(2p)\Big)^\frac{p-2}{p-1} \int_t^{t+1} \|n(\cdot,t)-\onz\|_{L^2(\Omega)}^\frac{1}{p-1} dt \\
	&\le& \Big(C(2p)\Big)^\frac{p-2}{p-1} \int_t^{t+1} \|n(\cdot,t)-\onz\|_{L^2(\Omega)}^2 dt \\
	&\le& \Big(C(2p)\Big)^\frac{p-2}{p-1} C_1 \int_t^{t+1} \io |\nabla n|^2  
	\qquad \mbox{for all } t>T(p),
  \eas
  whence (\ref{65.11}) is implied by (\ref{65.1}).
\qed
The stabilization property implied by (\ref{65.1}) can now be turned into a preliminary statement on decay of $u$ by making
use of the energy inequality (\ref{energy1}).
\begin{lem}\label{lem4333}
  Let $(n,c,u)$ be an eventual energy solution of (\ref{0}).
  Then there exists $T>0$ such that
  \be{4333.1}
	\int_T^\infty \io |\nabla u|^2 < \infty.
  \ee
  In particular, for all $p\in [1,6]$ we have
  \be{4333.2}
	\int_t^{t+1} \|u(\cdot,t)\|_{L^p(\Omega)} dt \to 0
	\qquad \mbox{as } t\to\infty.
  \ee
\end{lem}
\proof
  We combine Lemma \ref{lem65} with Definition \ref{defi_ees} to find $T_1>0$ such that
  $C_1:=\int_{T_1}^\infty \io |\nabla n|^2$ is finite, and such that 
  \bas
	\frac{1}{2} \io |u(\cdot,t)|^2
	+ \int_{t_0}^t \io |\nabla u|^2 
	\le \frac{1}{2} \io |u(\cdot,t_0)|^2
	+ \int_{t_0}^t \io nu\cdot\nabla\Phi
	\qquad \mbox{for a.e.~$t_0>T_1$ and all } t>t_0.
  \eas
  Since $\nabla\cdot u=0$, and since with some $C_2>0$ we have
  $\io |u|^2 \le C_2 \io |\nabla u|^2$ for a.e.~$t>T_1$ by the Poincar\'e inequality, integrating by parts in the
  rightmost integral and using Young's inequality we see that for any such $t_0$ and $t$ we have
  \bas
	\int_{t_0}^t \io |\nabla u|^2
	&\le& \frac{1}{2} \io |u(\cdot,t_0)|^2
	- \int_{t_0}^t \io  \Phi u\cdot\nabla n \\
	&\le& \frac{1}{2} \io |u(\cdot,t_0)|^2
	+ \frac{1}{2C_2} \int_{t_0}^t \io |u|^2
	+ \frac{C_2 \|\Phi\|_{L^\infty(\Omega)}^2}{2} \int_{t_0}^t \io |\nabla n|^2 \\
	&\le& \frac{1}{2} \io |u(\cdot,t_0)|^2
	+ \frac{1}{2} \int_{t_0}^t \io |\nabla u|^2
	+ \frac{C_1 C_2 \|\Phi\|_{L^\infty(\Omega)}^2}{2}.
  \eas
  This implies that
  \bas
	\int_{T_1+1}^t \io |\nabla u|^2
	\le \ {\rm ess} \hspace*{-3mm} \inf_{\hspace*{-5mm} t_0\in (T_1,T_1+1)} \io |u(\cdot,t_0)|^2
	+ C_1 C_2 \|\Phi\|_{L^\infty(\Omega)}^2
	\qquad \mbox{for all } t>T_1+1
  \eas
  and hence establishes (\ref{4333.1}), from which in turn one can readily derive (\ref{4333.2}),
  once again because $W^{1,2}(\Omega) \hra L^p(\Omega)$ for $p\le 6$.
\qed
The latter implies smallness of $u(\cdot,t_\star)$ at some conveniently large $t_\star>0$ 
in some of the spaces $L^{3+\eps}(\Omega)$, $\eps>0$,
which are supercritical with respect to the current knowledge on the global existence of smooth small-data
solutions to the unforced three-dimensional Navier-Stokes system (\cite{wiegner}). 
Thanks to the decay property of $n-\onz$ formulated in (\ref{65.11}), this actually entails a certain eventual 
regularity and decay of $u$ also in the present situation. 
More precisely, by means of a contraction mapping argument we can achieve the following.
\begin{lem}\label{lem57}
  Let $(n,c,u)$ be an eventual energy solution of (\ref{0}).
  Then for all $p>3$ we have
  \be{57.1}
	\|u\|_{L^\infty((t,\infty);L^p(\Omega))} \to 0
	\qquad \mbox{as } t\to\infty.
  \ee
\end{lem}
\proof
  First, since $p>3$ we can pick any $p_0\in (3,p)$ such that $p_0\le 6$ to achieve that then
  $\gamma:=\frac{3}{2}(\frac{1}{p_0}-\frac{1}{p})$ satisfies $\gamma<\frac{1}{2}-\frac{3}{2p}$, whence in particular
  \bas
	C_1:=\int_0^1 (1-\sigma)^{-\frac{1}{2}-\frac{3}{2p}} \sigma^{-2\gamma} d\sigma
  \eas
  is finite.
  We next recall known facts on the regularizing action of the Stokes semigroup (\cite{giga1986}) to fix positive
  constants $C_2, C_3$ and $C_4$ satisfying
  \be{57.2}
	\|e^{-tA}\varphi\|_{L^p(\Omega)}
	\le C_2 t^{-\gamma} \|\varphi\|_{L^{p_0}(\Omega)}
	\qquad \mbox{for all } \varphi\in C_0^\infty(\Omega) \cap L^2_\sigma(\Omega)
	\quad \mbox{and } t\in (0,2)
  \ee
  and 
  \be{57.3}
	\|e^{-tA}\proj [\nabla \cdot \varphi] \|_{L^p(\Omega)}
	\le C_3 t^{-\frac{1}{2}-\frac{3}{2p}}
	\|\varphi\|_{L^\frac{p}{2}(\Omega)}
	\qquad \mbox{for all } \varphi\in C_0^\infty(\Omega) \cap L^2_\sigma(\Omega)
	\quad \mbox{and } t\in (0,2)
  \ee
  as well as
  \be{57.4}
	\Big\|e^{-tA} \proj\Big[ \varphi - \Mint_\Omega \varphi\Big] \Big\|_{L^p(\Omega)} 
	\le C_4 \|\varphi\|_{L^p(\Omega)}
	\qquad \mbox{for all } \varphi\in C_0^\infty(\Omega) \cap L^2_\sigma(\Omega)
	\quad \mbox{and } t\in (0,2),
  \ee
  where (\ref{57.3}) in particular implies that for each $t\in (0,2)$ the operator $e^{-tA} \proj[\nabla\cdot ()]$
  admits a continuous extension to all of $L^\frac{p}{2}_\sigma(\Omega)$ with norm controlled according to (\ref{57.3}).\\
  We finally invoke the Cauchy-Schwarz inequality to fix $C_5>0$ such that
  \be{57.44}
	\|\varphi\mult\psi\|_{L^\frac{p}{2}(\Omega)} \le C_5\|\varphi\|_{L^p(\Omega)} \|\psi\|_{L^p(\Omega)}
	\qquad \mbox{for all $\varphi$ and $\psi$ belonging to } \in L^p(\Omega).
  \ee
  We now take any positive 
  $\delta\le \delta_0:=\big\{ 3\cdot 2^{\frac{1}{2}-\frac{3}{2p}-\gamma} \cdot C_1 C_3 C_5 \big\}^{-1}$
  and let $\delta_1:=\frac{\delta}{3C_2}$ and $\delta_2:=\frac{\delta}{3\cdot 2^\gamma C_4}$,
  and note that as a consequence of Lemma \ref{lem4333} and Lemma \ref{lem65}, for any such $\delta$ we can pick
  $T(\delta)>2$ such that
  \be{57.5}
	\int_t^{t+1} \|u(\cdot,s)\|_{L^{p_0}(\Omega)}^2 ds \le \delta_1^2
	\qquad \mbox{for all } t>T(\delta)-2
  \ee
  and
  \be{57.6}
	\|\nabla\Phi\|_{L^\infty(\Omega)} \cdot \int_t^{t+2} \|n(\cdot,s)-\onz\|_{L^p(\Omega)} ds \le \delta_2
	\qquad \mbox{for all } t>T(\delta)-2.
  \ee
  To see that these choices ensure that
  \be{57.7}
	\|u(\cdot,t_0)\|_{L^p(\Omega)} \le \delta
	\qquad \mbox{for all } t_0>T(\delta),
  \ee
  we fix any such $t_0$ and then infer from (\ref{57.5}) that there exists $t_\star \in (t_0-2,t_0-1)$ such that
  $\|u(\cdot,t_\star)\|_{L^{p_0}(\Omega)} \le \delta_1$. 
  We now follow a standard reasoning to construct, independently of $u$, another weak solution $\hu$ of the initial value
  problem associated with the Navier-Stokes system $\hu_t + A\hu = -\proj[\nabla \cdot (\hu\mult\hu)] + \proj[n\nabla\Phi]$
  in $\Omega\times (t_\star,t_\star+2)$ with $\hu(\cdot,t_\star)=u(\cdot,t_\star)$ and some favorable additional properties,
  finally implying by a uniqueness argument that actually $u=\hu$ and that hence $u$ itself has these properties.
  To this end, in the Banach space
  \bas
	X:=\Big\{ \varphi \in C^0((t_\star,t_\star+2];L^p_\sigma(\Omega)) \ \Big| \ 
	\|\varphi\|_X := \sup_{t\in (t_\star,t_\star+2)} (t-t_\star)^\gamma \|\varphi(\cdot,t)\|_{L^p(\Omega)} < \infty
	\Big\}
  \eas
  we consider the mapping $\Psi$ defined by
  \bas
	(\Psi\varphi)(\cdot,t)
	&:=& e^{-(t-t_\star)A} u(\cdot,t_\star)
	- \int_{t_\star}^t e^{-(t-s)A} \proj \Big[ \nabla \cdot (\varphi(\cdot,s)\mult \varphi(\cdot,s)) \Big] ds \\
	& & + \int_{t_\star}^t e^{-(t-s)A} \proj [n(\cdot,s)\nabla \Phi] ds,
	\qquad t\in (t_\star,t_\star+2],
  \eas
  for $\varphi$ belonging to the closed subset
  \bas
	S:= \Big\{ \varphi\in X \ \Big| \ \|\varphi\|_X \le \delta \Big\}
  \eas
  of $X$. 
  Then for $\varphi\in S$ we can use (\ref{57.2})-(\ref{57.4}) and (\ref{57.6}) to estimate
  \bas
	\|(\Psi\varphi)(\cdot,t)\|_{L^p(\Omega)}
	&\le& C_2(t-t_\star)^{-\frac{3}{2}(\frac{1}{p_0}-\frac{1}{p})} \|u(\cdot,t_\star)\|_{L^{p_0}(\Omega)}
	+ C_3 \int_{t_\star}^t (t-s)^{-\frac{1}{2}-\frac{3}{2p}} 
	\Big\|\varphi(\cdot,s)\mult \varphi(\cdot,s)\Big\|_{L^\frac{p}{2}(\Omega)} ds \\
	& & + C_4 \int_{t_\star}^t \Big\| [n(\cdot,s)-\onz] \nabla\Phi\Big\|_{L^p(\Omega)} ds \\[2mm]
	&\le& C_2 \delta_1 (t-t_\star)^{-\frac{3}{2}(\frac{1}{p_0}-\frac{1}{p})}
	+ C_3 C_5 \int_{t_\star}^t (t-s)^{-\frac{1}{2}-\frac{3}{2p}} \|\varphi(\cdot,s)\|_{L^p(\Omega)}^2 ds \\
	& & + C_4 \|\nabla\Phi\|_{L^\infty(\Omega)} 
	\int_{t_\star}^t \|n(\cdot,s)-\onz\|_{L^p(\Omega)} ds \\[2mm]
	&\le& C_2 \delta_1 (t-t_\star)^{-\frac{3}{2}(\frac{1}{p_0}-\frac{1}{p})}
	+ C_3 C_5 \delta^2 \int_{t_\star}^t (t-s)^{-\frac{1}{2}-\frac{3}{2p}} s^{-2\gamma} ds \\
	& & + C_4 \|\nabla\Phi\|_{L^\infty(\Omega)}
	\int_{t_\star}^{t_\star+2} \|n(\cdot,s)-\onz\|_{L^p(\Omega)} ds \\[2mm]
	&\le& C_2 \delta_1 (t-t_\star)^{-\frac{3}{2}(\frac{1}{p_0}-\frac{1}{p})}
	+ C_3 C_5 \delta^2 \int_{t_\star}^t (t-s)^{-\frac{1}{2}-\frac{3}{2p}} s^{-2\gamma} ds \\
	& & + C_4 \delta_2
	\qquad \mbox{for all } t\in (t_\star,t_\star+2],
  \eas
  so that according to our choice of $\gamma$ we obtain
  \bas
	(t-t_\star)^\gamma \|(\Psi\varphi)(\cdot,t)\|_{L^p(\Omega)}
	&\le& C_2 \delta_1
	+ C_3 C_5 \delta^2 (t-t_\star)^{\gamma-\frac{1}{2}-\frac{3}{2p}-2\gamma+1}
	\int_0^1 (1-\sigma)^{-\frac{1}{2}-\frac{3}{2p}} \sigma^{-2\gamma} d\sigma
	+ C_4 \delta_2 (t-t_\star)^\gamma \\
	&=& C_2 \delta_1 + C_1 C_3 C_5 \delta^2 (t-t_\star)^{\frac{1}{2}-\frac{3}{2p}-\gamma} 
	+ C_4 \delta_2 (t-t_\star)^\gamma \\
	&\le& C_2 \delta_1 + 2^{\frac{1}{2}-\frac{3}{2p}-\gamma} C_1 C_3 C_5 \delta^2
	+ 2^\gamma C_4 \delta_2 \\
	&\le& \frac{\delta}{3} + \frac{\delta}{3} + \frac{\delta}{3}=\delta
	\qquad \mbox{for all } t\in (t_\star,t_\star+2],
  \eas
  from which it readily follows that $\Psi S \subset S$.
  Likewise, for $\varphi\in S$ and $\psi\in S$ we can use (\ref{57.3}) and (\ref{57.44})
  to find that
  \bas
	\Big\| (\Psi\varphi-\Psi\psi)(\cdot,t)\Big\|_{L^p(\Omega)}
	&=& \Bigg\| \int_{t_\star}^t e^{-(t-s)A} \proj \bigg[
	-\nabla\cdot \Big\{\varphi(\cdot,s)\mult [\varphi(\cdot,s)-\psi(\cdot,s)] \Big\} \\
	& & \hspace*{29mm}
	- \nabla \cdot \Big\{ [\varphi(\cdot,s)-\psi(\cdot,s)] \mult \psi(\cdot,s) \Big\} \bigg] ds \Bigg\|_{L^p(\Omega)} 
	\\[2mm]
	&\le& C_3 \int_{t_\star}^t (t-s)^{-\frac{1}{2}-\frac{3}{2p}} 
	\bigg\{ \Big\| \varphi(\cdot,s)\mult [\varphi(\cdot,s)-\psi(\cdot,s)] \Big\|_{L^\frac{p}{2}(\Omega)} \\
	& & \hspace*{32mm}
	+ \Big\| [\varphi(\cdot,s)-\psi(\cdot,s)]\mult \psi(\cdot,s) \Big\|_{L^\frac{p}{2}(\Omega)} \bigg\} ds \\[2mm]
	&\le& C_3 C_5 \int_{t_\star}^t (t-s)^{-\frac{1}{2}-\frac{3}{2p}}
	\Big\{\|\varphi(\cdot,s)\|_{L^p(\Omega)} + \|\psi(\cdot,s)\|_{L^p(\Omega)} \Big\} \times \\[2mm]
	& & \hspace*{38mm}
	\times \|\varphi(\cdot,s)-\psi(\cdot,s)\|_{L^p(\Omega)} ds \\[2mm]
	&\le& C_3 C_5 \int_{t_\star}^t (t-s)^{-\frac{1}{2}-\frac{3}{2p}} \cdot 2\delta s^{-\gamma} 
	\cdot s^{-\gamma} \|\varphi-\psi\|_X ds
	\qquad \mbox{for all } t\in (t_\star,t_\star+2],
  \eas
  implying that
  \bas
	(t-t_\star)^\gamma
	\Big\| (\Psi\varphi-\Psi\psi)(\cdot,t)\Big\|_{L^p(\Omega)}
	\le 2 \cdot 2^{\frac{1}{2}-\frac{3}{2p}-\gamma} C_1 C_3 C_5 \delta \|\varphi-\psi\|_X
	\qquad \mbox{for all } t\in (t_\star,t_\star+2].
  \eas
  As $2 \cdot 2^{\frac{1}{2}-\frac{3}{2p}-\gamma} C_1 C_3 C_5 \delta \le \frac{2}{3}<1$, this proves that $\Phi$ acts
  as a contraction on $S$ and hence possesses a unique fixed point $\hu$.
  By standard arguments (\cite{sohr}), it follows that $\hu$ in fact is a weak solution of the Navier-Stokes subsystem
  of (\ref{0}) in $\Omega\times (t_\star,t_\star+2)$ subject to the initial condition $\hu(\cdot,t_\star)=u(\cdot,t_\star)$.
  Since furthermore our choice of $\gamma$ ensures that with $q:=\frac{2p}{p-3}$ we have 
  $q\gamma<q\cdot (\frac{1}{2}-\frac{3}{2p})=1$, it follows from the inclusion $\hu\in S$ that
  \bas
	\int_{t_\star}^{t_\star+2} \|\hu(\cdot,t)\|_{L^p(\Omega)}^q dt 
	\le \delta^q \int_{t_\star}^{t_\star+2} (t-t_\star)^{-q\gamma} dt < \infty.
  \eas
  Using that $p$ and $q$ satisfy the Serrin condition $\frac{2}{q}+\frac{3}{p}=1$, a well-known uniqueness property 
  of the Navier-Stokes equations  (\cite{sohr}) entails that $\hu$ must coincide with $u$ in 
  $\Omega\times (t_\star,t_\star+2)$. In particular, since $t_0\in (t_\star+1,t_\star+2)$, this implies that
  \bas
	\|u(\cdot,t_0)\|_{L^p(\Omega)}
	= \|\hu(\cdot,t_0)\|_{L^p(\Omega)} 
	\le \delta (t_0-t_\star)^{-\gamma} \le \delta
  \eas
  and thereby establishes (\ref{57.7}), which in turn proves (\ref{57.1}), because $\delta\in (0,\delta_0]$ was arbitrary.
\qed
\subsection{Eventual H\"older regularity of $u$ and $\nabla u$}
We next plan to derive some higher order regularity properties of a given eventual energy solution.
Here we first combine the boundedness feature of $u$ implied by Lemma \ref{lem57} with the integrability properties
of the forcing term $n\nabla\Phi$ in the Navier-Stokes equations in (\ref{0}), as obtained from Lemma \ref{lem65},
to achieve spatio-temporal $L^p$ bounds for $u, \nabla u, D^2 u$ and $u_t$ for any $p\ge 1$ by means of a bootstrap
argument based on maximal Sobolev regularity in the Stokes evolution system.\abs
For use in this and also the following sections, we fix a function $\xi_0\in C^\infty(\R)$ such that
\be{xi}
	0\le\xi_0\le 1 \ \mbox{in } \R,
	\quad
	\xi_0\equiv 0 \ \mbox{ in } (-\infty,\mbox{$\frac{1}{2}$}]
	\quad \mbox{and} \quad
	\xi_0\equiv 1 \ \mbox{ in } [1,\infty),
\ee
and for $t_0>1$, we introduce
\be{xitz}
	\xitz(t):=\xi_0(t-t_0), \qquad t\in\R.
\ee
\begin{lem}\label{lem58}
  Let $(n,c,u)$ be an eventual energy solution of (\ref{0}).
  Then for all $p\ge 1$ there exist $T>0$ and $C>0$ such that
  \be{58.1}
	\|u\|_{L^p((t,t+1);W^{2,p}(\Omega))} + \|u_t\|_{L^p(\Omega\times (t,t+1))} \le C
	\qquad \mbox{for all } t>T.
  \ee
\end{lem}
\proof
  We first claim that there exist $T_1>0$ and $C_1>0$ such that
  \be{58.2}
	\|u\|_{L^2((t_0,t_0+1);W^{2,2}(\Omega))} \le C_1
	\qquad \mbox{for all } t>T_1.
  \ee
  To see this, let us apply Lemma \ref{lem65}, Lemma \ref{lem57} and Lemma \ref{lem4333}
  to fix $T_1>1$ and positive constants $C_2, C_3$ and $C_4$ such that
  \be{58.3}
	\|n\|_{L^2((t,t+2);L^2(\Omega))} \le C_2
	\qquad \mbox{for all } t>T_1-1
  \ee
  and
  \be{58.4}
	\|u\|_{L^\infty((t,t+2);L^4(\Omega))} \le C_3
	\qquad \mbox{for all } t>T_1-1
  \ee
  and
  \be{58.5}
	\|\nabla u\|_{L^2((t,t+2);L^2(\Omega))} \le C_4
	\qquad \mbox{for all } t>T_1-1.
  \ee
  Then (\ref{58.5}) in particular implies that for any choice of $t_0>T_1$ we can find $t_\star\in (t_0-1,t_0)$ fulfilling
  $\|\nabla u(\cdot,t_\star)\|_{L^2(\Omega)} \le C_3$,
  and upon an interpolation using the H\"older inequality and the Gagliardo-Nirenberg inequality, (\ref{58.5})
  combined with (\ref{58.4}) shows that with some $C_5>0$ and $C_6>0$ we have
  \bas
	\int_{t_\star}^{t_\star+2} \Big\|\proj [(u\cdot \nabla)u](\cdot,t)\Big\|_{L^2(\Omega)}^2 dt
	&\le& C_5 \int_{t_\star}^{t_\star+2} \|\nabla u(\cdot,t)\|_{L^4(\Omega)}^2 \|u(\cdot,t)\|_{L^4(\Omega)}^2 dt \nn\\
	&\le& C_3^2 C_5 \int_{t_\star}^{t_\star+2} \|\nabla u(\cdot,t)\|_{L^4(\Omega)}^2 dt \\
	&\le& C_6 \int_{t_\star}^{t_\star+2} \|u(\cdot,t)\|_{W^{2,2}(\Omega)}^\frac{8}{5} 
	\|u(\cdot,t)\|_{L^4(\Omega)}^\frac{2}{5} dt \\
	&\le& C_3^\frac{2}{5} C_6 \int_{t_\star}^{t_\star+2} \|u(\cdot,t)\|_{W^{2,2}(\Omega)}^\frac{8}{5} dt.
  \eas
  As moreover
  \bas
	\int_{t_\star}^{t_\star+2} \Big\| \proj [n(\cdot,t)\nabla\Phi]\Big\|_{L^2(\Omega)}^2 dt
	&\le& \|\nabla\Phi\|_{L^\infty(\Omega)}^2 \int_{t_\star}^{t_\star+2} \|n(\cdot,t)\|_{L^2(\Omega)}^2 dt \\
	&\le& C_2^2 \|\nabla\Phi\|_{L^\infty(\Omega)}^2
  \eas
  by (\ref{58.3}), it follows from a well-known maximal Sobolev regularity property of the Stokes evolution equation
  (\cite{giga_sohr}) 
  and a corresponding uniqueness argument (\cite{sohr}) 
  that there exist $C_7>0$ and $C_8>0$ satisfying
  \bas
	\int_{t_\star}^{t_\star+2} \|u(\cdot,t)\|_{W^{2,2}(\Omega)}^2 dt 
	&\le& C_7 \cdot \bigg\{ \io |\nabla u(\cdot,t_\star)|^2
	+ \int_{t_\star}^{t_\star+2} \Big\| -\proj[(u\cdot\nabla)u](\cdot,t) + \proj[n(\cdot,t)\nabla\Phi] 
	\Big\|_{L^2(\Omega)}^2 dt \bigg\} \\
	&\le& C_8 \cdot \bigg\{ 1 + \int_{t_\star}^{t_\star+2} \|u(\cdot,t)\|_{W^{2,2}(\Omega)}^\frac{8}{5} dt \bigg\}.
  \eas
  Since $\frac{8}{5}<2$, Young's inequality becomes applicable here to warrant that (\ref{58.2}) indeed holds if we let
  $C_1>0$ be apporopriately large, because by construction we have $(t_0,t_0+1)\subset (t_\star,t_\star+2)$.\abs
  In order to prove the lemma, upon a recursive argument it is hence sufficient to show that whenever $p\ge 2$ is such
  that
  \be{58.6}
	\|u\|_{L^p((t,t+2);W^{2,p}(\Omega))} \le C_9
	\qquad \mbox{for all } t>T_2
  \ee
  with certain $T_2>1$ and $C_9>0$, there exist $T_3>1$ and $C_{10}>0$ fulfilling
  \be{58.7}
	\|u\|_{L^\frac{3p}{2}((t_0,t_0+1);W^{2,\frac{3p}{2}}(\Omega))}
	+ \|u_t\|_{L^\frac{3p}{2}(\Omega\times (t_0,t_0+1))}
	\le C_{10}
	\qquad \mbox{for all } t>T_3.
  \ee
  To see that this implication actually holds, under the assumption therein we invoke Lemma \ref{lem57} and \ref{lem65} to
  fix $T_3>T_2$, $C_{11}>0$ and $C_{12}>0$ such that with $q:=\max\{6p,\frac{6p}{2p-3}\}$ we have
  \be{58.89}
	\|n\|_{L^\infty((t-1,t+1);L^\frac{3p}{2}(\Omega))} \le C_{11}
	\qquad \mbox{for all } t>T_3
  \ee
  and
  \be{58.8}
	\|u\|_{L^\infty((t-1,t+1);L^q(\Omega))} \le C_{12}
	\qquad \mbox{for all } t>T_3,
  \ee
  and given $t_0>T_3$ we let $\xitz$ be as defined in (\ref{xitz}).
  Then the function $v: \Omega\times (t_0-1,\infty) \to \R^3$ defined by $v(x,t):=\xitz(t) u(x,t)$,
  $(x,t)\in \Omega\times (t_0-1,\infty)$, is a weak solution in $L^\infty_{loc}([t_0-1,\infty);L^2_\sigma(\Omega))
  \cap L^2_{loc}([t_0-1,\infty);W_0^{1,2}(\Omega) \cap W^{2,2}(\Omega))$ of 
  \be{58.88}
	v_t + Av = h(x,t):=\xitz(t) \proj[-(u\cdot\nabla)u + n\nabla \Phi] + \xitz'(t) u
	\qquad \mbox{in } \Omega\times (t_0-1,\infty)
  \ee
  with $v(\cdot,t_0-1)\equiv 0$.
  To estimate the inhomogeneity $h$ herein, we first note that the boundedness of the Helmholtz projection
  in $L^\frac{3p}{2}(\Omega)$, (\ref{58.89}) and (\ref{58.8}) imply that
  there exist positive constants $C_{13}, C_{14}$ and $C_{15}$ such that
  \bea{58.9}
	\int_{t_0-1}^{t_0+1} \Big\|\xitz(t) \proj[n(\cdot,t)\nabla\Phi]\Big\|_{L^\frac{3p}{2}(\Omega)}^\frac{3p}{2} dt
	&\le& C_{13} \int_{t_0-1}^{t_0+1} \|n(\cdot,t)\nabla\Phi\|_{L^\frac{3p}{2}(\Omega)}^\frac{3p}{2} dt \nn\\
	&\le& C_{13} \|\nabla\Phi\|_{L^\infty(\Omega)}^\frac{3p}{2} 
	\int_{t_0-1}^{t_0+1} \|n(\cdot,t)\|_{L^\frac{3p}{2}(\Omega)}^\frac{3p}{2} dt \nn\\[2mm]
	&\le& C_{14}
  \eea
  and
  \be{58.10}
	\int_{t_0-1}^{t_0+1} \big\|\xitz'(t) u(\cdot,t)\big\|_{L^\frac{3p}{2}(\Omega)}^\frac{3p}{2} dt
	\le C_{15}.
  \ee
  Moreover, by means of another interpolation on the basis of the H\"older inequality and the 
  Gagliardo-Nirenberg inequality we may use (\ref{58.89})
  and then (\ref{58.6}) to find positive constants $C_{16}, C_{17}, C_{18}$ and $C_{19}$ such that
  \bea{58.11}
	\int_{t_0-1}^{t_0+1} \Big\|\xitz(t) \proj[(u\cdot\nabla)u](\cdot,t)\Big\|_{L^\frac{3p}{2}(\Omega)}^\frac{3p}{2} dt
	&\le& C_{16} \int_{t_0-1}^{t_0+1} \|\nabla u(\cdot,t)\|_{L^{2p}(\Omega)}^\frac{3p}{2}
	\|u(\cdot,t)\|_{L^{6p}(\Omega)}^\frac{3p}{2} dt \nn\\
	&\le& C_{17} \int_{t_0-1}^{t_0+1} \|\nabla u(\cdot,t)\|_{L^{2p}(\Omega)}^\frac{3p}{2} dt \nn\\
	&\le& C_{18} \int_{t_0-1}^{t_0+1} \|u(\cdot,t)\|_{W^{2,p}(\Omega)}^p 
	\|u(\cdot,t)\|_{L^\frac{6p}{2p-3}(\Omega)}^\frac{p}{2} dt \nn\\
	&\le& C_{19} \int_{t_0-1}^{t_0+1} \|u(\cdot,t)\|_{W^{2,p}(\Omega)}^p dt \nn\\
	&\le& C_{19} C_9^p.
  \eea
  As a consequence of (\ref{58.9}), (\ref{58.10}) and (\ref{58.11}), once more by maximal Sobolev regularity estimates,
  now applied to (\ref{58.88}), we obtain $C_{20}>0$ and $C_{21}>0$ satisfying
  \bas
	& & \hspace*{-40mm}
	\int_{t_0}^{t_0+1} \|u(\cdot,t)\|_{W^{2,\frac{3p}{2}}(\Omega)}^\frac{3p}{2} dt
	+ \int_{t_0}^{t_0+1} \|u_t(\cdot,t)\|_{L^\frac{3p}{2}(\Omega)}^\frac{3p}{2} dt \\
	&\le& \int_{t_0-1}^{t_0+1} \|v(\cdot,t)\|_{W^{2,\frac{3p}{2}}(\Omega)}^\frac{3p}{2} dt
	+ \int_{t_0-1}^{t_0+1} \|v_t(\cdot,t)\|_{L^\frac{3p}{2}(\Omega)}^\frac{3p}{2} dt \\
	&\le& C_{20} \int_{t_0-1}^{t_0+1} \|h(\cdot,t)\|_{L^\frac{3p}{2}(\Omega)}^\frac{3p}{2} dt \nn\\[2mm]
	&\le& C_{21}.
  \eas
  This establishes (\ref{58.7}) and thereby completes the proof.
\qed
Since the exponent $p$ in Lemma \ref{lem58} can be chosen arbitrarily large, an immediate consequence is the following.
\begin{cor}\label{cor588}
  Suppose that $(n,c,u)$ is an eventual energy solution of (\ref{0}).
  Then one can find $\alpha\in (0,1)$, $T>0$ and $C>0$ such that
  \bas
	\|u\|_{C^{1+\alpha,\alpha}(\bar\Omega\times [t,t+1])} \le C 
	\qquad \mbox{for all } t>T.
  \eas
\end{cor}
\proof
  According to a well-known embedding result (\cite{amann}), for any $\alpha>0$ and $\beta>0$ fulfilling $\alpha+2\beta<2$
  there exists $p>1$ such that for each bounded interval $J\subset\R$, 
  the space of functions $\varphi$ on $\Omega\times (0,T)$ having finite norm
  $\|\varphi\|_{L^p(J;W^{2,p}(\Omega))} + \|\varphi_t\|_{L^p(J;L^p(\Omega))}$ is continuously embedded into
  $C^{\alpha,\beta}(\bar\Omega\times \bar J)$.
  Therefore, the claim is an immediate consequence of Lemma \ref{lem58} when applied to conveniently large $p\ge 2$.
\qed
\subsection{Eventual $L^p$ regularity of $c,\nabla c$ and $D^2 c$. H\"older regularity of $c$ and $\nabla c$}
By pursuing a similar overall strategy, we can derive the counterpart of the statement in Lemma \ref{lem58} for
the second component $c$.
As compared to the situation in the previous section, however,
the different structure of the inhomogeneity $h(x,t)=-nf(c)-u\cdot\nabla c$ in $c_t-\Delta c=h(x,t)$, and especially  
its dependence on $c$, require modifications in the argument.
\begin{lem}\label{lem500}
  Let $(n,c,u)$ be an eventual energy solution of (\ref{0}).
  Then for all $p\ge 1$ there exist $T>0$ and $C>0$ such that
  \be{500.1}
	\|c\|_{L^p((t,t+1);W^{2,p}(\Omega))} + \|c_t\|_{L^p(\Omega\times (t,t+1))} \le C
	\qquad \mbox{for all } t>T.
  \ee
\end{lem}
\proof
  Let us first make sure that with some $T_1>0$ and $C_1>0$ we have
  \be{500.2}
	\|c\|_{L^4((t_0,t_0+1);W^{2,4}(\Omega))} \le C
	\qquad \mbox{for all } t>T.
  \ee
  To this end, we observe that in view of Lemma \ref{lem65}, Definition \ref{defi_ees}, (\ref{cinfty}) and
  Corollary \ref{cor588}
  there exist $T_1>1$ and positive constants $C_2, C_3$ and $C_4$ such that
  \be{500.3}
	\|n\|_{L^4(\Omega\times (t-1,t+1))} \le C_2
	\qquad \mbox{for all } t>T_1
  \ee
  and
  \be{500.4}
	\|c\|_{L^4((t-1,t+1);W^{1,4}(\Omega))} \le C_3
	\qquad \mbox{for all } t>T_1
  \ee
  as well as
  \be{500.5}
	\|u\|_{L^\infty(\Omega\times (t-1,t+1))} \le C_4
	\qquad \mbox{for all } t>T_1.
  \ee
  Now for fixed $t_0>T_1$ we let $\xitz$ be as given by (\ref{xitz}) and consider the problem
  \be{500.6}
	\left\{ \begin{array}{l}
	z_t-\Delta z = h(x,t):=-\xitz(t) \cdot \big\{ n f(c)+ u\cdot\nabla c \big\} + \xitz'(t) c,
	\qquad x\in\Omega, \ t>t_0-1, \\[1mm]
	\frac{\partial z}{\partial\nu}=0, \qquad x\in\pO, \ t>t_0+1, \\[1mm]
	z(x,t_0-1)=0, \qquad x\in\Omega.
	\end{array} \right.
  \ee
  Here using (\ref{500.3})-(\ref{500.5}) and (\ref{cinfty}) we see that for some $C_5>0$
  we have
  \bea{500.7}
	& & \hspace*{-35mm}
	\int_{t_0-1}^{t_0+1} \Big\| \xitz(t) \Big\{ n(\cdot,t) f(c(\cdot,t)) + u(\cdot,t)\cdot\nabla c(\cdot,t)
	\Big\}\Big\|_{L^4(\Omega)}^4 dt \nn\\
	&\le& C_5 \int_{t_0-1}^{t_0+1} \Big\{ \|n(\cdot,t)\|_{L^4(\Omega)}^4 + 
	\|u(\cdot,t)\|_{L^\infty(\Omega)}^4 \|\nabla c(\cdot,t)\|_{L^4(\Omega)}^4 \Big\} dt \nn\\[2mm]
	&\le& C_5 \cdot ( C_2^4 + C_4^4 C_3^4)
  \eea
  and
  \be{500.8}
	\int_{t_0-1}^{t_0+1} \|\xitz'(t) c(\cdot,t)\|_{L^4(\Omega)}^4 dt
	\le \|\xi_0'\|_{L^\infty(\R)}^4 \cdot C_3^4.
  \ee
  According to well-known results on maximal Sobolev regularity properties of the Neumann heat semigroup
  (\cite{giga_sohr}), we thus infer from (\ref{500.7}) and (\ref{500.8}) that (\ref{500.6}) possesses a unique strong
  solution $z\in L^4((t_0-1,t_0+1);W^{2,4}(\Omega))$ with $z_t\in L^4(\Omega\times (t_0-1,t_0+1))$ which satisfies
  \bea{500.9}
	\int_{t_0-1}^{t_0+1} \|z(\cdot,t)\|_{W^{2,4}(\Omega)}^4 dt
	+ \int_{t_0-1}^{t_0+1} \|z_t(\cdot,t)\|_{L^4(\Omega)}^4 dt
	\le C_6 \int_{t_0-1}^{t_0+1} \|h(\cdot,t)\|_{L^4(\Omega)}^4 dt \le C_7
  \eea
  with some $C_6>0$ and $C_7>0$.
  Since clearly both $z$ and the function $\Omega\times (t_0-1,t_0+1)\ni(x,t) \mapsto \xitz(t) c(x,t)$
  are weak solutions of (\ref{500.6}) in the class of functions from $L^2((t_0-1,t_0+1);W^{1,2}(\Omega))$, it follows
  from a corresponding uniqueness property that $z(x,t)=\xitz(t)c(x,t)$ for a.e.~$(x,t)\in\Omega\times (t_0-1,t_0+1)$,
  whereupon (\ref{500.9}) implies (\ref{500.2}), because $\xitz \equiv 1$ in $(t_0,t_0+1)$.\\
  Let us next verify that if
  \be{500.10}
	\|c\|_{L^p((t,t+1);W^{2,p}(\Omega))} \le C_8
	\qquad \mbox{for all } t>T_2
  \ee
  with some $p>\frac{3}{2}, T_2>1$ and $C_8>0$, then there exist $T_3>T_2$ and $C_9>0$ such that
  \be{500.11}
	\|c\|_{L^{2p}((t_0,t_0+1);W^{2,2p}(\Omega))}
	+ \|c_t\|_{L^{2p}(\Omega\times (t_0,t_0+1))} 
	\le C_9 
	\qquad \mbox{for all } t>T_3.
  \ee
  Indeed, assuming (\ref{500.10}) we once again invoke Lemma \ref{lem65} to obtain
  $T_3>T_2$ and $C_{10}>0$ such that $T_3>T_1$ and
  \be{500.12}
	\|n\|_{L^{2p}(\Omega\times (t_0-1,t_0+1))} \le C_{10}
	\qquad \mbox{for all } t>T_3.
  \ee
  Then given $t_0>T_3$ we define $z$ in the same manner as before and the see that in (\ref{500.6}) we can use
  (\ref{500.12}), (\ref{500.5}) and (\ref{cinfty}) to find $C_{11}>0$ fulfilling
  \bas
	& & \hspace*{-30mm}
	\int_{t_0-1}^{t_0+1} \Big\| \xitz(t) \Big\{ n(\cdot,t) f(c(\cdot,t)) + u(\cdot,t)\cdot\nabla c(\cdot,t)
	\Big\} \Big\|_{L^{2p}(\Omega)}^{2p} dt \\
	&\le& C_{11} \int_{t_0-1}^{t_0+1} \Big\{ \|n(\cdot,t)\|_{L^{2p}(\Omega)}^{2p} + 
	\|u(\cdot,t)\|_{L^\infty(\Omega)}^{2p} \|\nabla c(\cdot,t)\|_{L^{2p}(\Omega)}^{2p} \Big\} dt \nn\\
	&\le& C_{10}^{2p} C_{11} + C_4^{2p} C_{11} \int_{t_0-1}^{t_0+1} \|\nabla c(\cdot,t)\|_{L^{2p}(\Omega)}^{2p} dt,
  \eas
  where combining the Gagliardo-Nirenberg inequality with (\ref{500.10}) and (\ref{cinfty}) provides
  $C_{12}>0$ and $C_{13}>0$ such that
  \bas
	\int_{t_0-1}^{t_0+1} \|\nabla c(\cdot,t)\|_{L^{2p}(\Omega)}^{2p} dt
	&\le& C_{12} \int_{t_0-1}^{t_0+1} \|c(\cdot,t)\|_{W^{2,p}(\Omega)}^p \|c(\cdot,t)\|_{L^\infty(\Omega)}^p dt \\[2mm]
	&\le& C_{13}.		
  \eas
  Since again from (\ref{cinfty}) we obtain $C_{14}>0$ such that
  \bas
	\int_{t_0-1}^{t_0+1} \|\xitz'(t) c(\cdot,t)\|_{L^{2p}(\Omega)}^{2p} dt
	\le 2\|\xi_0'\|_{L^\infty(\R)}^{2p} \|c\|_{L^\infty(\Omega\times (t_0-1,t_0+1))}^{2p}
	\le C_{14},
  \eas
  we may argue as above to conclude from maximal Sobolev regularity estimates that there exist $C_{15}>0$ and $C_{16}>0$
  such that
  \bas
	\int_{t_0-1}^{t_0+1} \|z(\cdot,t)\|_{W^{2,2p}(\Omega)}^{2p} dt
	+ \int_{t_0-1}^{t_0+1} \|z_t(\cdot,t)\|_{L^{2p}(\Omega)}^{2p} dt
	\le C_{15} \int_{t_0-1}^{t_0+1} \|h(\cdot,t)\|_{L^{2p}(\Omega)}^{2p} dt \le C_{16}.
  \eas
  Again since $z=c$ a.e.~in $\Omega\times (t_0,t_0+1)$ by definition of $\xitz$, 
  this shows (\ref{500.11}) and thereby completes the proof
  of (\ref{500.1}) upon iteration.
\qed
Again, this implies H\"older estimates as follows.
\begin{cor}\label{cor520}
  Let $(n,c,u)$ be an eventual energy solution of (\ref{0}). 
  Then there exist $\alpha\in (0,1), T>0$ and $C>0$ such that
  \bas
	\|c\|_{C^{1+\alpha,\alpha}(\bar\Omega\times [t,t+1])}
	\le C
	\qquad \mbox{for all } t>T.
  \eas
\end{cor}
\proof 
  In precisely the same manner as Corollary \ref{cor588} was derived from Lemma \ref{lem58}, this follows
  from Lemma \ref{lem500} by application of a standard embedding result (\cite{amann}).
\qed
\subsection{Eventual $L^p$ regularity of $n,\nabla n$ and $D^2 n$. H\"older regularity of $n$ and $\nabla n$}
The above estimates, inter alia asserting a uniform pointwise bound for $\nabla c$ and a 
space-time bound for $\Delta c$ in any $L^p$ norm for $p<\infty$,
now provide sufficient information for the derivation of the following analogue of Lemma \ref{lem500} for $n$.
\begin{lem}\label{lem510}
  Given an eventual energy solution $(n,c,u)$ of (\ref{0}),
  for each $p\ge 1$ one can find $T>0$ and $C>0$ such that
  \be{510.1}
	\|n\|_{L^p((t,t+1);W^{2,p}(\Omega))} + \|n_t\|_{L^p(\Omega\times (t,t+1))} \le C
	\qquad \mbox{for all } t>T.
  \ee
\end{lem}
\proof
  The proof is quite similar to that of Lemma \ref{lem500}, so that we may confine ourselves with outlining the main steps
  only.
  We first show that there exist $T_1>0$ and $C_1>0$ such that
  \be{510.2}
	\|n\|_{L^2((t_0,t_0+1);W^{2,2}(\Omega))} \le C_1
	\qquad \mbox{for all } t_0>T_1.
  \ee
  To see this, we recall Lemma \ref{lem65}, Corollary \ref{cor520}, Lemma \ref{lem500} and Corollary \ref{cor588}
  to fix $T_1>1$ and $C_2>0$ such that
  \bea{510.22}
	& & \hspace*{-10mm}
	\|n\|_{L^2((t-1,t+1);W^{1,2}(\Omega))}
	+ \|n\|_{L^4(\Omega\times (t-1,t+1))}
	+ \|\nabla c\|_{L^\infty(\Omega\times (t-1,t+1))}
	+ \|\Delta c\|_{L^4(\Omega\times (t-1,t+1))} \nn\\[2mm]
	& & \hspace*{60mm}
	+ \|u\|_{L^\infty(\Omega\times (t-1,t+1))}
	\le C_1
	\qquad \mbox{for all } t>T_1.
  \eea
  Then given $t_0>T_1$, with $\xitz$ taken from (\ref{xitz}) we see that the source term $h$ in
  \be{510.3}
	\left\{ \begin{array}{l}
	w_t - \Delta w = h(x,t):=-\xitz(t) \cdot \Big\{ \nabla n \cdot \nabla c + n\Delta c + u\cdot\nabla n \Big\}
	+ \xitz'(t) n,
	\qquad x\in \Omega, \ t>t_0-1, \\[1mm]
	\frac{\partial w}{\partial\nu}=0, \qquad x\in \pO, \ t>t_0-1, \\[1mm]
	w(x,t_0-1)=0, \qquad x\in \pO,
	\end{array} \right.
  \ee
  satisfies
  \bas
	\int_{t_0-1}^{t_0+1} \|h(\cdot,t)\|_{L^2(\Omega)}^2 dt
	&\le& 
	\|\nabla c\|_{L^\infty(\Omega\times (t_0-1,t_0+1))}^2 \int_{t_0-1}^{t_0+1} \|\nabla n(\cdot,t)\|_{L^2(\Omega)}^2 dt \\
	& & \hspace*{0mm}
	+ \bigg( \int_{t_0-1}^{t_0+1} \|n(\cdot,t)\|_{L^4(\Omega)}^2 \bigg)^\frac{1}{2} \cdot
	\bigg(\int_{t_0-1}^{t_0+1} \|\Delta c(\cdot,t)\|_{L^4(\Omega)}^2 \bigg)^\frac{1}{2} \\
	& & \hspace*{0mm}
	+ \|u\|_{L^\infty(\Omega\times (t_0-1,t_0+1))} \int_{t_0-1}^{t_0+1} \|\nabla n(\cdot,t)\|_{L^2(\Omega)}^2 dt \nn\\
	& & + \|\xi_0'\|_{L^\infty(\R)}^2 \int_{t_0-1}^{t_0+1} \|n(\cdot,t)\|_{L^2(\Omega)}^2 dt \nn\\[2mm]
	&\le& C_2
  \eas
  for some $C_2>0$, where we have used the Cauchy-Schwarz inequality.
  Therefore, (\ref{510.2}) is a consequence of a maximal Sobolev regularity inequality along with a uniqueness argument
  applied to (\ref{510.3}).\\
  Next, assuming that with some $p>1, T_2>1$ and $C_3>0$ we have
  \be{510.4}
	\|n\|_{L^p((t-1,t+1);W^{2,p}(\Omega))} \le C_3
	\qquad \mbox{for all } t>T_2,
  \ee
  we can infer the existence of $T_3>0$ and $C_4>0$ fulfilling
  \be{510.5}
	\|n\|_{L^\frac{3p}{2}((t_0,t_0+1);W^{2,\frac{3p}{2}}(\Omega))}
	+ \|n_t\|_{L^\frac{3p}{2}(\Omega\times (t_0,t_0+1))}
	\le C_4
	\qquad \mbox{for all } t_0>T_3.
  \ee
  Indeed, thanks to Lemma \ref{lem65} and Lemma \ref{lem500} we may pick $T_3>\max\{T_1,T_2\}$ and $C_5>0$ such that
  \bas
	\|n\|_{L^{3p}(\Omega\times (t-1,t+1))} 
	+ \|n\|_{L^\infty((t-1,t+1);L^3(\Omega))}
	+ \|\Delta c\|_{L^{3p}(\Omega\times (t-1,t+1))}
	\le C_5
	\qquad \mbox{for all } t>T_3,
  \eas
  whence by (\ref{510.4}) and the Gagliardo-Nirenberg inequality we obtain $C_6>0$ satisfying
  \bas
	\int_{t-1}^{t+1} \|n(\cdot,s)\|_{W^{1,\frac{3p}{2}}(\Omega)}^\frac{3p}{2} ds
	&\le& C_6 \int_{t-1}^{t+1} \|n(\cdot,s)\|_{W^{2,p}(\Omega)}^p \|n(\cdot,s)\|_{L^3(\Omega)}^\frac{p}{2} ds
	\le C_6 C_3^p C_5^\frac{p}{2}
	\qquad \mbox{for all } t>T_3.
  \eas
  Accordingly, if $t_0>T_3$ and $\xitz$ is as in (\ref{xitz}), then in (\ref{510.3}) we can once more use 
  (\ref{510.22}) and the Cauchy-Schwarz inequality to estimate
  \bas
	\int_{t_0-1}^{t_0+1} \|h(\cdot,t)\|_{L^\frac{3p}{2}(\Omega)}^\frac{3p}{2} dt
	&\le& \|\nabla c\|_{L^\infty(\Omega\times (t_0-1,t_0+1))}^\frac{3p}{2} 
	\int_{t_0-1}^{t_0+1} \|\nabla n(\cdot,t)\|_{L^\frac{3p}{2}(\Omega)}^\frac{3p}{2} dt \\
	& & + \bigg( \int_{t_0-1}^{t_0+1} \|n(\cdot,t)\|_{L^{3p}(\Omega)}^{3p} dt \bigg)^\frac{1}{2}
	\cdot \bigg( \int_{t_0-1}^{t_0+1} \|\Delta c(\cdot,t)\|_{L^{3p}(\Omega)}^{3p} \bigg)^\frac{1}{2} \\
	& & + \|u\|_{L^\infty(\Omega\times (t_0-1,t_0+1))}^\frac{3p}{2}
	\int_{t_0-1}^{t_0+1} \|\nabla n(\cdot,t)\|_{L^\frac{3p}{2}(\Omega)}^\frac{3p}{2} dt \\
	& & + \|\xi_0'\|_{L^\infty(\R)}^\frac{3p}{2}
	\int_{t_0-1}^{t_0+1} \|n(\cdot,t)\|_{L^\frac{3p}{2}(\Omega)}^\frac{3p}{2} dt \\[3mm]
	&\le& C_7
  \eas
  with some $C_7>0$. Another application of maximal Sobolev regularity theory thus yields (\ref{510.5}) and hence proves
  (\ref{510.1}), because $p>1$ was arbitrary.
\qed
Once more, this entails a certain H\"older regularity.
\begin{cor}\label{cor530}
  Let $(n,c,u)$ be an eventual energy solution of (\ref{0}). 
  Then there exist $\alpha\in (0,1), T>0$ and $C>0$ such that
  \bas
	\|n\|_{C^{1+\alpha,\alpha}(\bar\Omega\times [t,t+1])}
	\le C
	\qquad \mbox{for all } t>T.
  \eas
\end{cor}
\proof 
  In precisely the same manner as Corollary \ref{cor588} was derived from Lemma \ref{lem58}, this follows
  from Lemma \ref{lem510} by application of a standard embedding result (\cite{amann}).
\qed
\subsection{Estimates in $C^{2+\alpha,1+\frac{\alpha}{2}}$}
Straightforward applications of standard Schauder estimates for the Stokes evolution equation and the heat equation,
respectively, finally yield eventual smoothness of the solution components $u$ as well as $n$ and $c$, respectively.
\begin{lem}\label{lem589}
  For any eventual energy solution $(n,c,u)$ (\ref{0}) one can find
  $\alpha\in (0,1), T>0$ and $C>0$ such that
  \be{589.1}
	\|u\|_{C^{2+\alpha,1+\frac{\alpha}{2}}(\bar\Omega\times [t,t+1])} \le C
	\qquad \mbox{for all } t>T.
  \ee
\end{lem}
\proof
  According to Corollary \ref{cor530} and Corollary \ref{cor588}, 
  it is possible to fix $\alpha'\in (0,1)$, $T_1>0$ and $C_2>0$
  such that
  \be{589.2}
	\|n\|_{C^{\alpha',\frac{\alpha'}{2}}(\bar\Omega\times [t,t+1])} \le C_1
	\qquad \mbox{for all } t>T_1
  \ee
  and 
  \be{589.3}
	\|u\|_{C^{1+\alpha',\alpha'}(\bar\Omega\times [t,t+1])} \le C_2
	\qquad \mbox{for all } t>T_1.
  \ee
  We next set $T:=T_1+1$ and let $t_0>T$ be given. Then with $\xitz$ taken from (\ref{xitz}), we again use that
  $v(\cdot,t):=\xitz(t) u(\cdot,t)$, $t>t_0-1$, is a solution of
  \be{589.33}
	\left\{ \begin{array}{l}
	v_t + Av = h(x,t):=\xitz(t) \proj[(u\cdot\nabla)u + n\nabla \Phi] + \xitz'(t) u,
	\qquad x\in\Omega, \ t>t_0-1, \\[1mm]
	v(x,t_0)=0, \qquad x\in\Omega,
	\end{array} \right.
  \ee
  which hence, in particular, satisfies the associated first-order compatibility condition
  at $t=t_0-1$ (\cite{solonnikov2007}).\\
  Now from (\ref{589.2}), (\ref{589.3}) and the smoothness of $\xi_0$ we readily obtain $\alpha\in (0,1)$
  and $C_3>0$ fulfilling
  \bas
	\|h\|_{C^{\alpha,\frac{\alpha}{2}}(\bar\Omega\times [t_0-1,t_0+1])} \le C_3,
  \eas
  so that regularity estimates from Schauder theory for the Stokes evolution equation (\cite{solonnikov2007}) ensure
  that (\ref{589.33}) possesses a classical solution 
  $\tilde v\in C^{2+\alpha,1+\frac{\alpha}{2}}(\bar\Omega\times [t_0-1,t_0+1])$ satisfying
  \bas
	\|\tilde v\|_{C^{2+\alpha,1+\frac{\alpha}{2}}(\bar\Omega\times [t_0-1,t_0+1])} \le C_4
  \eas
  with some $C_4>0$ which is independent of $t_0$.
  As clearly $\tilde v\equiv v$ by an evident uniqueness property of (\ref{589.33}), this proves (\ref{589.1}).
\qed
\begin{lem}\label{lem556}
  Assume that $(n,c,u)$ be an eventual energy solution of (\ref{0}).
  Then there exist $\alpha\in (0,1), T>0$ and $C>0$ such that
  \be{556.1}
	\|n\|_{C^{2+\alpha,1+\frac{\alpha}{2}}(\bar\Omega\times [t,t+1])}
	+ \|c\|_{C^{2+\alpha,1+\frac{\alpha}{2}}(\bar\Omega\times [t,t+1])} \le C
	\qquad \mbox{for all } t>T.
  \ee
\end{lem}
\proof
  We first combine Corollary \ref{cor520} with Corollary \ref{cor530} and Corollary \ref{cor588} 
  to infer the existence of $\alpha'\in (0,1)$, $T_1>0$ and $C_1>0$ such that
  \bas
	\|nf(c)\|_{C^{\alpha',\frac{\alpha'}{2}}(\bar\Omega\times [t,t+1])} 
	+ \|u\cdot\nabla c\|_{C^{\alpha',\frac{\alpha'}{2}}(\bar\Omega\times [t,t+1])}
	\le C_1
	\qquad \mbox{for all } t>T_1.
  \eas
  Standard parabolic Schauder estimates applied to the second equation in (\ref{0})
  (\cite{LSU}) thus provide $C_2>0$ fulfilling
  \be{556.2}
	\|c\|_{C^{2+\alpha',1+\frac{\alpha'}{2}}(\bar\Omega\times [t,t+1])} \le C_2
	\qquad \mbox{for all } t>T_1+1.
  \ee
  Again in view of Corollary \ref{cor530} and Corollary \ref{cor588}, this in turn warrants that for some 
  $\alpha''\in (0,1), T_2>0$ and $C_3>0$ we have
  \bas
	\|\nabla\cdot (n\nabla c)\|_{C^{\alpha'',\frac{\alpha''}{2}}(\bar\Omega\times [t,t+1])}
	+ \|u\cdot\nabla n\|_{C^{\alpha'',\frac{\alpha''}{2}}(\bar\Omega\times [t,t+1])}
	\le C_3
	\qquad \mbox{for all } t>T_2,
  \eas
  whereupon Schauder theory says that
  \bas
	\|n\|_{C^{2+\alpha'',1+\frac{\alpha''}{2}}(\bar\Omega\times [t,t+1])} \le C_4
	\qquad \mbox{for all } t>T_2+1.
  \eas
  Along with (\ref{556.2}), this proves (\ref{556.1}) with $\alpha:=\min\{\alpha',\alpha''\}, T:=\max\{T_1,T_2\}$
  and some suitably large $C>0$.
\qed
\mysection{Stabilization of $n$ and $u$. Proof of Theorem \ref{theo_eventual}}\label{sect8}
On the basis of the eventual uniform continuity properties implied by the estimates in the previous section, 
we can now turn the weak stabilization properties of $n$ and $u$ 
from Lemma \ref{lem65} and Lemma \ref{lem57} into convergence with regard to the norm in $L^\infty(\Omega)$.
In the proofs of our results in Lemma \ref{lem600} and Lemma \ref{lem666} in this direction 
we shall make use of the following statement, the elementary proof of which may be omitted here.
\begin{lem}\label{lem610}
  Let $T\in\R$, and assume that $h: (T,\infty) \to [0,\infty)$ is uniformly continuous and such that
  $\int_t^{t+1} h(s) ds \to 0$ and $t\to\infty$. Then $h(t)\to 0$ as $t\to\infty$.
\end{lem}
A first application thereof shows that Corollary \ref{cor520} and Lemma \ref{lem65} entail the following.
\begin{lem}\label{lem600}
  Assume that $(n,c,u)$ be an eventual energy solution of (\ref{0}).
  Then with $\onz:=\mint_\Omega n_0>0$, we have
  \be{600.1}
	\|n(\cdot,t)-\onz\|_{L^\infty(\Omega)} \to 0
	\qquad \mbox{as } t\to\infty.
  \ee
\end{lem}
\proof
  By means of Corollary \ref{cor520} and the Arzel\`a-Ascoli theorem, we can fix $T_1>0$ such that 
  \be{600.2}
	(n(\cdot,t))_{t>T_1}
	\mbox{ is relatively compact in } C^0(\bar\Omega), 
  \ee
  and that 
  the function $\varphi:(T_1,\infty)\to\R$ defined by $h(t):=\|n(\cdot,t)-\onz\|_{L^2(\Omega)}^2$, $t>T_1$,
  is uniformly continuous. 
  Moreover, from Lemma \ref{lem65} we know that with some $T_2>T_1$ we have $\int_{T_2}^\infty \io |\nabla n|^2<\infty$,
  so that since the Poincar\'e inequality combined with (\ref{mass_ees}) provides $C_1>0$ such that
  \bas
	\int_t^{t+1} \|n(\cdot,s)-\onz\|_{L^2(\Omega)}^2 ds \le C_1 \int_t^{t+1} \io |\nabla n|^2 
	\qquad \mbox{for all } t>T_2,
  \eas
  we obtain that $\int_t^{t+1} h(s)ds \to 0$ as $t\to\infty$.
  In view of Lemma \ref{lem610},
  this asserts that $h(t)\to 0$ and hence $n(\cdot,t)\to \onz$ in $L^2(\Omega)$ as $t\to\infty$. 
  Along with (\ref{600.2}), this proves (\ref{600.1}).
\qed
Next, combining Corollary \ref{cor588} with Lemma \ref{lem65} yields decay of $u$ in $L^\infty(\Omega)$.
\begin{lem}\label{lem666}
  If $(n,c,u)$ is an eventual energy solution of (\ref{0}), then
  \be{666.1}
	\|u(\cdot,t)\|_{L^\infty(\Omega)} \to 0
	\qquad \mbox{as } t\to\infty.
  \ee
\end{lem}
\proof
  We first observe that as a consequence of Lemma \ref{lem4333}, 
  \be{666.2}
	\int_t^{t+1} \io |u|^2 \to 0
	\qquad \mbox{as } t\to\infty,
  \ee
  meaning that if with some conveniently large $T_1>0$ we let $h(t):=\|u(\cdot,t)\|_{L^2(\Omega)}$, $t>T_1$, then
  $\int_t^{t+1} h(s)ds \to 0$ as $t\to\infty$.
  We now follow the reasoning from Lemma \ref{lem600} and invoke Corollary \ref{cor588} in choosing $T_2>T_1$ such that
  \be{666.3}
	(u(\cdot,t))_{t>T_2}
	\mbox{ is relatively compact in $C^0(\bar\Omega)$ and $h$ is uniformly continuous for $t>T_2$,}
  \ee
  where the latter along with (\ref{666.2}) entails that $h(t)\to 0$ as $t\to\infty$ according to Lemma \ref{lem610}.
  Thus knowing that $u(\cdot,t)\to 0$ in $L^2(\Omega)$ as $t\to\infty$, using the compactness property stated
  in (\ref{666.3}) we readily end up with (\ref{666.1}).
\qed
Proving our final result on large time behavior of arbitrary eventual energy solutions now reduces
to collecting the above convergence and regularity properties. \abs
\proofc of Theorem \ref{theo_eventual} \quad
  The regularity properties in (\ref{reg}) are immediate consequences of Lemma \ref{lem556} and Lemma \ref{lem589}.
  In order to verify (\ref{conv}), we note that
  on choosing $\tu:=u$ and $F(s):=s$ for $s\ge 0$ we obtain from Definition \ref{defi_ees}, (\ref{mass_ees}) 
  and (\ref{cinfty})
  that $(n,c,F)\in \set$ with $m:=\io n_0$, $M:=\|c_0\|_{L^\infty(\Omega)}$ and some $L>0$ and $\tz>0$.
  Therefore, Lemma \ref{lem46} in particular implies that $c(\cdot,t)\to 0$ in $L^\infty(\Omega)$ as $t\to\infty$.
  Combined with the results of Lemma \ref{lem600} and Lemma \ref{lem666}, this proves the claimed stabilization 
  properties.
\qed
\mysection{Existence of an eventual energy solution. Proof of Theorem \ref{theo_exist_ees}}\label{sect9}
Following \cite{win_ct_nasto_exist},
we regularize the original problem (\ref{0}) by fixing families
of approximate initial data $n_{0\eps}, c_{0\eps}$ and $u_{0\eps}$, $\eps\in (0,1)$, with the properties that
\be{I1}
	\left\{ \begin{array}{l}
	n_{0\eps} \in C_0^\infty(\Omega),
	\quad n_{0\eps}\ge 0 \mbox{ in $\Omega$ and}
	\quad
	\io n_{0\eps}=\io n_0
	\quad \mbox{for all } \eps\in (0,1) 
	\qquad \mbox{and} \\[1mm]
	n_{0\eps} \to n_0
	\quad \mbox{in } L\log L(\Omega)
	\qquad \mbox{as } \eps\searrow 0,
	\end{array} \right.
\ee
that
\be{I2}
	\left\{ \begin{array}{l}
	c_{0\eps} \ge 0 \mbox{ in $\Omega$ is such that }
	\sqrt{c_{0\eps}} \in C_0^\infty(\Omega)
	\quad \mbox{and} \quad
	\|c_{0\eps}\|_{L^\infty(\Omega)} \le \|c_0\|_{L^\infty(\Omega)}
	\quad \mbox{for all } \eps\in (0,1) 
	\qquad \mbox{and} \\[1mm]
	\sqrt{c_{0\eps}} \to \sqrt{c_0}
	\quad \mbox{a.e.~in $\Omega$ and in } W^{1,2}(\Omega)
	\qquad \mbox{as } \eps\searrow 0,
	\end{array} \right.
\ee
and that
\be{I3}
	\left\{ \begin{array}{l}
	u_{0\eps} \in C_0^\infty(\Omega) \cap L^2_\sigma(\Omega)
	\quad \mbox{with} \quad
	\|u_{0 \eps}\|_{L^2(\Omega)}=\|u_0\|_{L^2(\Omega)}
	\quad \mbox{for all } \eps\in (0,1) 
	\qquad \mbox{and} \\[1mm]
	u_{0\eps} \to u_0
	\quad \mbox{in } L^2(\Omega)
	\qquad \mbox{as } \eps\searrow 0.
	\end{array} \right.
\ee
Then for $\eps\in (0,1)$, we consider
\be{0eps}
    \left\{ \begin{array}{rcll}
    n_{\eps t} + \ueps\cdot\nabla \neps &=& \Delta \neps - \nabla \cdot (\neps F'_\eps(\neps)\chi(\ceps)\nabla \ceps),
	\qquad & x\in\Omega, \ t>0,\\[1mm]
    c_{\eps t} + \ueps\cdot\nabla \ceps &=& \Delta \ceps-F_\eps(\neps)f(\ceps), \qquad & x\in\Omega, \ t>0,   \\[1mm]
    u_{\eps t}   &=& \Delta \ueps - \nabla \Peps + (\yeps \ueps\cdot\nabla) \ueps + \neps \nabla \phi, 
	\qquad & x\in\Omega, \ t>0, \\[1mm]
    \nabla \cdot \ueps &=& 0, \qquad & x\in\Omega, \ t>0, \\[1mm]
	& & \hspace*{-32mm}
	\frac{\partial\neps}{\partial\nu}=\frac{\partial\ceps}{\partial\nu}=0, \quad \ueps=0,
	\qquad & x\in\pO, \ t>0, \\[1mm]
	& & \hspace*{-32mm}
	\neps(x,0)=n_{0\eps}(x), \quad \ceps(x,0)=c_{0 \eps}(x), \quad \ueps(x,0)=u_{0\eps}(x),
	\qquad & x\in\Omega,
    \end{array} \right.
\ee
where
\be{Feps}
	F_\eps(s):=\frac{1}{\eps} \ln (1+\eps s)
	\qquad \mbox{for } s\ge 0,
\ee
and where
\be{yeps}
	\yeps v:=(1+\eps A)^{-1} v
	\qquad \mbox{for } v\in L^2_\sigma(\Omega).
\ee
The following lemma summarizes some of the results for (\ref{0eps}) obtained 
in \cite[Lemma 2.2, Lemma 2.3, Lemma 3.9, Lemma 3.6]{win_ct_nasto_exist}.
\begin{lem}\label{lem90}
  For each $\eps\in (0,1)$, there exist uniquely determined functions
  \be{90.1}
	\neps \in C^{2,1}(\bar\Omega\times [0,\infty)),
	\quad
	\ceps \in C^{2,1}(\bar\Omega\times [0,\infty))
	\quad \mbox{and} \quad
	\ueps \in C^{2,1}(\bar\Omega\times [0,\infty);\R^3)
  \ee
  which are such that $\neps>0$ and $\ceps>0$ in $\bar\Omega\times (0,\tme)$,
  and such that with some $P_\eps\in C^{1,0}(\Omega\times (0,\tme))$, the quadruple $(\neps,\ceps,\ueps,P_\eps)$
  solves (\ref{0eps}) classically in $\Omega\times (0,\tme)$.\\
  These solutions satisfy 
  \be{mass_eps}
	\io \neps(\cdot,t)=\io n_0
	\qquad \mbox{for all } t>0
  \ee
  as well as
  \be{cinfty_eps}
	\|c_\eps(\cdot,t)\|_{L^\infty(\Omega)} \le \|c_0\|_{L^\infty(\Omega)}
	\qquad \mbox{for all } t>0,
  \ee
  and there exist $\kappa>0$ an $K>0$ such that
  \be{90.11}
	\frac{d}{dt} \F[\neps,\ceps,\ueps](t)
	+ \frac{1}{K} \bigg\{ \io \frac{|\nabla\neps|^2}{\neps} 
	+ \io \frac{|\nabla \ceps|^4}{\ceps^3} + \io |\nabla\ueps|^2 \bigg\}
	\le \ K
	\qquad \mbox{for all } t>0.
  \ee
\end{lem}
Furthermore, Theorem 1.1 in \cite{win_ct_nasto_exist} asserts that
these solutions approach a global weak solution of (\ref{0}) in the following sense.
\begin{lem}\label{lem_limit}
  There exist $(\eps_j)_{j\in\N} \subset (0,1)$ and a global weak solution 
  $(n,c,u)$ of (\ref{0}) such that $\eps_j\searrow 0$ as $j\to\infty$ and 
  $(\neps,\ceps,\ueps) \to (n,c,u)$ a.e.~in $\Omega\times (0,\infty)$ as $\eps=\eps_j\searrow 0$.
  For this solution, we moreover have
  \bea{reg_w}
	& & n\in L^\infty((0,\infty);L^1(\Omega))
	\quad \mbox{with} \quad
	n^\frac{1}{2} \in L^2_{loc}([0,\infty);W^{1,2}(\Omega)), \nn\\
	& & c\in L^\infty(\Omega\times (0,\infty))
	\quad \mbox{with} \quad
	c^\frac{1}{4} \in L^4_{loc}([0,\infty);W^{1,4}(\Omega)), 
	\qquad \mbox{and} \nn\\
	& & u\in L^\infty_{loc}([0,\infty);L^2_\sigma(\Omega))
	\cap L^2_{loc}([0,\infty);W_0^{1,2}(\Omega)),
  \eea
  and there exist $\kappa>0$ and $K>0$ such that
  (\ref{energy1}) and (\ref{energy}) hold with $T:=0$.
\end{lem}
It remains to verify that the component $n$ of this limit function has the additional regularity properties
$\nabla n\in L^2_{loc}(\bar\Omega\times (T,\infty))$ and
$n\in L^4_{loc}(\bar\Omega\times (T,\infty))$ 
required in Definition \ref{defi_ees}.
Thanks to all our previous analysis,
without substantial further efforts
these will result from the fact that for each of the approximate solutions $(\neps,\ceps,\ueps)$,
the triple $(\neps,\ceps,F_\eps)$ lies in ${\cal S}_{m,M,L,0}$ with suitable $m>0, M>0$ and $L>0$:
\begin{lem}\label{lem66}
  There exist $T>0$ and $C>0$ such that for all $\eps\in (0,1)$, the solution of (\ref{0eps}) satisfies
  \be{66.1}
	\int_T^\infty \io |\nabla \neps|^2 \le C
  \ee
  and	
  \be{66.02}
	\int_T^{T+\tilde T} \io \neps^4 \le C\cdot (\tilde T+1)
	\qquad \mbox{for all } \tilde T>0.
  \ee
\end{lem}
\proof
  In order to prepare an application of Lemma \ref{lem64}, we recall that $\io \neps(\cdot,t)=m:=\io n_0$
  and $\|\ceps(\cdot,t)\|_{L^\infty(\Omega)} \le M:=\|c_0\|_{L^\infty(\Omega)}$ for all $t>0$ and $\eps\in (0,1)$
  according to (\ref{mass_eps}) and (\ref{cinfty_eps}).
  Next, Lemma \ref{lem90} combined with Lemma \ref{lem333} says that there exist positive constants $C_1$ and $C_2$,
  independent of $\eps\in (0,1)$, such that with $\kappa>0$ as given there, the function $y$ defined by 
  $y(t):=\F[\neps,\ceps,\ueps](t)$, $t\ge 0$, satisfies
  $y'(t) + C_1 y(t) \le C_2$ for all $t\ge 0$, which shows that $|y(t)| \le C_3$ for all $t\ge 0$ with some $C_3>0$
  independent of $\eps$.
  Therefore, integrating (\ref{90.11}) in time yields
  \be{66.11}
	\int_t^{t+1} \io \bigg\{ \frac{|\nabla\neps|^2}{\neps} + |\nabla\ceps|^4 \bigg\} 
	\le K \max\{1,M^3 \} \cdot (K+2C_3)
	\qquad \mbox{for all $t>0$ and } \eps\in (0,1),
  \ee
  from which we infer that $(\neps,\ceps,F_\eps)\in {\cal S}_{m,M,L,0}$ for all $\eps\in (0,1)$
  if we let $L:=K \max\{1,M^3 \} \cdot (K+2C_3)$.
  Therefore, Lemma \ref{lem46} applies so as to assert the doubly uniform decay property
  \be{66.2}
	\sup_{\eps\in (0,1)} \|\ceps(\cdot,t)\|_{L^\infty(\Omega\times (\tz,\infty))} \to 0	
	\qquad \mbox{as } \tz\to\infty.
  \ee
  In particular, if we let $\eta>0$ and $\tau>0$ denote the numbers obtained from Lemma \ref{lem64} upon the specific choice
  $p:=2$, we can fix $\tz>0$ such that
  \bas
	\|\ceps(\cdot,t)\|_{L^\infty(\Omega\times (\tz,\infty))} \le \eta
	\qquad \mbox{for all } \eps\in (0,1).
  \eas
  Combining this with (\ref{66.11}) and using the outcome of Lemma \ref{lem64}, we thus infer that 
  \be{66.01}
	\|\neps(\cdot,t)\|_{L^3(\Omega)} \le C_1
	\quad \mbox{for all } t>T
	\qquad \mbox{and} \qquad
	\int_T^\infty \io |\nabla \neps|^2 \le C_1
  \ee
  with $T:=\tz+\tau$ and some $C_1>0$ possibly depending on $m$ and $L$ but not on $\eps\in (0,1)$.
  Since using the Gagliardo-Nirenberg inequality we obtain $C_2>0$ such that for each $\tilde T>0$ we have
  \bas
	\int_T^{T+\tilde T} \|\neps(\cdot,t)\|_{L^4(\Omega)}^4 dt
	\le \int_T^{T+\tilde T} \Big\{ C_2 \|\nabla \neps(\cdot,t)\|_{L^2(\Omega)}^2 \|\neps(\cdot,t)\|_{L^3(\Omega)}^2
	+ C_2 \|\neps(\cdot,t)\|_{L^3(\Omega)}^4 \Big\} dt,
  \eas
  both (\ref{66.1}) and (\ref{66.02}) result from (\ref{66.01}).
\qed
We thereby immediately obtain our final result on global existence of an eventual energy solution.\abs
\proofc of Theorem \ref{theo_exist_ees}.\quad
  We let $(n,c,u)$ denote the limit function gained in Lemma \ref{lem_limit}. Then as a consequence of 
  Lemma \ref{lem66} we know that there exist $T>0$ and $C_1>0$ such that
  \bas
	\int_T^\infty \io |\nabla n|^2 < \infty
  \eas
  and
  \bas
	\int_T^{T+\tilde T} \io n^4 \le C_1\cdot (\tilde T+1)
	\qquad \mbox{for all } \tilde T>0.
  \eas
  In conjunction with (\ref{reg_w}), this shows that $(n,c,u)$ enjoys all the regularity properties required
  in Definition \ref{defi_ees} and therefore implies that
  $(n,c,u)$ indeed is an eventual energy solution of (\ref{0}).
\qed
\mysection{Appendix}
For convenience, let us include a short proof of a stability property of Steklov averages which is essential
to our arguments in both Lemma \ref{lem50} and also Lemma \ref{lem31}.
\begin{lem}\label{lem60}
  Let $-\infty<T_0<T_1<\infty$, and for $h\in (0,1)$ define
  \be{S}
	S_h[\varphi](x,t):=\frac{1}{h} \int_t^{t+h} \varphi(x,s) ds,
	\qquad x\in\Omega, \ t\in (T_0,T_1),
  \ee
  for $\varphi\in L^1(\Omega\times (T_0,T_1))$.\abs
  i) \ If $\varphi\in L^p(\Omega\times (T_0,T_1+1))$ for some $p\in [1,\infty)$, then
  \be{60.1}
	\|S_h[\varphi]\|_{L^p(\Omega\times (T_0,T_1))} \le \|\varphi\|_{L^p(\Omega\times (T_0,T_1+h))}
	\qquad \mbox{for all } h\in (0,1).
  \ee
  ii) \ Suppose that for some $p\in [1,\infty)$, $(\varphi_h)_{h\in (0,1)} \subset L^p(\Omega\times (T_0,T_1+1))$
  and $\varphi \in L^p(\Omega\times (T_0,T_1+1))$ are such that
  \be{60.2}
	\|\varphi_h-\varphi\|_{L^p(\Omega\times (T_0,T_1+h))} \to 0
	\qquad \mbox{as } h\searrow 0.
  \ee
  Then 
  \be{60.3}
	S_h[\varphi_h] \to \varphi
	\quad \mbox{in } L^p(\Omega\times (T_0,T_1))
	\qquad \mbox{as } h\searrow 0.
  \ee
\end{lem}
\proof
  We may assume that $T_0=0$ and write $T:=T_1$.\\
  i) \ Using the H\"older inequality and the Fubini theorem, we directly see that
  \bas
	& & \hspace*{-16mm}
	\int_0^T \io \Big|S_h[\varphi](x,t)\Big|^p dxdt 
	= \frac{1}{h^p} \io \int_0^T \bigg| \int_t^{t+h} \varphi(x,s)ds \bigg|^p dtdx \\
	&\le& \frac{1}{h} \io \int_0^T \int_t^{t+h} |\varphi(x,s)|^p dsdtdx \\
	&=& \frac{1}{h} \io \bigg\{ \int_h^T \int_{s-h}^s |\varphi(x,s)|^p dtds
	+ \int_0^h \int_0^s |\varphi(x,s)|^p dtds
	+ \int_T^{T+h} \int_{s-h}^T |\varphi(x,s)|^p dtds \bigg\} dx \\
	&=& \frac{1}{h} \io \bigg\{ h\int_h^T |\varphi(x,s)|^p ds
	+ \int_0^h s |\varphi(x,s)|^p ds
	+ \int_T^{T+h} (T+h-s) |\varphi(x,s)|^p ds \bigg\} dx \\
	&\le& \frac{1}{h} \io \bigg\{ h\int_h^T |\varphi(x,s)|^p ds 
	+ h\int_0^h |\varphi(x,s)|^p ds
	+ h\int_T^{T+h} |\varphi(x,s)|^p ds \bigg\} dx \\
	&=& \int_0^{T+h} \io |\varphi(x,s)|^p dxds
  \eas
  for all $h\in (0,1)$.\\
  ii) \ By linearity of $S_h$ and i), we have
  \bas
	\big\| S_h[\varphi_h] - \varphi\big\|_{L^p(\Omega\times (0,T))}
	&=& \big\| S_h[\varphi_h-\varphi] + (S_h[\varphi]-\varphi) \big\|_{L^p(\Omega\times (0,T))} \\
	&\le& \big\|S_h[\varphi_h-\varphi]\big\|_{L^p(\Omega\times (0,T))}
	+ \big\|S_h[\varphi]-\varphi\big\|_{L^p(\Omega\times (0,T))} \\
	&\le& \|\varphi_h-\varphi\|_{L^p(\Omega\times (0,T+h))}
	+ \big\| S_h[\varphi]-\varphi\big\|_{L^p(\Omega\times (0,T))}
	\qquad \mbox{for all } h\in (0,1).
  \eas
  Since $S_h[\varphi]\to \varphi$ in $L^p(\Omega\times (0,T))$ by a well-known result 
  (see e.g.~\cite[Lemma I.3.2]{DiBenedetto}), (\ref{60.2}) therefore entails (\ref{60.3}).
\qed
We also separately state the following interpolation lemma which is used in several places.
\begin{lem}\label{lem511}
  There exists $C>0$ such that whenever $J\subset \R$ is an interval, the following holds.\\
  i) \ Any function
  $n\in L^\infty(J;L^1(\Omega)) \cap L^2(J;W_0^{1,2}(\Omega))$ belongs to
  $L^\frac{8}{3}(\Omega\times J)$ and satisfies
  \be{511.2}
	\int_J \io n^\frac{8}{3}
	\le C \|\nabla n\|_{L^2(\Omega\times J)}^2 \|n\|_{L^\infty(J;L^1(\Omega))}^\frac{2}{3}.
  \ee
  ii) \ If
  $u\in L^\infty(J;L^2(\Omega;\R^3)) \cap L^2(J;W_0^{1,2}(\Omega;\R^3))$, then
  $u\in L^\frac{10}{3}(\Omega\times J;\R^3)$ with
  \be{511.1}
	\int_J \io |u|^\frac{10}{3}
	\le C \|\nabla u\|_{L^2(\Omega\times J)}^2 \|u\|_{L^\infty(J;L^2(\Omega))}^\frac{4}{3}.
  \ee
\end{lem}
\proof
  Both statements can be obtained upon straightforward interpolation using the Gagliardo-Nirenberg inequality.
\qed
{\bf Acknowledgement.} \quad
The author would like to thank Johannes Lankeit and Ken Abe 
for numerous useful comments which lead to significant improvements in the manuscript.
\end{document}